\documentclass{amsart}
\usepackage{amssymb,amsmath,amsthm,latexsym,booktabs,todonotes, color, comment, eucal,uncial,eufrak}
\usepackage[normalem]{ulem}
\newcommand{\stkout}[1]{\ifmmode\text{\sout{\ensuremath{#1}}}\else\sout{#1}\fi}
\theoremstyle{definition}
\newtheorem{definition}{Definition}
\newtheorem{notation}[definition]{Notation}

\newtheorem{remark}[definition]{Remark}
\newtheorem{example}[definition]{Example}

\theoremstyle{plain}
\newtheorem{lemma}[definition]{Lemma}
\newtheorem{proposition}[definition]{Proposition}
\newtheorem{theorem}[definition]{Theorem}
\newtheorem{corollary}[definition]{Corollary}

\begin{document}
\title{Equatorially balanced $C_4$-face-magic Labelings 
on Klein Bottle Grid Graphs}
\author{Stephen J. Curran}
\address{Department of Mathematics\\
         University of Pittsburgh at Johnstown\\
         450 Schoolhouse Rd\\
         Johnstown, PA 15904\\
         USA}
         \email{sjcurran@pitt.edu}
         \author{Richard M. Low}
\address{Department of Mathematics\\
         San Jose State University\\
         1 Washington Sq\\
         San Jose, CA 95192\\
         USA}
\email{richard.low@sjsu.edu}
\author{Stephen C. Locke}
\address{Department of Mathematical Sciences\\
         Florida Atlantic University\\
         777 Glades Rd\\
         Boca Raton, FL 33431\\
         USA}
\email{lockes@fau.edu}

\date{June 3, 2022.  }

\begin{abstract}
For a graph $G = (V, E)$ embedded in the Klein bottle, let $\mathcal{F}(G)$ denote the set of faces of $G$.
Then, $G$ is called a $C_k$-face-magic Klein bottle  graph if there exists a bijection $f: V(G) \to \{1, 2, \dots, |V(G)|\}$
such that for any $F \in \mathcal{F}(G)$ with $F \cong C_k$, the sum of all the vertex labelings along $C_k$ is a constant $S$.
Let $x_v =f(v)$ for all $v\in V(G)$. We call $\{x_v : v\in V(G)\}$ a $C_k$-face-magic Klein bottle labeling on $G$.
We consider the $m \times n$ grid graph, denoted by $\mathcal{K}_{m,n}$, embedded in the Klein bottle in the natural way.
We show that for $m,n\ge 2$, $\mathcal{K}_{m,n}$ admits a $C_4$-face-magic Klein bottle labeling if and only if
$n$ is even. We say that a $C_4$-face-magic Klein bottle labeling $\{x_{i,j}: (i,j) \in V(\mathcal{K}_{m,n}) \}$ on $\mathcal{K}_{m,n}$
is equatorially balanced if $x_{i,j} + x_{i,n+1-j} = \tfrac{1}{2} S$ for all $(i,j) \in V(\mathcal{K}_{m,n})$.
We show that when $m$ is odd, a $C_4$-face-magic Klein bottle labeling on $\mathcal{K}_{m,n}$
must be equatorially balanced. Also when $m$ is odd, we show that (up to symmetries on the Klein bottle)
the number of $C_4$-face-magic Klein bottle labelings on the $m \times 4$ Klein bottle grid graph
is $2^m \, (m-1)! \, \tau(m)$, where $\tau(m)$ is the number of positive divisors of $m$.

Furthermore, let $m\ge 3$ be an odd integer and $n \ge 6$ be an even integer.
Then, the minimum number of distinct $C_4$-face-magic Klein bottle
labelings $X$ on $\mathcal{K}_{m,n}$ (up to symmetries on a Klein bottle) is either $(5\cdot 2^m)(m-1)!$ if $n \equiv 0\pmod{4}$,
or $(6\cdot 2^m)(m-1)!$ if $n \equiv 2\pmod{4}$.
\end{abstract}

\maketitle

\section{Introduction} \label{S:introduction}

Graph labelings were formally introduced in the 1970s by Kotzig and Rosa \cite{Kotzig}.
Graph labelings have been applied to graph decomposition problems, radar pulse code designs, X-ray
crystallography and communication network models.
The interested reader should consult J.A. Gallian's comprehensive dynamic survey on graph labelings \cite{Gallian}
for further investigation.

We refer the reader to Bondy and Murty \cite{Bondy} for concepts and notation not explicitly defined in this paper.
All graphs in this paper are simple and connected.
The concept of a $C_k$-face-magic labeling was first applied to planar graphs.
For a planar (toroidal, cylindrical, or Klein bottle) graph $G = (V, E)$ embedded in the plane (torus, cylinder, or Klein bottle),
let $\mathcal{F}(G)$ denote the set of faces of $G$.
Then, $G$ is called a \textit{$C_k$-face-magic planar (toroidal, cylindrical, or Klein bottle)} graph if there
exists a bijection $f: V(G) \to \{1, 2, \dots, |V(G)|\}$ such that for any $F \in \mathcal{F}(G)$ with $F \cong C_k$,
the sum of all the vertex labelings along $C_k$ is a constant $S$. Here, the constant $S$ is called
a \textit{$C_k$-face-magic value} of $G$.
More generally, $C_k$-face-magic planar graph labelings are a special case of $(a, b, c)$-magic labelings
introduced by Lih \cite{Lih}.
For assorted values of $a, b$ and $c$, Baca and others \cite{Baca1, Baca2, Baca3, Hsieh, Kasif, Kath, Lih} have analyzed
the problem for various classes of graphs.
Wang \cite{Wang} showed that the toroidal grid graphs $C_m \times C_n$ are antimagic for all integers $m,n\ge 3$.
Butt et al. \cite{Buttetal} investigated face antimagic labelings on toroidal and Klein bottle grid graphs.
Curran, Low and Locke \cite{CurranLow, CurranLowLocke} investigated $C_4$-face-magic toroidal labelings
on $C_m \times C_n$ and $C_4$-face-magic cylindrical labelings on $P_m \times C_n$.
Curran and Low \cite{CurranLow} showed that (up to symmetries on a torus) there are three
$C_4$-face-magic toroidal labelings on $C_4 \times C_4$ that are row-sum balanced, column-sum balanced, and diagonal-sum balanced.
Curran, Low and Locke \cite{CurranLowLocke} showed that for all positive integers $n$,
$C_{4n} \times C_{4n}$ has a $C_4$-face-magic toroidal labeling that is row-sum balanced, column-sum balanced, and diagonal-sum balanced.
In this paper, we investigate $C_4$-face-magic Klein bottle labelings on an $m \times n$
Klein bottle grid graph.
We show that for $m,n\ge 2$, $\mathcal{K}_{m,n}$ has a $C_4$-face-magic Klein bottle labeling if and only if
$n$ is even.
When $m$ is odd, we show that (up to symmetries on the Klein bottle)
the number of $C_4$-face-magic Klein bottle labelings on the $m \times 4$ Klein bottle grid graph
is $2^m \, (m-1)! \, \tau(m)$, where $\tau(m)$ is the number of positive divisors of $m$.
Furthermore, let $m\ge 3$ be an odd integer and $n \ge 6$ be an even integer.
Then, the minimum number of distinct $C_4$-face-magic Klein bottle
labelings $X$ on $\mathcal{K}_{m,n}$ (up to symmetries on a Klein bottle) is either $(5\cdot 2^m)(m-1)!$ if $n \equiv 0\pmod{4}$,
or $(6\cdot 2^m)(m-1)!$ if $n \equiv 2\pmod{4}$.

\section{Preliminaries} \label{S:preliminaries}

\begin{notation}
Let $m$ be a positive integer. For convenience,
we let $\stkout{m} =\lfloor m/2 \rfloor$, $\stkout{m}^{+} =\lfloor (m+1)/2 \rfloor$, and $\stkout{m}^{-} =\lfloor (m-1)/2 \rfloor$.
\end{notation}

\begin{definition}
For a graph $G = (V, E)$ embedded on the Klein bottle (plane or torus or cylinder),
let $\mathcal{F}(G)$ denote the set of faces of $G$.
Then, $G$ is called a \textit{$C_n$-face-magic Klein bottle (planar or toroidal or cylindrical)}
graph if there exists a bijection $f: V(G) \to \{1, 2, \dots, |V(G)|\}$ such that
for any $F \in \mathcal{F}(G)$ with $F \cong C_n$, the sum of all the vertex labelings
along $C_n$ is a constant $S$. We call $S$ the \textit{$C_4$-face-magic value}.
\end{definition}

\begin{theorem}[\cite{CurranLowLocke}, Theorem 1] \label{thmPmTimesPmIsFaceMagic}
Let $m, n \geq 2$. Then, $P_m \times P_n$ is $C_4$-face-magic planar.
\end{theorem}

\begin{lemma}[\cite{CurranLowLocke}, Lemma 2]  \label{lemCylinderLabelNIsEven}
Let $m$ and $n$ be integers such that $m \geq 3$ and $n \geq 2$. Suppose $P_m \times C_n$ is a
$C_4$-face-magic cylindrical graph with the natural embedding of $P_m \times C_n$ on the cylinder. Then, $n$ is even.
\end{lemma}

\begin{definition}
Let $m$ and $n$ be integers such that $m,n\geq 2$.
The {\it $m \times n$ Klein bottle grid graph}, denoted by $\mathcal{K}_{m,n}$, is the graph whose vertex set is
\begin{align*}
V\left( \mathcal{K}_{m,n} \right)
=\left\{ \left( i,j\right) :1\leq i\leq m,1\leq j\leq n\right\},
\end{align*}
and whose edge set consists of the following edges:
\begin{itemize}
\item There is an edge from $(i,j)$ to $(i,j+1)$, for $1\leq i \leq m$ and $1\leq j \leq n-1$.
\item There is an edge from $(i,n)$ to $(i,1)$, for $1\leq i \leq m$.
\item There is an edge from $(i,j)$ to $(i+1,j)$, for $1\leq i \leq m-1$ and $1\leq j \leq n$.
\item There is an edge from $(m,j)$ to $(1,n+1-j)$, for $1\leq j \leq n$.
\end{itemize}
The graph $\mathcal{K}_{m,n}$ has a natural embedding on the Klein bottle.
\end{definition}

\begin{proposition} \label{propKleinCheckerboardNIsEven}
Let $m$ and $n$ be integers such that $m,n\geq 2$. Suppose that $\mathcal{K}_{m,n}$
is a $C_4$-face-magic Klein bottle graph. Then, $n$ is even.
\end{proposition}

\begin{proof}
Let $\{ x_{i,j} : (i,j) \in V( \mathcal{K}_{m,n}) \}$ be a $C_4$-face-magic Klein bottle labeling on $\mathcal{K}_{m,n}$
Suppose $m\geq 3$. This labeling yields a $C_4$-face-magic cylindrical labeling on $P_m \times C_n$.
By Lemma \ref{lemCylinderLabelNIsEven}, $n$ is even.

For the purposes of contradiction, assume that $m=2$ and $n$ is odd.
Let $n=2 n_1 +1$ for some positive integer $n_1$.
When we set the two $C_4$-face sums given below equal to each other
\begin{equation*}
x_{1,n_1} + x_{1,n_1 +1} + x_{2,n_1} + x_{2,n_1 +1} = S =  x_{2,n_1} + x_{2,n_1 +1} + x_{1,n_1 +1} + x_{1,n_1 +2},
\end{equation*}
we obtain $x_{1,n_1} = x_{1,n_1 +2}$. This is a contradiction. Hence, $n$ is even. This completes the proof.
\end{proof}

\begin{lemma} \label{lemmaKleinCheckerboardC4FaceValue}
Let $m$ be an integer such that $m\geq 2$ and let $n$ be an even positive integer.
Let $\{ x_{i,j} : (i,j)\in V(\mathcal{K}_{m,n}) \}$ be a $C_4$-face-magic Klein bottle labeling on $\mathcal{K}_{m,n}$.
Then, the $C_4$-face-magic value of this labeling is $S= 2(mn+1)$.
\end{lemma}

\begin{proof}
Consider the sum
\begin{align*}
mnS = &\sum_{i=1}^{m-1} \sum_{j=1}^{n} ( x_{i,j} + x_{i+1,j} + x_{i,j+1} + x_{i+1,j+1}) \\
 &+ \; \sum_{i=1}^{m} (  x_{m,j} + x_{m,j+1} + x_{1,n-j} + x_{1,n-j+1}) \; \\
 = \; &4\biggl( \sum_{k=1}^{mn} k \biggr) = (2mn)(mn+1),
\end{align*}
where the indices $i$ and $j$ are taken modulo $m$ and $n$, respectively.
Thus, $S=2(mn+1)$.
\end{proof}

\section{$C_4$-face-magic labelings on an $m \times n$ Klein bottle grid graph} \label{S:KleinbottleCheckerboard}

\begin{definition}
Let $m,n\ge 2$ be integers such that $n$ is even.
Let $X = \{ x_{i,j} : (i,j)\in V(\mathcal{K}_{m,n})\}$ be a $C_4$-face-magic Klein bottle labeling on $\mathcal{K}_{m,n}$.
Let $S= 2(mn+1)$ be the $C_4$-face-magic value of this labeling.
We say that $X$ is \textit{equatorially balanced} if, for all $1\le i \le m$ and $1\le j \le \stkout{n}$, we have
\begin{equation*}
x_{i,j} + x_{i,n+1-j} = \tfrac{1}{2} S = mn+1.
\end{equation*}
\end{definition}

\begin{lemma} \label{lemmaKleinCheckerboardSemiC4FaceValueBalance}
Let $m$ be an odd integer such that $m\geq 3$ and let $n$ be an even positive integer.
Let $X = \{ x_{i,j} : (i,j)\in V(\mathcal{K}_{m,n}) \}$ be a $C_4$-face-magic Klein bottle labeling on $\mathcal{K}_{m,n}$.
Then, $X$ is equatorially balanced.
\end{lemma}

\begin{proof}
Consider the $C_4$-face sums
\begin{align*}
S &= x_{i,1} + x_{i+1,1} + x_{i,n} + x_{i+1,n} \mbox{ and } \\
S &= x_{i+1,1} + x_{i+2,1} + x_{i+1,n} + x_{i+2,n},
\mbox{ where } i=1,2,\ldots,m.
\end{align*}
When we equate these sums and simplify, we obtain
\begin{equation*}
x_{i,1} + x_{i,n} = x_{i+2,1} + x_{i+2,n}, \mbox{ where } i=1,2,\ldots,m.
\end{equation*}
Since the indices $i$ and $i+2$ are taken modulo $m$, and $m$ is odd, we have
$x_{i_1,1} + x_{i_1,n} = x_{i_2,1} + x_{i_2,n}$, for all $1\leq i_1, i_2 \leq m$.
Thus, $S=  x_{i,1} + x_{i,n} + x_{i+1,1}  + x_{i+1,n} =2( x_{i,1} + x_{i,n})$, which implies $ x_{i,1} + x_{i,n} =\tfrac{1}{2} S$, for all $i=1,2,\ldots,m$.

Assume that for some positive integer $k$ such that $k< \stkout{n}$, we have
for all $1\leq i\leq m$ and $1\leq j\leq k$ that
$ x_{i,j} + x_{i,n+1-j} =\tfrac{1}{2} S$.
When we equate the $C_4$-face sums given below together
\begin{equation*}
 x_{i,k} + x_{i,k+1} + x_{i+1,k} + x_{i+1,k+1} = S = x_{i+1,k} + x_{i+1,k+1} + x_{i+2,k} + x_{i+2,k+1},
\end{equation*}
we obtain
\begin{equation*}
x_{i,k} + x_{i,k+1} = x_{i+2,k} + x_{i+2,k+1}, \mbox{ for all } 1\leq i \leq m-2.
\end{equation*}
A similar argument shows that
\begin{equation*}
x_{i,n-k} + x_{i,n-k+1} = x_{i+2,n-k} + x_{i+2,n-k+1}, \mbox{ for all } 1\leq i \leq m-2.
\end{equation*}
When we equate the $C_4$-face sums given below together
\begin{equation*}
 x_{m-1,k} + x_{m-1,k+1} + x_{m,k} + x_{m,k+1} = S = x_{m,k} + x_{m,k+1} + x_{1,n-k} + x_{1,n-k+1},
\end{equation*}
we obtain
\begin{equation*}
x_{m-1,k} + x_{m-1,k+1} = x_{1,n-k} + x_{1,n-k+1}.
\end{equation*}
Similarly, when we equate the $C_4$-face sums given below together
\begin{equation*}
 x_{m,k} + x_{m,k+1} + x_{1,n-k} + x_{1,n-k+1} = S = x_{1,n-k} + x_{1,n-k+1} + x_{2,n-k} + x_{2,n-k+1},
\end{equation*}
we obtain
\begin{equation*}
x_{m,k} + x_{m,k+1} = x_{2,n-k} + x_{2,n-k+1}.
\end{equation*}
Let $S_i =x_{i,k} + x_{i,k+1}$, for $i=1,2$.
Then,
\begin{align*}
S_1 = x_{2i+1,k} + x_{2i+1,k+1}, \mbox{ for } i=0,1,2,\ldots,\stkout{m}, \mbox{ and } \\
S_1 = x_{2i,n-k} + x_{2i,n-k+1}, \mbox{ for } i=1,2,\ldots,\stkout{m}.
\end{align*}
Similarly, we have
\begin{align*}
S_2 = x_{2i,k} + x_{2i,k+1}, \mbox{ for } i=1,2,\ldots,\stkout{m}, \mbox{ and } \\
S_2 = x_{2i+1,n-k} + x_{2i+1,n-k+1}, \mbox{ for } i=0,1,2,\ldots,\stkout{m}.
\end{align*}
Since $i$ and $i+1$ have opposite parity, we have $x_{i+1,k} +x_{i+1,k+1} = x_{i,n-k} + x_{i,n+1-k}$.
From the $C_4$-face sum $x_{i,k} + x_{i,k+1} + x_{i+1,k} + x_{i+1,k+1} =S$, we have
$x_{i,k} + x_{i,k+1} + x_{i,n-k} +x_{i,n+1-k} + S$.
From the inductive hypothesis, we have $x_{i,k} + x_{i,n+1-k} = \tfrac{1}{2} S$.
Thus, it follows that
\begin{equation*}
 x_{i,k+1} + x_{i,n-k}  = S - ( x_{i,k} + x_{i,n+1-k} ) = \tfrac{1}{2} S.
\end{equation*}
This completes the proof.
\end{proof}

\begin{proposition} \label{propKleinCheckerboardMOddNEvenC4FaceExists}
Let $m$ be an odd integer such that $m\geq 3$ and let $n$ be an even positive integer such that $n=2n_1$, for some integer $n_1$.
Define a labeling $X = \{ x_{i,j} : (i,j)\in V( \mathcal{K}_{m,n}) \}$  on $\mathcal{K}_{m,n}$ as follows:
\begin{itemize}
\item $x_{2i-1,2j-1} = (j-1)m+i$, for $i=1,2,\ldots,\stkout{m}^{+}$ and $j=1,2,\ldots,\stkout{n}_1^{+}$.
\item $x_{2i,2j-1} = (n+1-j)m+1-\stkout{m}^{+} -i$, for $i=1,2,\ldots,\stkout{m}$ and $j=1,2,\ldots,\stkout{n}_1^{+}$.
\item $x_{2i-1,2j} = (n_1 -j+1)m+1-i$, for $i=1,2,\ldots,\stkout{m}^{+}$ and $j=1,2,\ldots,\stkout{n}_1$.
\item $x_{2i,2j} = (n_1+j-1)m+\stkout{m}^{+} +i$, for $i=1,2,\ldots,\stkout{m}$ and $j=1,2,\ldots,\stkout{n}_1$.
\end{itemize}
Furthermore, for $1\leq i \leq m$ and $\stkout{n} +1 \leq j\leq n$, let
$x_{i,j} = mn+1-x_{i,n+1-j}$.
Then, $X$ is a $C_4$-face-magic Klein bottle labeling on $\mathcal{K}_{m,n}$.
\end{proposition}

\begin{proof}
Let $m= 2m_1 +1$ for some positive integer $m_1$. Then, $\stkout{m}^{+} = m_1 +1$ and $\stkout{m} = m_1$.
We observe that $\stkout{n}_{1}^{+} \, + \, \stkout{n}_{1} \, =n_1$.
Let $S$ be the $C_4$-face-magic value of $X$.
By Lemma \ref{lemmaKleinCheckerboardC4FaceValue}, $S=2(mn+1)$.
We verify the $C_4$-face sum
\begin{equation} \label{eqnC4FaceSumij}
x_{i,j} + x_{i+1,j} + x_{i+1,j} + x_{i+1,j+1} = S, \mbox{ \ \ for all } 1\leq i \leq m \mbox{ and } 1\leq j \leq n_1.
\end{equation}
We first verify equation (\ref{eqnC4FaceSumij}) for $j=n_1$.
By Lemma \ref{lemmaKleinCheckerboardSemiC4FaceValueBalance}, we have
\begin{equation*}
x_{i,n_1} +x_{i, n_1 +1} = \tfrac{1}{2} S, \mbox{ \ \ \ for all } 1\le i \le m.
\end{equation*}
Thus, for all $1\le i \le m$, we have
\begin{equation*}
x_{i,n_1} +x_{i, n_1 +1} +  x_{i+1,n_1} +x_{i+1, n_1 +1} = \tfrac{1}{2} S + \tfrac{1}{2} S = S.
\end{equation*}

\textbf{Case 1.} Assume $i$ is odd and $j$ is odd in equation (\ref{eqnC4FaceSumij}).
We verify the $C_4$-face sum
$x_{2i-1,2j-1} + x_{2i-1,2j} + x_{2i,2j-1} + x_{2i,2j} =S$, for $1\leq i \leq m_1 +1$ and $1\leq j \leq \stkout{n}_{1}$.
Assume that $1\leq i\leq m_1$ and $1\leq j \leq \stkout{n}_1$. Then, $x_{2i-1,2j-1} = (j-1)m+i$,
$x_{2i,2j-1} = (n+1-j)m +1 -\stkout{m}^{+} -i$,
$x_{2i-1,2j} = (n_1 -j+1)m +1-i$ and $x_{2i,2j} = (n_1 +j-1)m + \stkout{m}^{+} +i$.
Thus,
\begin{align*}
x_{2i-1,2j-1} &+ x_{2i-1,2j} +  x_{2i,2j-1} +  x_{2i,2j} \\
&= \bigl(  (j-1)m+i \bigr) +  \bigl(  (n_1 -j+1)m +1-i \bigr) \\
&+ \bigl(  (n+1-j)m +1 - \stkout{m}^{+} -i \bigr) + \bigl(  (n_1 +j-1)m + \stkout{m}^{+} +i \bigr) \\
&= 2(mn+1) =S.
\end{align*}

Assume $i=m_1 +1$ and $1\leq j \leq  \stkout{n}_1$.
We have $x_{m,2j-1} = (j-1)m + m_1 + 1$,
$x_{m,2j} = (n_1 -j +1)m - m_1$,
$x_{1,2j-1} =  (j-1)m +1 $, and
$x_{1,2j} =  (n_1 -j +1)m$.
Thus, $x_{m,2j-1} + x_{m,2j} = n_1 m +1 =   x_{1,2j-1} + x_{1,2j}$.
Since $x_{1,(n+1)-(2j-1)} = \tfrac{1}{2} S - x_{1,2j-1}$ and
$x_{1,(n+1)-2j} = \tfrac{1}{2}  S - x_{1,2j}$, we have
\begin{align*}
x_{m,2j-1} &+ x_{m,2j} + x_{1,(n+1)-(2j-1)} + x_{1,(n+1)-2j} \\
&= x_{m,2j-1} + x_{m,2j} + \bigl( \tfrac{1}{2} S - x_{1,2j-1} \bigr)  + \bigl( \tfrac{1}{2} S - x_{1,2j} \bigr) \\
&= \bigl(  n_1 m +1 \bigr) + S - \bigl(  n_1 m  +1  \bigr) =S.
\end{align*}
Similarly, since $x_{m,(n+1)-(2j-1)} = \tfrac{1}{2} S - x_{m,2j-1}$ and
$x_{m,(n+1)-2j} = \tfrac{1}{2}  S - x_{m,2j}$, we have
\begin{align*}
x_{1,2j-1} &+ x_{1,2j} + x_{m,(n+1)-(2j-1)} + x_{m,(n+1)-2j} \\
&= x_{1,2j-1} + x_{1,2j} + \bigl( \tfrac{1}{2} S - x_{m,2j-1} \bigr)  + \bigl( \tfrac{1}{2} S - x_{m,2j} \bigr) \\
&= \bigl(  n_1 m +1 \bigr) + S - \bigl(  n_1 m  +1  \bigr) =S.
\end{align*}

\textbf{Case 2.}  Assume $i$ is even and $j$ is odd in equation (\ref{eqnC4FaceSumij}).
We verify the $C_4$-face sum
$x_{2i,2j-1} + x_{2i,2j} + x_{2i+1,2j-1} + x_{2i+1,2j} =S$, for
$1\leq i \leq  m_1$ and $1\leq j \leq  \stkout{n}_1$.
We have $x_{2i,2j-1} = (n+1-j)m +1 - \stkout{m} -i$, $x_{2i,2j} = (n_1 +j-1)m + \stkout{m} +i$,
$x_{2(i+1)-1,2j-1} = (j-1)m+i+1$ and $x_{2(i+1)-1,2j} = (n_1 -j+1)m -i$.
Thus,
\begin{align*}
x_{2i,2j-1} &+ x_{2i,2j} +  x_{2(i+1)-1,2j-1} +  x_{2(i+1)-1,2j} \\
&= \bigl(  (n+1-j)m +1 - \stkout{m} -i \bigr) + \bigl(  (n_1 +j-1)m + \stkout{m} +i  \bigr) \\
&+ \bigl( (j-1)m+i+1 \bigr) + \bigl(  (n_1 -j+1)m -i \bigr) \\
&= 2(mn+1) =S.
\end{align*}

\textbf{Case 3.}  Assume $i$ is odd and $j$ is even in equation (\ref{eqnC4FaceSumij}).
We verify the $C_4$-face sum
$x_{2i-1,2j} + x_{2i-1,2j+1} + x_{2i,2j} + x_{2i,2j+1} =S$, for
$1\leq i \leq  m_1 +1$ and $1\leq j \leq  \stkout{n}_{1}^{-}$.
Assume $1\leq i \leq  m_1$ and
$1\leq j \leq \stkout{n}_{1}^{-}$.
Then, $x_{2i-1,2j} = (n_1 -j +1)m +1  -i$, $x_{2i-1,2j+1} = x_{2i-1,2(j+1)-1} = jm  +i$,
$x_{2i,2j} = (n_1 + j-1)m + \stkout{m}^{+} +i$ and $x_{2i,2j+1} = x_{2i,2(j+1)-1} = (n -j)m +1 - \stkout{m}^{+} - i$.
Thus,
\begin{align*}
x_{2i-1,2j} &+ x_{2i-1,2j+1} +  x_{2i,2j} +  x_{2i,2j+1} \\
&= \bigl( (n_1 -j +1)m +1  -i  \bigr) + \bigl(  jm  +i  \bigr) \\
&+ \bigl( (n_1 + j-1)m + \stkout{m}^{+} +i \bigr) + \bigl(   (n -j)m +1 - \stkout{m}^{+} - i \bigr) \\
&= 2(mn+1) =S.
\end{align*}

Next, assume $i= m_1 +1$ and $1\leq j \leq \stkout{n}_{1}^{-}$.
We have $x_{m,2j} = x_{2(m_1 +1) -1, 2j} = (n_1 -j +1)m - m_1$,
$x_{m,2j+1} = x_{2(m_1 +1)-1, 2(j+1)-1} = jm +m_1 +1$,
$x_{1,2j} = (n_1 -j +1)m$, and
$x_{1,2j+1} = x_{1,2(j+1)-1} = jm+1$.
Then, $x_{m,2j}  + x_{m,2j+1} =  (n_1  +1)m +1  =  x_{1,2j}  + x_{1,2j+1}$.
Since $x_{1,(n+1)-(2j)} = \tfrac{1}{2} S - x_{1,2j}$ and
$x_{1,(n+1)-(2j+1)} = \tfrac{1}{2}  S - x_{1,2j+1}$, we have
\begin{align*}
x_{m,2j} &+ x_{m,2j+1} + x_{1,(n+1)-2j} + x_{1,(n+1)-(2j+1)} \\
&= x_{m,2j} + x_{m,2j+1} + \bigl( \tfrac{1}{2} S - x_{1,2j} \bigr)  + \bigl( \tfrac{1}{2} S - x_{1,2j+1} \bigr) \\
&= \bigl(  (n_1 +1 )m +1 \bigr) + S - \bigl(  (n_1 +1) m  +1  \bigr) =S.
\end{align*}
Similarly, since $x_{m,(n+1)-2j} = \tfrac{1}{2} S - x_{m,2j}$ and
$x_{m,(n+1)-(2j+1)} = \tfrac{1}{2}  S - x_{m,2j+1}$, we have
\begin{align*}
x_{1,2j} &+ x_{1,2j+1} + x_{m,(n+1)-2j} + x_{m,(n+1)-(2j+1)} \\
&= x_{1,2j-1} + x_{1,2j} + \bigl( \tfrac{1}{2} S - x_{m,2j} \bigr)  + \bigl( \tfrac{1}{2} S - x_{m,2j+1} \bigr) \\
&= \bigl(  (n_1 +1) m +1 \bigr) + S - \bigl(  (n_1 +1) m  +1  \bigr) =S.
\end{align*}

\textbf{Case 4.}  Assume $i$ is even and $j$ is even in equation (\ref{eqnC4FaceSumij}).
We verify the $C_4$-face sum
$x_{2i,2j} + x_{2i,2j+1} + x_{2i+1,2j} + x_{2i+1,2j+1} =S$, for
$1\leq i \leq  m_1$ and $1\leq j \leq  \stkout{n}_{1}^{-}$.
We have $x_{2i,2j} = (n_1 +j -1)m + \stkout{m}^{+}  + i$, $x_{2i,2j+1} = x_{2i,2(j+1)-1} = (n-j)m +1 - \stkout{m}^{+} -i$,
$x_{2i+1,2j} = x_{2(i+1)-1,2j} = (n_1 - j +1)m - i$
and $x_{2i+1,2j+1} = x_{2(i+1)-1,2(j+1)-1} = jm +i +1$.
Thus,
\begin{align*}
x_{2i,2j} &+ x_{2i,2j+1} +  x_{2i+1,2j} +  x_{2i+1,2j+1} \\
&= \bigl( (n_1 +j -1)m + \stkout{m}^{+}  + i \bigr) + \bigl(  (n-j)m +1 - \stkout{m}^{+} -i  \bigr) \\
&+ \bigl( (n_1 - j +1)m - i \bigr) + \bigl(   jm +i +1 \bigr) \\
&= 2(mn+1) =S.
\end{align*}
\vspace{0.1in}

Next, for $1 \leq i < m$ and $n_1 +1 \leq j < n$, we have
\begin{align*}
x_{i,j} &+ x_{i+1,j} + x_{i,j+1} + x_{i+1,j+1} \\
&= \bigl( \tfrac{1}{2} S - x_{i, n+1-j} \bigr)
+ \bigl( \tfrac{1}{2} S - x_{i+1, n+1-j} \bigr)
+ \bigl( \tfrac{1}{2} S - x_{i, n-j} \bigr) + \bigl( \tfrac{1}{2} S - x_{i+1, n-j} \bigr) \\
&= 2S - \bigl(  x_{i, n+1-j} +  x_{i+1, n+1-j} +  x_{i, n-j} +  x_{i+1, n-j} \bigr) =S.
\end{align*}

Finally, we need to show that every integer $k$, for $1\le k \le mn$, that $k$ is used exactly once in the labeling $X$.
From conditions (1) and (2) of Proposition \ref{propKleinCheckerboardMOddNEvenC4FaceExists}
and the identity $x_{i,n+1-j} = \tfrac{1}{2} S - x_{i,j}$ of Lemma \ref{lemmaKleinCheckerboardSemiC4FaceValueBalance},
we have
\begin{align*}
x_{2i-1,2j-1} &= (j-1)m+i, \mbox{ \ \ \ for } 1 \le i \le m_1 +1 \mbox{ and } 1 \le j \le \stkout{n}_{1}^{+}, \mbox{ and } \\
x_{2i,n+2-2j} &=  (j-1)m+ m_1 +1 + i, \mbox{ \ \ \ for } 1 \le  i \le m_1 \mbox{ and } 1 \le j \le \stkout{n}_{1}^{+}.
\end{align*}
This produces the list of integers $k$, for $1\le k \le m \stkout{n}_{1}^{+}$, where each integer $k$ is used
exactly once in the labeling $X$.

From conditions (3) and (4) of Proposition \ref{propKleinCheckerboardMOddNEvenC4FaceExists}
and the identity $x_{i,n+1-j} = \tfrac{1}{2} S - x_{i,j}$ of Lemma \ref{lemmaKleinCheckerboardSemiC4FaceValueBalance},
we have
\begin{align*}
x_{2i-1,2j} &= (n_1 -j +1 )m + 1 -i, \mbox{ \ \ \ for } 1 \le i \le m_1 +1 \mbox{ and } 1 \le j \le \stkout{n}_{1}, \mbox{ and } \\
x_{2i,n+1-2j} &=  (n_1 - j + 1)m -  m_1 - i, \mbox{ \ \ \ for } 1 \le  i \le m_1 \mbox{ and } 1 \le j \le \stkout{n}_{1}.
\end{align*}
This produces the list of integers $k$, for $m \stkout{n}_{1}^{+} + 1 \le k \le m n_1$, where each integer $k$ is used
exactly once in the labeling $X$.

When we include the other values of the labeling $X$ that result from using the identity
$x_{i,n+1-j} = \tfrac{1}{2} S - x_{i,j}$ of Lemma \ref{lemmaKleinCheckerboardSemiC4FaceValueBalance},
we produce the list of integers $k$, for $m n_1 + 1 \le k \le m n$, where each integer $k$ is used
exactly once in the labeling $X$.
This completes the proof.
\end{proof}

\begin{example}
Table~\ref{tableK7x10labelingexample} illustrates Proposition~\ref{propKleinCheckerboardMOddNEvenC4FaceExists}, when $m =7$ and $n = 10$.
\end{example}

\begin{table}[h]
\hspace{0.0in}\begin{tabular}{|c|c|c|c|c|c|c|c|c|c|}
\hline
1   &35  &8    &28  &15  &56  &43  &63  &36  &70  \\[0.05in]
\hline
66  &40  &59   &47  &52  &19  &24  &12  &31  &5   \\[0.05in]
\hline
2   &34  &9    &27  &16  &55  &44  &62  &37  &69  \\[0.05in]
\hline
65  &41  &58   &48  &51  &20  &23  &13  &30  &6  \\[0.05in]
\hline
3   &33  &10   &26  &17  &54  &45  &61  &38  &68   \\[0.05in]
\hline
64  &42  &57   &49  &50  &21  &22  &14  &29  &7  \\[0.05in]
\hline
4   &32  &11   &25  &18  &53  &46  &60  &39  &67   \\[0.05in]
\hline
\end{tabular}
\caption{A $C_4$-face-magic Klein bottle labeling on $\mathcal{K}_{7,10}$.}
\label{tableK7x10labelingexample}
\end{table}

\begin{proposition} \label{propKleinCheckerboardMEvenNEvenC4FaceExists}
Let $m$ and $n$ be even positive integers.
Let $n=2 n_1$ for some positive integer $n_1$.
Define a labeling $X = \{ x_{i,j} : (i,j)\in V( \mathcal{K}_{m,n}) \}$  on $\mathcal{K}_{m,n}$ as follows:
\begin{itemize}
\item $x_{2i-1,2j-1} = 2m(j-1)+i$, for $i=1,2,\ldots,\stkout{m}$ and $j=1,2,\ldots,\stkout{n}_1^{+}$.
\item $x_{2i-1,2j} = 2mj +1 -i$, for $i=1,2,\ldots,\stkout{m}$ and $j=1,2,\ldots,\stkout{n}_1$.
\item $x_{2i-1,2j-1} = m(n +1 - 2j)+i$, for $i=1,2,\ldots,\stkout{m}$ and $j=\stkout{n}_{1}^{+} +1,\stkout{n}_{1}^{+} +2 ,\ldots,n_1$.
\item $x_{2i-1,2j} = m(n +1-2j)+1 -i$, for $i=1,2,\ldots,\stkout{m}$ and $j=\stkout{n}_{1} +1,\stkout{n}_{1} +2 ,\ldots,n_1$.
\end{itemize}
Furthermore for $1\leq i \leq \stkout{m}$ and $1 \leq j\leq n$, let
$x_{2i,j} = mn+1-x_{2i-1,j}$.
Then, $X$ is a $C_4$-face-magic Klein bottle labeling on $\mathcal{K}_{m,n}$.
\end{proposition}

\begin{proof}
Since $x_{2i,j}=mn+1-x_{2i-1,j}$, for $1\leq i\leq \stkout{m}$
and $1\leq j\leq n$, we know that $%
x_{2i,j}+x_{2i-1,j}+x_{2i,j+1}+x_{2i-1,j+1}=2mn+2$, for $1\leq i\leq \stkout{m}$
and $1\leq j\leq n-1$. Also, $x_{2i-1,n}+x_{2i,n}+x_{2i-1,1}+x_{2i,1}=2mn+2$.
\vspace{0.1in}

Let $1\leq i\leq \stkout{m}-1$ and $1\leq j\leq n$. Note that for $j$ odd, $%
x_{2i+1,j}=x_{2i-1,j}+1$, and for $j$ even, $x_{2i+1,j}=x_{2i-1,j}-1$.
Hence, for $j<n$,
$x_{2i,j}+x_{2i+1,j}+x_{2i,j+1}+x_{2i+1,j+1}=x_{2i,j}+x_{2i-1,j}+x_{2i,j+1}+x_{2i-1,j+1}=2mn+2$.
This also holds for $j=n$, with the understanding that $j+1=1$ in this
case.
\vspace{0.1in}

It is now necessary to calculate the sums involving $x_{1,j}$ and $x_{m,j}$,
for $1\leq j<n$.
If $j$ is odd, $x_{1,j}-x_{1,n+1-j}=1-m$. If $j$ is even, $x_{1,j}-x_{1,n+1-j}=m-1$.
Thus, for all integers $1 \le j < n$, we have
\begin{equation*}
(x_{1,j}-x_{1,n+1-j}) + (x_{1,j+1}-x_{1,n-j}  ) =0.
\end{equation*}
If $j$ is odd, $x_{m,j}=mn+1-x_{m-1,j}=mn+1-\left( x_{1,j}+ \, \stkout{m} \, -1\right)
=mn- \, \stkout{m} \, -x_{1,j}+2$. If $j$ is even, $x_{m,j}=mn+1-x_{m-1,j}=mn+1-\left(
x_{1,j}- \, \stkout{m} \, +1\right) =mn+ \, \stkout{m} \, -x_{1,j}$.
Thus, for all integers $1 \le j < n$, we have
\begin{equation*}
(x_{m,n+1-j} + x_{m,n-j}) = 2mn+2 - x_{1,n+1-j}-x_{1,n-j}.
\end{equation*}
Hence,
\begin{align*}
x_{1,j} &+x_{1,j+1}+x_{m,n+1-j}+x_{m,n-j}  \\
       &=x_{1,j}+x_{1,j+1}+  (  2mn+2 - x_{1,n+1-j}-x_{1,n-j}  ) \\
        &=2mn+2+ (  x_{1,j} -x_{1,n+1-j} ) + (x_{1,j+1} -x_{1,n-j} ) =2mn+2.
\end{align*}
For the case when $j=n$, we have
\begin{align*}
x_{1,1}+x_{1,n}+x_{m,1}+x_{m,n} = 1+m + (\tfrac{1}{2} S - \, \stkout{m} \, ) + ( \tfrac{1}{2} S - \, \stkout{m} \, -1) =S = 2mn+2.
\end{align*}

Finally, we need to show that every integer $k$, for $1\le k \le mn$, that $k$ is used exactly once in the labeling $X$.
A relabeling of the indices $i$ and $j$ in the expressions for $x_{2i-1,j}$ yields the following formulas:
\begin{align}
x_{2i-1,2j-1} &= (2j-2)m +i, \mbox{ \ \ \ for } 1\le i \le \stkout{m} \mbox{ and }
1 \le j \le \stkout{n}_1^{+},  \label{eqnReindexFormulaForX1}\\
x_{2( \, \stkout{m} \, +1-i)-1,2j} &= (2j)m - \, \stkout{m} \, +i, \mbox{ \ \ \ for }
1 \le i \le \stkout{m}   \mbox{ and }  1 \le j \le \stkout{n}_1, \\
x_{2i-1,2(n_1 +1-j)-1} &= (2j-1)m +i, \mbox{  for } 1 \le i \le \stkout{m}
\mbox{ and } 1 \le j \le \stkout{n}_1, \mbox{ and} \\
x_{2( \stkout{m} \, +1-i)-1,2(n_1 +1-j)} &= (2j-1)m -  \stkout{m}  +i,  \mbox{ for }
1 \le i \le \stkout{m}  \mbox{ and } 1 \le j \le \stkout{n}_1^{+}. \label{eqnReindexFormulaForX4}
\end{align}
Let $1 \le j \le \stkout{n}_1$. The sequence of labels
\begin{align*}
\bigl(  x_{2i-1,2j-1} : 1 \le i \le \stkout{m} \bigr) &= \bigl(  (2j-2)m +i  : 1 \le i \le \stkout{m} \bigr) \\
  &= \bigl( (2j-2)m +1, (2j-2)m +2, \ldots, (2j-2)m + \, \stkout{m} \bigr),
\end{align*}
followed by the sequence of labels
\begin{align*}
\bigl(  &x_{2( \, \stkout{m} \, +1-i)-1,2(n_1 +1-j)}  : 1 \le i \le \stkout{m} \bigr)
   = \bigl(  (2j-1)m - \, \stkout{m} \, +i  : 1 \le i \le \stkout{m} \bigr) \\
  &= \bigl( (2j-2)m  + \, \stkout{m}+1, (2j-2)m  + \, \stkout{m}+2, \ldots, (2j-1)m  \bigr),
\end{align*}
followed by the sequence of labels
\begin{align*}
\bigl(  x_{2i-1,2(n_1 +1-j)-1}  &: 1 \le i \le \stkout{m} \bigr)
   = \bigl(  (2j-1)m +i,  : 1 \le i \le \stkout{m} \bigr) \\
  &= \bigl( (2j-1)m  +1, (2j-1)m  +2, \ldots, (2j-1)m + \,  \, \stkout{m} \bigr),
\end{align*}
and then followed by the sequence of labels
\begin{align*}
\bigl(  &x_{2( \, \stkout{m} \, +1-i)-1,2j} : 1 \le i \le \stkout{m} \bigr)
   = \bigl(  (2j)m - \, \stkout{m} \, +i  : 1 \le i \le \stkout{m} \bigr) \\
  &= \bigl( (2j-1)m  + \, \stkout{m}+1, (2j-1)m  + \, \stkout{m}+2, \ldots, (2j)m  \bigr),
\end{align*}
yields the sequence of integers $\bigl( k : (2j-1)m+1 \le k \le (2j)m \bigr)$.
When we concatenate these sequences of integers together for $1\le j \le \, \stkout{n}_1$,
we obtain the sequence $\bigl( k : 1\le k \le 2m\, \stkout{n}_1  \bigr)$.
When $n_1$ is even, this is the sequence $\bigl( k : 1\le k \le m n_1  \bigr)$.
However, when $n_1$ is odd, this is the sequence $\bigl( k : 1\le k \le m (n_1 -1)  \bigr)$.
When we let $j=\stkout{n}_1^{+}= \tfrac{1}{2} (n_1 +1)$ in equations (\ref{eqnReindexFormulaForX1})
and (\ref{eqnReindexFormulaForX4}), the sequence of labels
\begin{align*}
\bigl(  x_{2i-1,n_1} : 1 \le i \le \stkout{m} \bigr) &= \bigl(  (n_1-1)m +i  : 1 \le i \le \stkout{m} \bigr) \\
  &= \bigl( (n_1-1)m +1, (n_1-1)m +2, \ldots, (n_1-1)m + \, \stkout{m} \bigr),
\end{align*}
followed by the sequence of labels
\begin{align*}
\bigl(  &x_{2( \, \stkout{m} \, +1-i)-1,n_1 +1}  : 1 \le i \le \stkout{m} \bigr)
   = \bigl(  (n_1 -1)m - \, \stkout{m} \, +i  : 1 \le i \le \stkout{m} \bigr) \\
  &= \bigl( (n_1 -1)m  + \, \stkout{m}+1, (n_1 -1)m  + \, \stkout{m}+2, \ldots, (n_1)m  \bigr),
\end{align*}
produces the sequence of integers $\bigl( k : m(n_1 -1) +1 \le k \le mn_1 \bigr)$.
Therefore, the set of labels $\{ x_{2i-1,j} : 1\le i \le \, \stkout{m} \mbox{ and } 1\le j \le n\}$
is the set of integers $\{ k : 1\le k \le mn_1\}$.
From the identity $x_{2i,j} = mn+1 -x_{2i-1,j}$, the
set of labels $\{ x_{2i,j} : 1\le i \le \, \stkout{m} \mbox{ and } 1\le j \le n\}$
is the set of integers $\{ k : mn_1 +1\le k \le mn\}$. This completes the proof.
\end{proof}

\begin{example}
Table~\ref{tableK6x10labelingexample} illustrates Proposition~\ref{propKleinCheckerboardMEvenNEvenC4FaceExists}, when $m = 6$ and $n = 10$.
\begin{table}[h]
\hspace{0.0in}\begin{tabular}{|c|c|c|c|c|c|c|c|c|c|}
\hline
1   &12  &13   &24  &25  &30  &19  &18  &7   &6   \\[0.05in]
\hline
60  &49  &48   &37  &36  &31  &42  &43  &54  &55  \\[0.05in]
\hline
2   &11  &14   &23  &26  &29  &20  &17  &8   &5   \\[0.05in]
\hline
59  &50  &47   &38  &35  &32  &41  &44  &53  &56 \\[0.05in]
\hline
3   &10  &15   &22  &27  &28  &21  &16  &9   &4    \\[0.05in]
\hline
58  &51  &46   &39  &34  &33  &40  &45  &52  &57 \\[0.05in]
\hline
\end{tabular}
\caption{A $C_4$-face-magic Klein bottle labeling on $\mathcal{K}_{6,10}$.}
\label{tableK6x10labelingexample}
\end{table}
\end{example}

\begin{theorem} \label{thmC4FaceMagicLabelingExists}
Let $m$ and $n$ be positive integers such that $m,n \geq 2$.
Then, $\mathcal{K}_{m,n}$ has a $C_4$-face-magic Klein bottle labeling if and only if $n$ is even.
\end{theorem}

\begin{proof}
\noindent ($\Rightarrow$) Suppose that $\mathcal{K}_{m,n}$ is a $C_4$-face-magic Klein bottle graph.
Then by Proposition \ref{propKleinCheckerboardNIsEven}, $n$ is even.

\

\noindent ($\Leftarrow$) Assume $n$ is even. On the one hand, if $m$ is odd, by Proposition \ref{propKleinCheckerboardMOddNEvenC4FaceExists},
$\mathcal{K}_{m,n}$ has a $C_4$-face-magic Klein bottle labeling.
On the other hand, if $m$ is even, by Proposition \ref{propKleinCheckerboardMEvenNEvenC4FaceExists},
$\mathcal{K}_{m,n}$ has a $C_4$-face-magic Klein bottle labeling.
\end{proof}

\section{Equivalent $C_4$-face-magic Klein bottle labelings} \label{S:equivalentlabelings}

\begin{definition}
Consider the natural embedding of $\mathcal{K}_{m,n}$ on the Klein bottle.
We say that two $C_4$-face-magic Klein bottle labelings
on $\mathcal{K}_{m,n}$ are {\it Klein bottle equivalent} if
there is a homeomorphism of the Klein bottle that maps $\mathcal{K}_{m,n}$ onto itself such that
the first $C_4$-face-magic Klein bottle labeling on $\mathcal{K}_{m,n}$ is mapped to the
second $C_4$-face-magic Klein bottle labeling on $\mathcal{K}_{m,n}$.
\end{definition}

\begin{definition} \label{notnKleinBottleSymmetries}
We consider the following symmetries of the Klein bottle that result in a graph automorphism on $\mathcal{K}_{m,n}$.
Let $U$ be the translation on $\mathcal{K}_{m,n}$ given by
\[
U(i,j) = \left\{\begin{array}{ll}
              (i+1,j), &\mbox{if \ } 1\leq i \leq m-1 \mbox{ and } 1\leq j\leq n, \\
              (1,n+1-j), &\mbox{if \ }  i = m \mbox{ and } 1\leq j\leq n. \\
\end{array} \right. \hspace{4.5in}
\]
Let $H$ be the translation on $\mathcal{K}_{m,n}$ given by $H(i,j) = (i,j+\stkout{n})$.
Let $F$ be the reflection on $\mathcal{K}_{m,n}$ given by $F(i,j) = (i,n+1-j)$.
We observe that $U^m = F$. Also, the symmetries of the Klein bottle that result in a graph automorphism
of $\mathcal{K}_{m,n}$ is the group of graph automorphisms
$KBLS(m,n) = \langle U,H,F\rangle = \{ U^i, FU^i, HU^i, FHU^i : 0\leq i \leq m-1\} = \{ U^i, HU^i  : 0\leq i \leq 2m-1\}$.
Given a $C_4$-face-magic Klein bottle labeling $X =\{ x_{i,j} : (i,j)\in V(\mathcal{K}_{m,n})\}$ and $A\in KBLS(m,n)$,
we let $A(X) =\{ x_{A(i,j)} : (i,j)\in V(\mathcal{K}_{m,n})\}$.
We call $KBLS(m,n)$ the {\it Klein bottle labeling symmetry group on } $\mathcal{K}_{m,n}$.
\end{definition}

\begin{definition} \label{defnGraphToSearchForLabelings}
Let $m$ be an odd positive integer such that $m\geq 3$, $n$ be an even positive integer such that $n\geq 4$,
$\stkout{n}=n/2$, and $n_0 =\stkout{n} -1$.
Let  $a_1,a_2,\dots,a_{n_0}$ be positive integers.
We define the \textit{labeling search graph corresponding to} $(a_1,a_2,\ldots,a_{n_0})$,
denoted by $G(a_1,a_2,\ldots,a_{n_0})$, to be the graph whose vertex set is
\begin{align*}
V( G(a_1,a_2,\ldots,a_{n_0}) ) = \{ \{ q, mn+1-q\} : 1\leq q \leq  m \stkout{n} \}
\end{align*}
and there is an edge from $\{ x_1, x_2\}$ to $\{ y_1,y_2\}$, where $\{ x_1, x_2 \} \ne  \{ y_1, y_2\}$,  if and only if there exist
$z_1\in \{ x_1, x_2\}$, $z_2\in \{ y_1, y_2\}$,
and $1\leq j \leq n_0$ such that $z_1 + z_2 = a_j$.
\vspace{0.1in}

Let $P = ( \{ x_{j,1}, x_{j,2}\} : 1\le j \le n_1 )$ be a path on $n_1$ distinct vertices in $G(a_1,a_2,\ldots,a_{n_0})$.
We say that $P$ is an \textit{$(a_1,a_2,\ldots,a_{n_0})$-admissible path }
provided that, for all $1 \le j \le n_1$, there exist $z_j \in \{ x_{j,1}, x_{j,2}\}$,
such that, for all $1\le j \le n_0$, $z_j + z_{j+1} =a_j$.
We call $z_1$ and $z_{n_1}$ \textit{end vertex labels} of $P$.
For all $1 < j < n_1$, we call $z_j$ an \textit{interior vertex label} of $P$.
We say that $z_1$ is an \textit{end vertex label} of $P$ \textit{incident to the edge labeled} $a_1$.
We say that $z_{n_1}$ is an \textit{end vertex label} of $P$ \textit{incident to the edge labeled} $a_{n_0}$.
For all $1 < j < n_1$, we say that $z_j$ is an \textit{interior vertex label}
of $P$ \textit{incident to the edges labeled} $a_{j-1}$ and $a_j$.
We call the sequence $(z_j: 1\le j \le n_1)$ the \textit{vertex label sequence} of $P$.
\vspace{0.1in}

Let $H$ be a spanning subgraph of $G(a_1,a_2,\ldots,a_{n_0})$.
We say that $H$ is an \textit{$(a_1,a_2,\ldots,a_{n_0})$-admissible path partition} of $G(a_1,a_2,\ldots,a_{n_0})$
if $H$ is the disjoint union of $m$ distinct $(a_1,a_2,\ldots,a_{n_0})$-admissible paths in $G(a_1,a_2,\ldots,a_{n_0})$.
\end{definition}

\begin{remark} \label{remGraphGaHasNoLoops}
The condition $\{ x_1, x_2 \} \ne  \{ y_1, y_2\}$ in Definition \ref{defnGraphToSearchForLabelings}
on vertices $\{ x_1, x_2 \}$ and $\{ y_1, y_2\}$ ensures that the graph
$G(a_1,a_2,\ldots,a_{n_0})$ does not have any loops.
Suppose integer $a_j$ is even, for some $1\le j \le n_0$.
Let $a_j =2k$ for some positive integer $k$.
Then, vertex $\{ k, mn+1-k\}$ does not have a loop labeled $a_j$ incident to it in $G(a_1,a_2,\ldots,a_{n_0})$.
\end{remark}

\begin{lemma} \label{lemmaSequenceOfPairSums}
Let $m\geq 3$ be an odd integer and $n\geq 4$ be an even integer, $\stkout{n}=n/2$, and $n_0 =\stkout{n} -1$.
Let $S=2(mn+1)$.
Let $X= \{ x_{i,j} : (i,j)\in V(\mathcal{K}_{m,n})\}$ be a labeling on $\mathcal{K}_{m,n}$
such that each integer $k$ with $1\le k \le mn$ is used exactly once.
For all $1\le j \le n_0$, let $a_j =x_{1,j} + x_{1,j+1}$.
Then, $X= \{ x_{i,j} : (i,j)\in V(\mathcal{K}_{m,n})\}$ is a $C_4$-face-magic Klein bottle labeling on $\mathcal{K}_{m,n}$
if and only if
\begin{itemize}
\item for all $1\leq i \leq \stkout{m}^{+}$ and $1\leq j \le n_0$, we have $x_{2i-1,j} + x_{2i-1,j+1}=a_j$,
\item for all $1\leq i \leq \stkout{m}^{+}$ and $1\leq j \le \, \stkout{n}$, we have $x_{2i-1,n+1-j} = \tfrac{1}{2} S - x_{2i-1,j}$,
\item for all $1\leq i \leq \stkout{m}$ and $1\leq j \le n_0$, we have $x_{2i,n-j} + x_{2i,n-j+1}=a_j$, and
\item for all $1\leq i \leq \stkout{m}$ and $1\leq j \le \, \stkout{n}$, we have $x_{2i,j} = \tfrac{1}{2} S - x_{2i,n+1-j}$.
\end{itemize}
\end{lemma}

\begin{proof}
\noindent ($\Rightarrow$)
Assume that $X= \{ x_{i,j} : (i,j)\in V(\mathcal{K}_{m,n})\}$ is a $C_4$-face-magic Klein bottle labeling on $\mathcal{K}_{m,n}$.
When we equate the two $C_4$-face sums given below equal to each other
\begin{align*}
x_{2i-1,j} &+ x_{2i-1,j+1} + x_{2i,j} + x_{2i,j+1} = S \\
   &=  x_{2i,j} + x_{2i,j+1} +  x_{2i+1,j} + x_{2i+1,j+1},
\end{align*}
we obtain
\begin{equation*}
x_{2i-1,j} + x_{2i-1,j+1}  =   x_{2i+1,j} + x_{2i+1,j+1}, \mbox{ \ \ \ for } 1\leq i  \leq \stkout{m}.
\end{equation*}
Thus, $x_{2i-1,j} + x_{2i-1,j+1}  =   x_{1,j} + x_{1,j+1} = a_j$, for  $1\leq i  \leq \stkout{m}^{+}$.
By Lemma \ref{lemmaKleinCheckerboardSemiC4FaceValueBalance},
for all $1\leq i \leq \stkout{m}^{+}$ and $1\leq j \le \, \stkout{n}$, we have $x_{2i-1,n+1-j} = \tfrac{1}{2} S - x_{2i-1,j}$.

When we equate the two $C_4$-face sums given below equal to each other
\begin{equation*}
x_{m,j} + x_{m,j+1} + x_{1,n-j} + x_{1,n-j+1} = S =  x_{1,n-j} + x_{1,n-j+1} +  x_{2,n-j} + x_{2,n-j+1},
\end{equation*}
we obtain
\begin{equation*}
x_{2,n-j} + x_{2,n-j+1}  =   x_{m,j} + x_{m,j+1} = a_j.
\end{equation*}

When we equate the two $C_4$-face sums given below equal to each other
\begin{align*}
x_{2i,n-j} &+ x_{2i,n-j+1} + x_{2i+1,n-j} + x_{2i+1,n-j+1} = S \\
           &=  x_{2i+1,n-j} + x_{2i+1,n-j+1} +  x_{2i+2,n-j} + x_{2i+2,n-j+1},
\end{align*}
we obtain
\begin{equation*}
x_{2i,n-j} + x_{2i,n-j+1}  =   x_{2i+2,n-j} + x_{2i+2,n-j+1}, \mbox{ \ \ \ for } 1\leq i  \leq \stkout{m} -1.
\end{equation*}
Thus, $x_{2i,n-j} + x_{2i,n-j+1}  =   x_{2,n-j} + x_{2,n-j+1} = a_j$, for  $1\leq i  \leq \stkout{m}$.
By Lemma \ref{lemmaKleinCheckerboardSemiC4FaceValueBalance},
for all $1\leq i \leq \stkout{m}$ and $1\leq j \le \, \stkout{n}$, we have $x_{2i,j} = \tfrac{1}{2} S - x_{2i,n+1-j}$.
\vspace{0.1in}

\noindent ($\Leftarrow$)
Let $X= \{ x_{i,j} : (i,j)\in V(\mathcal{K}_{m,n})\}$ be a labeling on $\mathcal{K}_{m,n}$
such that each integer $k$ with $1\le k \le mn$ is used exactly once;
for all $1\leq i \leq \stkout{m}^{+}$ and $1\leq j \le n_0$, we have
\begin{equation} \label{eqnNecessaryConditionKleinLabeling1}
x_{2i-1,j} + x_{2i-1,j+1}=a_j;
\end{equation}
for all $1\leq i \leq \stkout{m}^{+}$ and $1\leq j \le n_1$, we have
\begin{equation} \label{eqnNecessaryConditionKleinLabeling2}
x_{2i-1,n+1-j} = \tfrac{1}{2} S - x_{2i-1,j};
\end{equation}
for all $1\leq i \leq \stkout{m}$ and $1\leq j \le n_0$, we have
\begin{equation} \label{eqnNecessaryConditionKleinLabeling3}
x_{2i,n-j} + x_{2i,n-j+1}=a_j, \mbox{ and }
\end{equation}
for all $1\leq i \leq \stkout{m}$ and $1\leq j \le n_1$, we have
\begin{equation} \label{eqnNecessaryConditionKleinLabeling4}
x_{2i,j} = \tfrac{1}{2} S - x_{2i,n+1-j}.
\end{equation}
We show that $X= \{ x_{i,j} : (i,j)\in V(\mathcal{K}_{m,n})\}$ is a $C_4$-face-magic Klein bottle labeling on $\mathcal{K}_{m,n}$.

\textbf{Case 1.} Let $1\leq i \leq \stkout{m}$ and $1\leq j \le n_0$. We show that
\begin{equation*}
x_{2i-1,j} + x_{2i-1,j+1} + x_{2i,j} + x_{2i,j+1} = S.
\end{equation*}
From equation (\ref{eqnNecessaryConditionKleinLabeling1}), we have $x_{2i-1,j} + x_{2i-1,j+1}=a_j$.
From equation (\ref{eqnNecessaryConditionKleinLabeling4}), we have $x_{2i,j} = \tfrac{1}{2} S - x_{2i,n+1-j}$
and $x_{2i,j+1} = \tfrac{1}{2} S - x_{2i,n-j}$.
Also, from equation (\ref{eqnNecessaryConditionKleinLabeling3}), we have $x_{2i,n-j} + x_{2i,n-j+1}=a_j$.
Thus,
\begin{align*}
  x_{2i-1,j} + x_{2i-1,j+1} + x_{2i,j} + x_{2i,j+1} &=
x_{2i-1,j} + x_{2i-1,j+1} \\
 &+ (\tfrac{1}{2} S - x_{2i,n+1-j})
  + (\tfrac{1}{2} S - x_{2i,n-j}) \\
  &= a_j +S -a_j =S.
\end{align*}
A similar argument shows that
\begin{equation*}
x_{2i,j} + x_{2i,j+1} + x_{2i+1,j} + x_{2i+1,j+1} = S.
\end{equation*}

\textbf{Case 2.} Let $1\leq i  \le m$. We show that
\begin{align*}
 x_{i,1} + x_{i+1,1} +  x_{i,n} + x_{i+1,n} &= S \mbox{ and }  \\
 x_{i,\stkout{n}} + x_{i+1,\stkout{n}} +  x_{i,\stkout{n} \, + 1} + x_{i+1,\stkout{n} \, +1} &= S,
\end{align*}
where we take the index $i$ modulo $m$.
From equations (\ref{eqnNecessaryConditionKleinLabeling2}) and (\ref{eqnNecessaryConditionKleinLabeling4}), we have
$x_{i,1} = \tfrac{1}{2} S - x_{i,n}$ and $x_{i+1,1} = \tfrac{1}{2} S - x_{i+1,n}$.
Thus,
\begin{align*}
  x_{i,1} + x_{i+1,1} + x_{i,n} + x_{i+1,n} =   (\tfrac{1}{2} S - x_{i,n})
  + (\tfrac{1}{2} S - x_{i+1,n}) + x_{i,n} + x_{i+1,n} =S.
\end{align*}
A similar argument shows that
\begin{align*}
 x_{i,\stkout{n}} + x_{i+1,\stkout{n}} +  x_{i,\stkout{n} \, + 1} + x_{i+1,\stkout{n} \, +1} = S.
\end{align*}

\textbf{Case 3.} Let $1\leq j \le n_0$.
From equation (\ref{eqnNecessaryConditionKleinLabeling1}), we have $x_{1,j} + x_{1,j+1}=a_j$ and $x_{m,j} + x_{m,j+1}=a_j$.
From equation (\ref{eqnNecessaryConditionKleinLabeling2}), we have
$x_{1,n+1-j} = \tfrac{1}{2} S - x_{1,j}$ and $x_{1,n-j} = \tfrac{1}{2} S - x_{1,j+1}$.
Thus,
\begin{align*}
  x_{m,j} + x_{m,j+1} + x_{1,n+1-j} + x_{1,n-j} &=
x_{m,j} + x_{m,j+1} \\
 &+ (\tfrac{1}{2} S - x_{1,j})
  + (\tfrac{1}{2} S - x_{1,j+1}) \\
  &= a_j +S -a_j =S.
\end{align*}
A similar argument shows that, for $\stkout{n} \, +1 \le j <n$, we have
\begin{align*}
  x_{m,j} + x_{m,j+1} + x_{1,n+1-j} + x_{1,n-j} = S.
\end{align*}

Since every integer $k$ with $1\le k \le mn$ is used in the labeling $X$,
$X$ is a $C_4$-face-magic Klein bottle labeling on $\mathcal{K}_{m,n}$.
\end{proof}

\begin{definition} \label{defnLabelingGraphOfX}
Let $m\geq 3$ be an odd integer and $n\geq 4$ be an even integer, $\stkout{n}=n/2$, and $n_0 =\stkout{n} -1$.
Let $X= \{ x_{i,j} : (i,j)\in V(\mathcal{K}_{m,n})\}$ be a $C_4$-face-magic Klein bottle labeling on $\mathcal{K}_{m,n}$.
We define the {\it labeling graph of X}, denoted by $\mathcal{L}(X)$, to be the graph with vertex set
\begin{align*}
V(\mathcal{L}(X)) = \{ \{q, mn+1-q\} : 1\leq q \leq m \stkout{n} \}
\end{align*}
and consists of the following edges:
\begin{itemize}
\item For all $1\leq i \leq \stkout{m}^{+}$ and $1\le j \le n_0$, there is an edge labeled $a_j =x_{2i-1,j} + x_{2i-1,j+1}$
from $\{ x_{2i-1,j}, mn+1-x_{2i-1,j} \}$ to $\{ x_{2i-1,j+1}, mn+1-x_{2i-1,j+1} \}$. The path
$( \{ x_{2i-1,j}, mn+1-x_{2i-1,j} \} : 1\le j \le \stkout{n})$ is an $(a_1,a_2,\ldots,a_{n_0})$-admissible path in
$G(a_1,a_2,\ldots,a_{n_0})$ with vertex label sequence $( x_{2i-1,j} : 1\le j \le \stkout{n})$.
\item For all $1\leq i \leq \stkout{m}$ and $1\le j \le n_0$, there is an edge labeled $a_j =x_{2i,n+1-j} + x_{2i,n-j}$
from $\{ x_{2i,n+1-j}, mn+1-x_{2i,n+1-j} \}$ to $\{ x_{2i,n-j}, mn+1-x_{2i,n-j} \}$.
The path $( \{ x_{2i,n+1-j}, mn+1-x_{2i,n-j} \} : 1\le j \le \stkout{n})$ is an $(a_1,a_2,\ldots,a_{n_0})$-admissible path in
$G(a_1,a_2,\ldots,a_{n_0})$ with vertex label sequence $( x_{2i,n+1-j} : 1\le j \le \stkout{n})$.
\end{itemize}
By Lemma \ref{lemmaSequenceOfPairSums}, the graph $\mathcal{L}(X)$ is well-defined.
Also, $\mathcal{L}(X)$ is an $(a_1,a_2,\ldots,a_{n_0})$-admissible path partition of $G(a_1,a_2,\ldots,a_{n_0})$.
\end{definition}

\begin{figure}
\hspace{0.0in}\begin{picture}(300,150)(-50,0)
\put(10,120){\circle*{7}}
\put(-5,100){$\{ \textbf{1},30\}$}
\put(10,120){\line(1,0){40}}
\put(25,130){22}

\put(50,120){\circle*{7}}
\put(35,100){$\{ 10, \textbf{21}\}$}
\put(50,120){\line(1,0){40}}
\put(65,130){28}

\put(90,120){\circle*{7}}
\put(75,100){$\{ \textbf{7},24\}$}
\put(90,120){\line(1,0){40}}
\put(105,130){34}

\put(130,120){\circle*{7}}
\put(115,100){$\{ 4,\textbf{27}\}$}
\put(130,120){\line(1,0){40}}
\put(145,130){40}

\put(170,120){\circle*{7}}
\put(155,100){$\{ \textbf{13},18\}$}

\put(10,70){\circle*{7}}
\put(-5,50){$\{ \textbf{2},29\}$}
\put(10,70){\line(1,0){40}}
\put(25,80){22}

\put(50,70){\circle*{7}}
\put(35,50){$\{ 11,\textbf{20}\}$}
\put(50,70){\line(1,0){40}}
\put(65,80){28}

\put(90,70){\circle*{7}}
\put(75,50){$\{ \textbf{8},23\}$}
\put(90,70){\line(1,0){40}}
\put(105,80){34}

\put(130,70){\circle*{7}}
\put(115,50){$\{ 5,\textbf{26}\}$}
\put(130,70){\line(1,0){40}}
\put(145,80){40}

\put(170,70){\circle*{7}}
\put(155,50){$\{ \textbf{14},17\}$}

\put(10,20){\circle*{7}}
\put(-5,0){$\{\textbf{3},28\}$}
\put(10,20){\line(1,0){40}}
\put(25,30){22}

\put(50,20){\circle*{7}}
\put(35,0){$\{ 12,\textbf{19}\}$}
\put(50,20){\line(1,0){40}}
\put(65,30){28}

\put(90,20){\circle*{7}}
\put(75,0){$\{ \textbf{9},22\}$}
\put(90,20){\line(1,0){40}}
\put(105,30){34}

\put(130,20){\circle*{7}}
\put(115,0){$\{ 6,\textbf{25}\}$}
\put(130,20){\line(1,0){40}}
\put(145,30){40}

\put(170,20){\circle*{7}}
\put(155,0){$\{ \textbf{15},16\}$}
\end{picture}
\caption{$(22,28,34,40)$-admissible path partition $K$ in $G(22,28,34,40)$.}
\label{figGraphProp33Case2K1at3by10}
\end{figure}
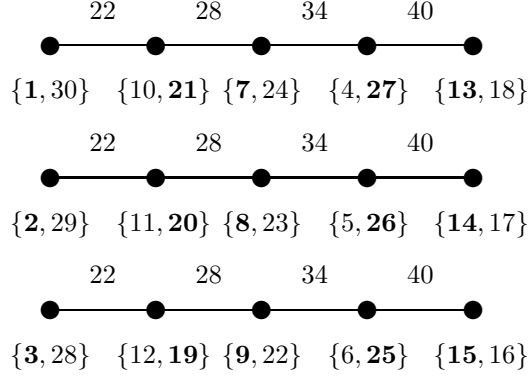

\begin{table}[h]
\hspace{0.0in}\begin{tabular}{|c|c|c|c|c|c|c|c|c|c|}
\hline
1   &21  &7    &27  &13  &18  &4  &24  &10  &30  \\[0.05in]
\hline
29  &11  &23   &5  &17  &14  &26  &8  &20  &2   \\[0.05in]
\hline
3   &19  &9   &25  &15  &16  &6  &22  &12  &28   \\[0.05in]
\hline
\end{tabular}
\caption{$C_4$-face-magic Klein bottle labeling $X$ on $\mathcal{K}_{3,10}$ such that $\mathcal{L}(X)=K$ is the $(22,28,34,40)$-admissible path partition $K$ in $G(22,28,34,40)$ given in Figure \ref{figGraphProp33Case2K1at3by10}.}
\label{tableK3by10here}
\end{table}

\begin{example}
In Table \ref{tableK3by10here} we have a $C_4$-face-magic Klein bottle labeling $X$ on the $3\times 10$ Klein bottle grid graph $\mathcal{K}_{3,10}$ such that
$\mathcal{L}(X)=K$ is the $(22,28,34,40)$-admissible path partition $K$ in the labeling search graph $G(22,28,34,40)$ given in Figure \ref{figGraphProp33Case2K1at3by10}.
\end{example}

\begin{lemma} \label{lemmaSpanningGraphOfPathsCriteria}
Let $m\geq 3$ be an odd integer, $n \geq 4$ be an even integer, $n_1 = \stkout{n}=n/2$, and $n_0 =\stkout{n} -1$.
Let  $a_1,a_2,\dots,a_{n_0}$ be positive integers.
Suppose $K$ is an $(a_1,a_2,\ldots,a_{n_0})$-admissible path partition of $G(a_1,a_2,\ldots,a_{n_0})$.
Then, $\mathcal{K}_{m,n}$ has a $C_4$-face-magic Klein bottle labeling
$X= \{ x_{i,j} : (i,j)\in V(\mathcal{K}_{m,n})\}$ on $\mathcal{K}_{m,n}$ such that $\mathcal{L}(X)=K$.
In addition, we have the following two possibilities:
\begin{itemize}
\item  Suppose $(a_1,a_2, \ldots, a_{n_0})$ is not a palindrome. Then, there are
$(m-1)!$ distinct Klein bottle nonequivalent
$C_4$-face-magic Klein bottle labelings $X= \{ x_{i,j} : (i,j)\in V(\mathcal{K}_{m,n})\}$
on $\mathcal{K}_{m,n}$ such that $\mathcal{L}(X)=K$.
\item Suppose $(a_1,a_2, \ldots, a_{n_0})$ is a palindrome. Then, there are
$2^{m-1} \, (m-1)!$ distinct Klein bottle nonequivalent
$C_4$-face-magic Klein bottle labelings $X= \{ x_{i,j} : (i,j)\in V(\mathcal{K}_{m,n})\}$
on $\mathcal{K}_{m,n}$ such that $\mathcal{L}(X)=K$.
\end{itemize}
\end{lemma}

\begin{proof}
Suppose $(a_1,a_2, \ldots, a_{n_0})$ is not a palindrome.
Let $T$ be the set of end vertex labels of paths in $K$ incident to the edge labeled $a_1$. In other words, we let
\begin{align*}
T = &\biggl\{ k\in \mathbb{Z} : 1\leq k \leq m n, \{k, mn+1-k\} \mbox{ is an end vertex of a path } P \mbox{ in } K, \\
    &\{\ell, mn+1-\ell\} \mbox{ is adjacent to }  \{k, mn+1-k\} \mbox{ in } P, k+\ell =a_1, \mbox{ and the } \\
    &\mbox{ sequence of labels on edges in the path } P \mbox{ from vertex } \{k, mn+1-k\} \mbox{ to } \\
    &\mbox{ the other end vertex of } P \mbox{ is } (a_1,a_2,\ldots,a_{n_0}). \biggr\}.
\end{align*}
Since there are $m$ paths in $K$ and the sequence $(a_1,a_2, \ldots, a_{n_0})$ is not a palindrome,
there are exactly $m$ positive integers in the set $T$.
For convenience, we label the elements of $T$ in ascending order.
Let
\begin{equation} \label{eqnStandardVertexLabelsNonPalindrome}
T=\{ k_1,k_2,\ldots,k_m \} \mbox{ \ \ \  such that }k_i < k_{i+1}, \mbox{ for } 1\le i <m.
\end{equation}
Let $\alpha$ be a permutation on the set $\{1,2,3,\ldots,m\}$.
We define a labeling $L(K,\alpha) = \{ x_{i,j} : (i,j)\in V(\mathcal{K}_{m,n})\}$ on $\mathcal{K}_{m,n}$ in the following manner:
\begin{itemize}
\item For $1\le i \le \stkout{m}^{+}$ and $1\le j \le \stkout{n}$, let
\begin{equation} \label{eqnStandardlabelingNoPalindrome1}
x_{2i-1,j} = (-1)^{j+1} k_{\alpha(i)}  + \sum_{\ell=1}^{j-1} (-1)^{\ell+1} a_{j-\ell}.
\end{equation}

\item For $1\le i \le \stkout{m}^{+}$ and $1\le j \le \stkout{n}$, let
\begin{equation} \label{eqnStandardlabelingNoPalindrome2}
x_{2i-1,n+1-j} = \tfrac{1}{2} S- x_{2i-1,j}.
\end{equation}

\item For $\stkout{m}^{+} +1\le i \le m$ and $1\le j \le \stkout{n}$, let
\begin{equation} \label{eqnStandardlabelingNoPalindrome3}
x_{2(i\, - \, \stkout{m}^{+}),n+1-j} = (-1)^{j+1}  k_{\alpha(i)}  + \sum_{\ell=1}^{j-1} (-1)^{\ell+1} a_{j-\ell}.
\end{equation}

\item For $\stkout{m}^{+} +1\le i \le m$ and $1\le j \le \stkout{n}$, let
\begin{equation} \label{eqnStandardlabelingNoPalindrome4}
x_{2(i \, - \, \stkout{m}^{+}),j} = \tfrac{1}{2} S- x_{2(i \, - \, \stkout{m}^{+}),n+1-j}.
\end{equation}
\end{itemize}

We show that for $1 \le i \le \, \stkout{m}^{+}$,
the sequence $( x_{2i-1,j} : 1 \le j \le \, \stkout{n} )$ is the vertex label sequence of
some $(a_1,a_2,\ldots,a_{n_0})$-admissible path in $K$.
Let $1\le i \le \, \stkout{m}^{+}$.
Then, $x_{2i-1,1} = k_{\alpha(i)}$.
For all $1\le j \le n_0$, we have $x_{2i-1,j} + x_{2i-1,j+1} = a_j$.
Hence for all $1\le i \le \, \stkout{m}^{+}$, the sequence of vertices
$( \{x_{2i-1,j}, \tfrac{1}{2} S - x_{2i-1,j} \} : 1\le j \le \, \stkout{n} )$ is an
$(a_1,a_2,\ldots,a_{n_0})$-admissible path in $K$ with vertex label sequence $( x_{2i-1,j} : 1\le j \le \, \stkout{n} )$.

Similarly, we show that for $\stkout{m}^{+} \, + 1 \le i \le  m$,
the sequence $( x_{2(i- \, \stkout{m}),n+1-j} : 1 \le j \le \, \stkout{n} )$ is the vertex label sequence of
some $(a_1,a_2,\ldots,a_{n_0})$-admissible path in $K$.
Let $\stkout{m}^{+} \, + 1 \le i \le  m$.
Then, $x_{2(i- \, \stkout{m}^{+}),n} = k_{\alpha(i)}$.
For all $1\le j \le n_0$, we have $x_{2(i- \, \stkout{m}^{+}),n+1-j} + x_{2(i- \, \stkout{m}^{+}),n-j} = a_j$.
Hence, for all $\stkout{m}^{+} \, + 1 \le i \le  m$, the sequence of vertices
$( \{x_{2(i- \, \stkout{m}^{+}),n+1-j}, \tfrac{1}{2} S - x_{2(i- \, \stkout{m}^{+}),n+1-j} \} : 1\le j \le \, \stkout{n} )$ is an
$(a_1,a_2,\ldots,a_{n_0})$-admissible path in $K$ with vertex label sequence
$( x_{2(i- \, \stkout{m}^{+}),n+1-j} : 1\le j \le \, \stkout{n} )$ .

Next we show that the labeling $L(K,\alpha)$ uses each integer in the set $\{ k\in \mathbb{Z}: 1\le k\le mn\}$ exactly once.
For $1\le i \le \stkout{m}^{+}$, we have $x_{2i-1,1} = k_{\alpha(i)}$, and for $\stkout{m}^{+} +1\le i \le m$,
we have $x_{2(i-\, \stkout{m}^{+}),n} =k_{\alpha(i)}$. Thus, each of the values $x_{2i-1,1}$, $1\le i \le \stkout{m}^{+}$,
and $x_{2(i-\, \stkout{m}^{+}),n}$ for $\stkout{m}^{+} +1\le i \le m$,
are distinct end vertex labels incident to an edge labeled $a_1$
of an $(a_1,a_2,\ldots,a_{n_0})$-admissible path in $K$.
For $1\le i \le \stkout{m}^{+}$, the sequence $(x_{2i-1,j} : 1\le j \le n)$ is the
sequence of labels $(x_{2i-1,j} : 1\le j \le \stkout{n})$  followed by the sequence of labels
$( \tfrac{1}{2} S - x_{2i-1,\stkout{n} \, +1-j} : 1\le j \le \stkout{n})$.
Furthermore, each label in the sequence $(x_{2i-1,j} : 1\le j \le n)$ appears exactly once as a label
for a vertex in the
path $P_i =( \{ x_{2i-1,j}, \tfrac{1}{2} S -  x_{2i-1,j}: 1\le j \le \stkout{n})$ in $K$.
Similarly, for $1 \le i \le \stkout{m}$, the sequence $(x_{2i,n+1-j} : 1\le j \le n)$ is the
sequence of labels $(x_{2i,n+1-j} : 1\le j \le \stkout{n})$ followed by the sequence of labels
$( \tfrac{1}{2} S - x_{2i,\stkout{n} \, +j} : 1\le j \le \stkout{n})$.
Furthermore, each label in the sequence $(x_{2i,n+1-j} : 1\le j \le n)$ appears exactly once as a label
for a vertex in the path
$P_{i+\, \stkout{m}^{+}} =( \{ x_{2i,n+1-j}, \tfrac{1}{2} S - x_{2i,n+1-j} \} : 1\le j \le \stkout{n})$ in $K$.
Since the set $\{ V(P_i) : 1 \le i \le m\}$ is a partition on the vertices of $G(a_1,a_2,\ldots,a_{n_0})$,
each integer $k$ in the set $\{ k\in \mathbb{Z}: 1\le k \le mn\}$ is used exactly once.
By Lemma \ref{lemmaSequenceOfPairSums}, $L(K,\alpha)$ is a $C_4$-face-magic Klein bottle labeling on $\mathcal{K}_{m,n}$.

We consider the graph automorphisms of $\mathcal{K}_{m,n}$ that are induced by a symmetry on the Klein bottle.
From Notation \ref{notnKleinBottleSymmetries}, these graph automorphisms of $\mathcal{K}_{m,n}$ are in the group $KBLS(m,n)$.
We consider the result of each graph automorphism in $KBLS(m,n)$ on $L(K,\alpha)$.
We let $i$ be an integer and apply the graph automorphism $U^{2i+1}$ from $KBLS(m,n)$ on $L(K,\alpha)$.
Observe that $\mathcal{L}( U^{2i+1} ( L(K,\alpha )))$ is an $(S-a_1,S-a_2,\ldots,S-a_{n_0 -1},S-a_{n_0})$-admissible
path partition of $G(S-a_1,S-a_2,\ldots,S-a_{n_0 -1},S-a_{n_0})$.
Thus, $\mathcal{L}( U^{2i+1} ( L(K,\alpha ))) \ne K$.
Next, we apply the graph automorphism $HU^{2i}$ to $L(K,\alpha)$.
Then, $\mathcal{L}( HU^{2i} ( L(K,\alpha )))$ is an $(S-a_{n_0}, S-a_{n_0 -1}, \ldots, S-a_2, S-a_1)$-admissible
path partition of $G(S-a_{n_0}, S-a_{n_0 -1}, \ldots, S-a_2, S-a_1)$.
Thus, $\mathcal{L}( HU^{2i} ( L(K,\alpha ))) \ne K$.
Similarly, we apply the graph automorphism $HU^{2i+1}$ to $L(K,\alpha)$.
Then, $\mathcal{L}( HU^{2i+1} ( L(K,\alpha )))$ is an $(a_{n_0}, a_{n_0 -1}, \ldots, a_2, a_1)$-admissible
path partition of $G(a_{n_0}, a_{n_0 -1}, \ldots, a_2, a_1 )$.
Thus, $\mathcal{L}( HU^{2i+1} ( L(K,\alpha ))) \ne K$.

On the other hand, when we apply the graph automorphism $U^{2i}$ to $L(K,\gamma)$,
where $\gamma$ is a permutation on the set $\{ 1,2,\ldots,n\}$,
we obtain $\mathcal{L}( U^{2i} ( L(K,\gamma ))) = K$.
Let $\beta$ be the cyclic permutation on the set $\{1,2,\ldots,m\}$ defined by $\beta(i) = i-1\pmod{m}$.
Let $X=\{ x_{i,j} : (i,j)\in V(\mathcal{K}_{m,n})$ be any $C_4$-face-magic Klein bottle labeling on $\mathcal{K}_{m,n}$
such that $\mathcal{L}(X)=K$.
Let $\gamma$ be the permutation on the set $\{1,2,\ldots,m\}$ such that
$x_{2i-1,1}= k_{\gamma(i)}$, for $1\le i\le \stkout{m}{+}$, and
$x_{2(i-\stkout{m}^{+}),n} = k_{\gamma(i)}$, for $\stkout{m}^{+} + 1 \le i \le m$.
Then $L(K,\gamma)=X$.
Let $t$ be the unique integer such that $0 \le t <m$ and
$\gamma(t) =1$. Let $\alpha =\beta^{t-1} \gamma$. Then, $\alpha(1)=1$.
Thus, there is a unique permutation $\alpha$ on the set $\{1,2,\ldots,m\}$ with $\alpha(1)=1$ and
$L(K,\alpha)=U^{-2(t-1)}(X)$.
Hence, for every $C_4$-face-magic labeling $X$ such that $\mathcal{L}(X)=K$,
there is a unique permutation $\alpha$ on the set $\{1,2,\ldots,m\}$  with $\alpha(1)=1$
such that $X$ is Klein bottle equivalent to $L(K,\alpha)$.
Furthermore, each distinct permutation $\alpha$ on $\{1,2,\ldots,m\}$ with
$\alpha(1)=1$ results in a distinct Klein bottle nonequivalent $C_4$-face-magic
Klein bottle labeling $L(K,\alpha)$ on $\mathcal{K}_{m,n}$ such that $\mathcal{L} (L(K,\alpha))=K$.
Since there are $(m-1)!$ permutation $\alpha$ on the set $\{1,2,\ldots,m\}$ with $\alpha(1)=1$,
there are $(m-1)!$ Klein bottle nonequivalent $C_4$-face-magic
Klein bottle labelings $X$ on $\mathcal{K}_{m,n}$ such that $\mathcal{L}(X)=K$.

Suppose the sequence $(a_1,a_2,\ldots,a_{n_0})$ is a palindrome.
We define the two-element sets of integers $\{ k,\ell\}$ such that
$k$ and $\ell$ are end vertex labels of an $(a_1,a_2,\ldots,a_{n_0})$-admissible path $P$ in $K$ in the following way:
Let
\begin{align*}
\widetilde{T} &= \biggl\{ \{ k,\ell \}  : k, \ell \in \mathbb{Z}, 1\leq k,\ell \leq m n, \{k, mn+1-k\}
\mbox{ and }   \{\ell, mn+1-\ell \} \mbox{ are end  } \\
  &\mbox{ vertices of an $(a_1,a_2,\ldots,a_{n_0})$-admissible path } P \mbox{ in } K,
    \mbox{ and there exist vertices }\\
    & \{p, mn+1-p\}  \mbox{ and } \{q, mn+1-q\}
    \mbox{ on path } P \mbox{ such that } k+p = a_1 = \ell +q \biggr\}.
\end{align*}
Since there are $m$ paths in $K$ and the sequence $(a_1,a_2, \ldots, a_{n_0})$ is a palindrome,
there are exactly $m$ two-element sets in $\widetilde{T}$.
For convenience, we label the elements of $\widetilde{T}$ in the following way:
Let
\begin{align}  \label{eqnStandardVertexLabelsPalindrome}
\widetilde{T} =\{ \{ k_{i,1}, k_{i,2} \} : 1 \le i \le m \}
&\mbox{ such that } k_{i,1} < k_{i,2}, \mbox{ for all } 1\le i \le m, \\
 &\mbox{ and } k_{i,1} < k_{i+1,1}, \mbox{ for all } 1\le i <m. \nonumber
\end{align}
Let $\alpha$ be a permutation on the set $\{1,2,3,\ldots,m\}$.
Let $\rho$ be a function from $\{ 1,2,\ldots,m\}$ to $\{1,2\}$.
We define a labeling $L(K,\alpha,\rho) = \{ x_{i,j} : (i,j)\in V(\mathcal{K}_{m,n})\}$ on $\mathcal{K}_{m,n}$ in the following manner:
\begin{itemize}
\item For $1\le i \le \stkout{m}^{+}$ and $1\le j \le \stkout{n}$, let
\begin{equation} \label{eqnStandardlabelingPalindrome1}
x_{2i-1,j} = (-1)^{j+1} k_{\alpha(i), \rho(i)}  + \sum_{\ell=1}^{j-1} (-1)^{\ell+1} a_{j-\ell}.
\end{equation}

\item For $1\le i \le \stkout{m}^{+}$ and $1\le j \le \stkout{n}$, let
\begin{equation} \label{eqnStandardlabelingPalindrome2}
x_{2i-1,n+1-j} = \tfrac{1}{2} S- x_{2i-1,j}.
\end{equation}

\item For $\stkout{m}^{+} +1\le i \le m$ and $1\le j \le \stkout{n}$, let
\begin{equation} \label{eqnStandardlabelingPalindrome3}
x_{2(i \, - \, \stkout{m}^{+}),n+1-j} = (-1)^{j+1}  k_{\alpha(i),\rho(i)}  + \sum_{\ell=1}^{j-1} (-1)^{\ell+1} a_{j-\ell}.
\end{equation}

\item For $\stkout{m}^{+} +1\le i \le m$ and $1\le j \le \stkout{n}$, let
\begin{equation} \label{eqnStandardlabelingPalindrome4}
x_{2(i \, - \, \stkout{m}^{+}),j} = \tfrac{1}{2} S- x_{2(i \, - \, \stkout{m}^{+}),n+1-j}.
\end{equation}
\end{itemize}
An argument similar to that in the case when $(a_1,a_2,\ldots,a_{n_0})$ is not a palindrome
establishes each of the following statements.
We leave the details of the proof to the reader.
Given a $C_4$-face-magic labeling $X= \{ x_{i,j} : (i,j) \in V(\mathcal{K}_{m,4})\}$ such that $\mathcal{L}(X)=K$,
we can find a unique graph automorphism $A\in KBLS(m,4)$, a unique permutation $\alpha$ on $\{1,2,\dots,m\}$ such that $\alpha(1)=1$,
and a unique function $\rho:\{1,2,\ldots,m\}\rightarrow \{1,2\}$ such that $\rho(1)=1$ for which $A(X)= L(K,\alpha,\rho)$.
Furthermore, each distinct choice of a permutation
$\alpha$ on $\{1,2,\dots,m\}$ such that $\alpha(1)=1$,
and each distinct choice of a function $\rho:\{1,2,\ldots,m\}\rightarrow \{1,2\}$
such that $\rho(1)=1$ yields a Klein bottle nonequivalent
$C_4$-face-magic Klein bottle labeling $L(K,\alpha,\rho)$ on $\mathcal{K}_{m,4}$ such that $\mathcal{L}(L(K,\alpha,\rho))=K$.
Thus, there $2^{m-1} \, (m-1)!$ Klein bottle nonequivalent $C_4$-face-magic
Klein bottle labelings $X$ on $\mathcal{K}_{m,n}$ such that $\mathcal{L}(X)=K$.
\end{proof}

\begin{definition} \label{defnStandardKleinBottleLabeling}
Let $m\geq 3$ be an odd integer, $n \geq 4$ be an even integer, $n_1 = \stkout{n}=n/2$, and $n_0 =\stkout{n} -1$.
Let  $a_1,a_2,\dots,a_{n_0}$ be positive integers.
\begin{itemize}
\item  Suppose $(a_1,a_2, \ldots, a_{n_0})$ is not a palindrome.
Let $K$ be the $(a_1,a_2, \ldots, a_{n_0})$-admis\-si\-ble path partition of $G(a_1,a_2, \ldots, a_{n_0})$ such that for each
$(a_1,a_2,\ldots,$ $a_{n_0})$-admissible path $P$ of $K$, there exists $k_i \in T$ such that
$\{ k_i, mn+1-k_i \}$ is an end vertex of $P$, where $k_i$ is defined in equation (\ref{eqnStandardVertexLabelsNonPalindrome}).
Let $\epsilon$ be the identity function on $\{1,2,\ldots,m\}$.
The labeling $L(K,\epsilon)$ as defined by equations (\ref{eqnStandardlabelingNoPalindrome1})
through (\ref{eqnStandardlabelingNoPalindrome4}) in the proof of Lemma \ref{lemmaSpanningGraphOfPathsCriteria}
is called the \textit{standard $C_4$-face-magic Klein bottle labeling} on $\mathcal{K}_{m,n}$
associated with $K$.

\item  Suppose $(a_1,a_2, \ldots, a_{n_0})$ is a palindrome.
Let $K$ be the $(a_1,a_2, \ldots, a_{n_0})$-admissible path partition of $G(a_1,a_2, \ldots, a_{n_0})$ such that for each
$(a_1,a_2,\ldots,$ $a_{n_0})$-admissible path $P$ of $K$, there exists $\{ k_{i,1}, k_{i,2} \} \in \widetilde{T}$ such that
$\{ k_{i,1}, mn+1-k_{i,1} \}$ is an end vertex of $P$, where $\{ k_{i,1}, k_{i,2} \}$
is defined in equation (\ref{eqnStandardVertexLabelsPalindrome}).
Let $\epsilon$ be the identity function on $\{1,2,\ldots,m\}$ and let $\mu:\{1,2,\ldots,m\}\rightarrow \{1,2\}$
be defined by $\mu(i)=1$, for all $1\le i\le m$.
The labeling $L(K,\epsilon,\mu)$ as defined by equations (\ref{eqnStandardlabelingPalindrome1})
through (\ref{eqnStandardlabelingPalindrome4}) in the proof of Lemma \ref{lemmaSpanningGraphOfPathsCriteria}
is called the
\textit{standard $C_4$-face-magic Klein bottle labeling} on $\mathcal{K}_{m,n}$
associated with $K$.
\end{itemize}
\end{definition}

\begin{definition} \label{defnAPPSet}
Let $m\ge 3$ be an odd integer, $n\ge 4$ be an even integer, and $n_0 = \, \stkout{n} \, -1$.
\begin{enumerate}
\item We let $NAPP(m,n)$ be the set of all non-palindromic sequences $(a_1,a_2,\ldots, \allowbreak a_{n_0})$ such that
there exists an  $(a_1,a_2,\ldots,a_{n_0})$-admissible path partition of $G(a_1,a_2,\ldots,\allowbreak a_{n_0})$.
\item We let $PAPP(m,n)$ be the set of all palindromic sequences $(a_1,a_2,\ldots,a_{n_0})$ such that
there exists an  $(a_1,a_2,\ldots,a_{n_0})$-admissible path partition of $G(a_1,$ $a_2,\ldots,a_{n_0})$.
\end{enumerate}

Given two non-palindromic sequences $(a_1,a_2,\dots,a_{n_0})$ and $(a'_1,a'_2,\dots,a'_{n_0})$ in $NAPP\allowbreak (m,\allowbreak n)$,
we say that $(a_1,a_2,\dots,a_{n_0})$ is \textit{path partition equivalent to}  $(a'_1,a'_2,\dots,\allowbreak a'_{n_0})$
if there exist an $(a_1,a_2,\dots,a_{n_0})$-admissible path partition $K$ of $G(a_1,a_2,\dots,\allowbreak a_{n_0})$,
an $(a'_1,a'_2,\dots,\allowbreak a'_{n_0})$-admissible path partition $K'$ of $G(a'_1,a'_2,\dots,\allowbreak a'_{n_0})$,
$A\in KBLS(m,n)$,
and $\alpha \in S_n$ for which $\alpha(1)=1$
such that $\mathcal{L}(A(L(K,\alpha))) = K'$.
Here, $S_n$ is the symmetry group of degree $n$.

Given two palindromic sequences $(a_1,a_2,\dots,a_{n_0})$ and $(a'_1,a'_2,\dots,a'_{n_0})$ in

\noindent $PAPP(m,n)$,
we say that $(a_1,a_2,\dots,a_{n_0})$ is \textit{path partition equivalent to}  $(a'_1,a'_2,\dots,$ $a'_{n_0})$
if there exist an $(a_1,a_2,\dots,a_{n_0})$-admissible path partition $K$ of $G(a_1,a_2,\dots,$ $a_{n_0})$,
an $(a'_1,a'_2,\dots,a'_{n_0})$-admissible path partition $K'$ of $G(a'_1,a'_2,\dots,a'_{n_0})$,
$A\in KBLS(m,n)$,
$\alpha \in S_n$ for which $\alpha(1)=1$,
and a function $\rho :\{1,2,\ldots,n\} \rightarrow \{1,2\}$ for which $\rho(1)=1$
such that $\mathcal{L}(A(L(K,\alpha,\rho))) = K'$.
\end{definition}

\begin{remark}
The path partition equivalence relation is an equivalence relation on both $NAPP(m,n)$ and $PAPP(m,n)$.
\end{remark}

\begin{lemma} \label{lemmaPathPartitionEquivalence}
Let $m\ge 3$ be an odd integer, $n\ge 4$ be an even integer, and $n_0 = \, \stkout{n} \, -1$.
\begin{enumerate}
\item Let $(a_1,a_2,\ldots,a_{n_0})$ and $(a'_1,a'_2,\dots,a'_{n_0})$
be non-palindromic sequences in $NAPP\allowbreak (m,\allowbreak n)$.
Then, $(a_1,a_2,\dots,a_{n_0})$ is path partition equivalent to $(a'_1,a'_2,$ $\dots,a'_{n_0})$
if and only if either $(a'_1,a'_2,\dots,a'_{n_0})= (a_1,a_2,\dots,a_{n_0})$,
$(a'_1,a'_2,$ $\dots,a'_{n_0})= (a_{n_0},a_{n_0 -1},\dots,a_{1})$,
$(a'_1,a'_2,\dots,a'_{n_0})= (S-a_1,S-a_2,\dots,S-a_{n_0})$, or
$((a'_1,a'_2,\dots,a'_{n_0})= (S-a_{n_0},S-a_{n_0 -1},\dots,S-a_{1})$.

\item Let $(a_1,a_2,\ldots,a_{n_0})$ and $(a'_1,a'_2,\dots,a'_{n_0})$
be palindromic sequences in $PAPP(m,n)$.
Then, $(a_1,a_2,  \allowbreak  \dots,  \allowbreak a_{n_0})$ is path partition equivalent to $(a'_1,a'_2,$ $\dots, \allowbreak a'_{n_0})$
if and only if either $(a'_1,a'_2,\dots,a'_{n_0})= (a_1,a_2,\dots,  \allowbreak a_{n_0})$ or
$(a'_1,a'_2,$ $\dots,a'_{n_0})= (S-a_1,S-a_2,\dots,S-a_{n_0})$.
\end{enumerate}
\end{lemma}

\begin{proof}
We prove the lemma only for the case when $(a_1,a_2,\ldots,a_{n_0})$ and $(a'_1,a'_2,\dots,$ $a'_{n_0})$
are non-palindromic sequences in $NAPP(m,n)$. The case when
$(a_1,a_2,\ldots,a_{n_0})$ and $(a'_1,a'_2,\dots,a'_{n_0})$
are palindromic sequences in $PAPP(m,n)$ is similar, and we leave the details of the proof for that case to the reader.

($\Rightarrow$) Let $(a_1,a_2,\ldots,a_{n_0})$ and $(a'_1,a'_2,\dots,a'_{n_0})$
be non-palindromic sequences in $NAPP\allowbreak (m,\allowbreak n)$ such that
$(a_1,a_2,\dots,a_{n_0})$ is path partition equivalent to $(a'_1,a'_2,\dots,$ $a'_{n_0})$.
Then, there exist an $(a_1,a_2,\dots,a_{n_0})$-admissible path partition $K$ of $G(a_1,a_2,$ $\dots,$ $a_{n_0})$,
an $(a'_1,a'_2,\dots,a'_{n_0})$-admissible path partition $K'$ of $G(a'_1,a'_2,\dots,a'_{n_0})$,
$A\in KBLS(m,n)$,
and $\alpha \in S_n$ for which $\alpha(1)=1$
such that $\mathcal{L}(A(L(K,\alpha))) = K'$.
Since $A\in KBLS(m,n)$, we have for some integer $i$, either $A=U^{2i}$, $A=U^{2i+1}$, $A=HU^{2i}$, $A=HU^{2i+1}$.

First, suppose $A=U^{2i}$ for some integer $i$.
Then, $\mathcal{L}(A(L(K,\alpha)))$ is an $(a_1,a_2,\dots,$ $a_{n_0})$-admissible path partition of $G(a_1,a_2,\dots,a_{n_0})$.
However, $K' = \mathcal{L}(A(L(K,\alpha)))$ is
an $(a'_1,a'_2,\dots,a'_{n_0})$-admissible path partition of $G(a'_1,a'_2,\dots,$ $a'_{n_0})$.
Thus, $(a'_1,a'_2,$ $\dots,a'_{n_0})= (a_1,a_2,\dots,a_{n_0})$.

Second, suppose $A=U^{2i+1}$ for some integer $i$.
Then, $\mathcal{L}(A(L(K,\alpha)))$ is an $(S-a_1,S-a_2,\dots,S-a_{n_0})$-admissible
path partition of $G(S-a_1,S-a_2,\dots,S-a_{n_0})$.
However, $K' = \mathcal{L}(A(L(K,\alpha)))$
is an $(a'_1,a'_2,\dots,a'_{n_0})$-admissible path partition of $G(a'_1,a'_2,\dots,$ $a'_{n_0})$.
Thus, $(a'_1,a'_2,\dots,a'_{n_0})= (S-a_1,S-a_2,\dots,S-a_{n_0})$.

Third, suppose $A=HU^{2i}$ for some integer $i$.
Then, $\mathcal{L}(A(L(K,\alpha)))$ is an $(S-a_{n_0},S-a_{n_0 -1},\dots,S-a_{1})$-admissible
path partition of $G(S-a_{n_0},S-a_{n_0 -1},\dots,S-a_{1})$.
However, $K' = \mathcal{L}(A(L(K,\alpha)))$
is an $(a'_1,a'_2,\dots,a'_{n_0})$-admissible path partition of $G(a'_1,a'_2,\dots,$ $a'_{n_0})$.
Thus, $(a'_1,a'_2,\dots,a'_{n_0})= (S-a_{n_0},S-a_{n_0 -1},\dots,S-a_{1})$.

Fourth, suppose $A=HU^{2i+1}$ for some integer $i$.
Then, $\mathcal{L}(A(L(K,\alpha)))$ is an $(a_{n_0},a_{n_0 -1},\allowbreak \dots,a_1)$-admissible
path partition of $G(a_{n_0},a_{n_0 -1},\dots,a_1)$.
However, $K' = \mathcal{L}(A(L(K,\alpha)))$
is an $(a'_1,a'_2,\dots,a'_{n_0})$-admissible path partition of $G(a'_1,$ $a'_2,\dots,$ $a'_{n_0})$.
Thus, $(a'_1,a'_2,\dots,a'_{n_0})\allowbreak = (a_{n_0},a_{n_0 -1},\dots,a_1)$.

($\Leftarrow$) Suppose either $(a'_1,a'_2,\dots,a'_{n_0})= (a_1,a_2,\dots,a_{n_0})$,
$(a'_1,a'_2,\dots,a'_{n_0})= (a_{n_0},$ $a_{n_0 -1},\allowbreak \dots,a_{n1})$,
$(a'_1,a'_2,\allowbreak \dots,\allowbreak a'_{n_0})= (S-a_1,S-a_2,\dots,S-a_{n_0})$, or
$(a'_1,a'_2,\dots,a'_{n_0})$ $= (S-a_{n_0},S-a_{n_0 -1},\dots,S-a_{1})$.
We want to show that $(a_1,a_2,\dots,a_{n_0})$ is path partition equivalent to $(a'_1,a'_2,\dots,a'_{n_0})$.
Since $(a_1,a_2,\dots,a_{n_0})$ is in $NAPP(m,n)$, there is an $(a_1,a_2,\dots,a_{n_0})$-admissible path partition $K$
of $G(a_1,a_2,\dots,a_{n_0})$.

First, suppose $(a'_1,a'_2,\dots,a'_{n_0})= (a_1,a_2,\dots,a_{n_0})$. Let $A=U^0 \in KBLS(m,n)$.
Then, $K$ is an $(a'_1,a'_2,\dots,\allowbreak a'_{n_0})$-admissible path partition of $G(a'_1,a'_2,\dots,a'_{n_0})$
such that $K = \mathcal{L}(A(L(K,\epsilon)))$.
Thus, $(a_1,a_2,\dots,a_{n_0})$ is path partition equivalent to $(a'_1,a'_2,\dots,a'_{n_0})$.

Second, suppose $(a'_1,a'_2,\dots,a'_{n_0})= (a_{n_0},a_{n_0 -1},\dots,a_1)$. Let $A=HU \in KBLS(m,n)$
and let $K' = \mathcal{L}(A(L(K,\epsilon)))$.
Then, $K'$ is an $(a'_1,a'_2,\dots,a'_{n_0})$-admissible path partition of $G(a'_1,a'_2,\dots,a'_{n_0})$
such that $K' = \mathcal{L}(A(L(K,\epsilon)))$.
Thus, $(a_1,a_2,$ $\dots,a_{n_0})$ is path partition equivalent to $(a'_1,a'_2,\dots,a'_{n_0})$.

The proofs of the remaining two cases when $(a'_1,a'_2,\dots,a'_{n_0})= (S-a_1,S-a_2,\dots,S-a_{n_0})$ or
$(a'_1,a'_2,\dots,\allowbreak a'_{n_0})= (S-a_{n_0},S-a_{n_0 -1},\dots,S-a_{1})$ are similar to the proof of the case
when $(a'_1,a'_2,\dots,a'_{n_0})= (a_{n_0},a_{n_0 -1},\dots,a_1)$.
If $(a'_1,a'_2,\dots,a'_{n_0})= (S-a_1,S-a_2,\dots,S-a_{n_0})$, we let $A=U \in KBLS(m,n)$.
If $(a'_1,a'_2,\dots,a'_{n_0})= (S-a_{n_0},S-a_{n_0 -1},\dots,S-a_{1})$, we let $A=H \in KBLS(m,n)$.
\end{proof}

\section{Results on an $m \times 4$ Klein bottle grid graph} \label{S:resultsonmby4gridgraph}

\begin{corollary} \label{corNumberOfLabelinsOnKM4}
Let $m\ge 3$ be an odd positive integer.
Let $a$ be a positive integer.
Suppose $G(a)$ has a perfect matching $K$.
Then, $\mathcal{K}_{m,4}$ has $2^{m-1} \, (m-1)!$ distinct Klein bottle nonequivalent
$C_4$-face-magic Klein bottle labelings $X= \{ x_{i,j} : (i,j)\in V(\mathcal{K}_{m,4})\}$ on $\mathcal{K}_{m,4}$
such that $\mathcal{L}(X)=K$.
\end{corollary}

\begin{proof}
Let $K$ be a perfect matching on $G(a)$.
Then, $K$ is an $(a)$-admissible path partition of $G(a)$.
Since the sequence $(a)$ is a palindrome,
by Lemma \ref{lemmaSpanningGraphOfPathsCriteria},
there are $2^{m-1} \, (m-1)!$ Klein bottle nonequivalent $C_4$-face-magic Klein bottle labelings $X$ on $\mathcal{K}_{m,4}$
such that $\mathcal{L}(X)=K$.
\end{proof}

\begin{notation} \label{notnDiscreteClosedInterval}
Let $a$ and $b$ be integers such that $a\le b$. The {\it discrete closed interval
from $a$ to $b$} is the set denoted by $[a,b] =\{ i\in \mathbb{Z} : a\le i \le b \}$.
\end{notation}

\begin{theorem} \label{thmNumberOfC4FaceMagicLabelsforK4n}
Let $m$ be an odd positive integer such that $m\geq 3$.
Then, the number of Klein bottle nonequivalent $C_4$-face-magic Klein bottle labelings on $\mathcal{K}_{m,4}$
is $2^m \, (m-1)! \, \tau(m)$, where $\tau(m)$ is the number of positive divisors of $m$.
\end{theorem}

\begin{proof}
Let $X = \{ x_{i,j} : (i,j)\in V(\mathcal{K}_{m,n}) \}$ be a $C_4$-face-magic Klein bottle labeling on $\mathcal{K}_{m,4}$
and let $a= x_{1,1} + x_{1,2}$.
Then, $\mathcal{L}(X)$ is an $(a)$-admissible path partition of $G(a)$.
Since $a=z_1 + z_2$ for integers $z_1$
and $z_2$ such that $1\le z_1,z_2 \le 4m$, we have $2\le a \le 8m$.
For convenience, let $S$ be the $C_4$-face-magic value of $X$.
By Lemma \ref{lemmaKleinCheckerboardC4FaceValue}, $S= 2(4m+1)$.
If $a \ge \tfrac{1}{2} S =4m+1$, then $U(X)$ is a $C_4$-face-magic Klein bottle labeling on $\mathcal{K}_{m,4}$
that is Klein bottle equivalent to $X$.
Furthermore, $\mathcal{L}(U(X))$ is an $(S-a)$-admissible path partition of $G(S-a)$, where $S-a\le \tfrac{1}{2} S= 4m+1$.
By replacing $a$ with $S-a$, we may assume that $2\le a \le 4m+1$.

If $G(a)$  has a perfect matching, then vertex $\{ 2m, 2m+1\}$ is adjacent to some vertex $\{ i, 4m+1-i\}$.
Thus, there is an integer $z\in \{i,4m+1-i \}$ such that $a \ge 2m+z \ge 2m+1$.
Also, all pairs of integers $z_1$ and $z_2$ such that $z_1 + z_2 =4m+1$ are paired in the set $\{ z_1,z_2\}$. Thus, $a\ne 4m+1$.
Hence, $2m+1 \le a \le 4m$.

Let $\{ z_1,z_2\}\in V(G(a))$.
If there are integers $y_1$ and $y_2$ such that
$1 \le y_1,y_2\le 4m$, $z_1 + y_1 =a$ and $z_2 + y_2 = a$, then vertex $\{ z_1,z_2\}$ in $G(a)$ has degree two.
However, the only integers $y_1$ and $y_2$ that satisfy this criterion are $y_1 = a -z_1$ and $y_2 = a-z_2$.
Thus, the degree of vertex $\{ z_1,z_2\}$ is at most two.
Hence, $G(a)$ is a disjoint union of paths and cycles. For convenience, let $b= 4m+1$.
We show that there are no cycles in $G(a)$.
For the purpose of contradiction, suppose $G(a)$ has a cycle on $k$ vertices for some positive integer $k$.
Let $( \{ y_{i,1}, y_{i,2}\} : 1\le i \le k )$ be a walk in $G(a)$ such that
$y_{i,2} + y_{i+1,1} =a$, for $1\le i \le k-1$ and $y_{k,2} + y_{1,1} =a$.
For convenience, for $j=1,2$, let $y_{k+1,j}=y_{1,j}$.
Since $y_{i,1} + y_{i,2} =b$ for $1\le i \le k+1$, and $y_{i,2} + y_{i+1,1} =a$ for $1\le i \le k$,
we have $y_{i,1} =y_{1,1} + (i-1)(a-b)$, for $1\le i \le k+1$.
Then, $y_{1,1} = y_{k+1,1} = y_{1,1} +k(a-b)$ implies that $a=b$. However, $a<b$. This is a contradiction.
Therefore, $G(a)$ is a disjoint union of paths.
Since $G(a)$ is a disjoint union of paths, $G(a)$ has a perfect matching
if and only if each path in $G(a)$ is a path on an even number of vertices.
We next establish that for all integers $k$ such that $1\le k \le 2m$,
$G(4m+1-k)$ has a perfect matching if and only if
either $k=d$ or $k=2d$ for some positive divisor $d$ of $m$. We consider five cases to verify this statement.

\textbf{ Case 1.}
Let $k=d$ for some positive divisor $d$ of $m$.
Let $m = m_2 d$ for some positive integer $m_2$.
Since $d$ is odd, we can write $d=2d_1 +1$ for some nonnegative integer $d_1$.
Let $\ell = \tfrac{1}{2} (4m+1-d) = 2m-d_1$.
The sequence of integers $\bigl( di+\ell : i\in [-(2m_2 -1), 2m_2] \bigr)$ satisfies $ di+\ell \equiv d_1 +1 \pmod{2d_1 +1}$,
and the corresponding descending sequence of integers
$\bigl( 4m+1-  (di+\ell) : i\in [-(2m_2 -1), 2m_2] \bigr)$ satisfies $4m+1- (di+\ell) \equiv d_1 +1 \pmod{2d_1 +1}$.
Thus, the sequence $\bigl( 4m+1-  (di+\ell) : i\in [-(2m_2 -1), 2m_2] \bigr)$ is the
sequence $\bigl( di+\ell : i\in [-(2m_2 -1), 2m_2] \bigr)$ in reverse order.
Hence, the sequence of vertices $\bigl( \{  di+\ell, 4m+1-  (di+\ell) \} : i\in [-(2m_2 -1), 0] \bigr)$
is a path in $G(4m+1-d)$ on $2m_2 = (2m)/d$ vertices.

Let $i_0 \in [1,d_1]$. The sequence of integers
$\bigl( di+ i_0 : i\in [0, 4m_2 -1] \bigr)$ satisfies $ di+ i_0 \equiv 1,2,\ldots, d_1 \pmod{2d_1 +1}$
for $i_0 =1,2,\ldots,d_1$,
and the corresponding descending sequence of integers
$\bigl( 4m+1-  (di+ i_0 ) : i\in [0, 4m_2 -1] \bigr)$ satisfies $4m+1- (di+ i_0) \equiv 2d_1 +1, 2d_1,\ldots, d_1 +2 \pmod{2d_1 +1}$
for $i_0 =1,2,\ldots,d_1$.
Thus, the sequence of vertices $\bigl( \{  di+ i_0, 4m+1-  (di+ i_0 ) \} : i\in [0, 4m_2 -1] \bigr)$
is a path in $G(4m+1-d)$ on $4m_2 = (4m)/d$ vertices, for all $i_0 \in [1,d_1]$.
Hence, $G(4m+1-d)$ is a disjoint union of $d_1$ paths
where each path has $(4m)/d$ vertices and one path on $(2m)/d$ vertices.
Since $G(4m+1-d)$ is a disjoint union of paths where each path has an even number of vertices,
$G(4m+1-d)$ has a unique perfect matching.

\textbf{Case 2.}
Let $k=2d$ for some positive divisor $d$ of $m$.
Let $m = m_2 d$ for some positive integer $m_2$.
Let $i_0 \in [1,d]$. The sequence of integers
$\bigl( (2d)i+ i_0 : i\in [0, 2m_2 -1] \bigr)$ satisfies $ (2d)i+ i_0 \equiv 1,2,\ldots, d \pmod{2d}$
for $i_0 =1,2,\ldots,d$,
and the corresponding descending sequence of integers
$\bigl( 4m+1-  ((2d)i+ i_0 ) : i\in [0, 2m_2 -1] \bigr)$ satisfies $4m+1- ((2d)i+ i_0) \equiv 2d, 2d -1,\ldots, d +1 \pmod{2d}$
for $i_0 =1,2,\ldots,d$.
Thus, the sequence of vertices $\bigl( \{  (2d)i+ i_0, 4m+1-  ((2d)i+ i_0 ) \} : i\in [0, 2m_2 -1] \bigr)$
is a path in $G(4m+1-2d)$ on $2m_2 = (2m)/d$ vertices, for all $i_0 \in [1,d]$.
Hence, $G(4m+1-2d)$ is a disjoint union of $d$ paths where each path has $(2m)/d$ vertices.
Since $G(4m+1-2d)$ is a disjoint union of paths where each path has an even number of vertices,
$G(4m+1-2d)$ has a unique perfect matching.

\textbf{Case 3.}
Let $k=4d$ for some positive divisor $d$ of $m$.
Let $m = m_2 d$ for some positive integer $m_2$.
Let $i_0 \in [1,2d]$. The sequence of integers
$\bigl( (4d)i+ i_0 : i\in [0, m_2 -1] \bigr)$ satisfies $ (4d)i+ i_0 \equiv 1,2,\ldots, 2d \pmod{4d}$
for $i_0 =1,2,\ldots,2d$,
and the corresponding descending sequence of integers
$\bigl( 4m+1-  ((4d)i+ i_0 ) : i\in [0, m_2 -1] \bigr)$ satisfies $4m+1- ((4d)i+ i_0) \equiv 4d, 4d -1,\ldots, 2d +1 \pmod{4d}$
for $i_0 =1,2,\ldots,2d$.
Thus, the sequence of vertices $\bigl( \{  (4d)i+ i_0, 4m+1-  ((4d)i+ i_0 ) \} : i\in [0, m_2 -1] \bigr)$
is a path in $G(4m+1-4d)$ on $m_2 = m/d$ vertices, for all $i_0 \in [1,2d]$.
Hence, $G(4m+1-4d)$ is a disjoint union of $2d$ paths where each path has $m/d$ vertices.
Since $G(4m+1-4d)$ has at least one path on an odd number of vertices,
$G(4m+1-4d)$ has no perfect matching.

\textbf{Case 4.}
Let $k$ be an even positive integer such that $k \ne 2d$ for all positive divisors $d$ of $m$.
Let $k = 2 k_1$ for some positive integer $k_1$.
By the division algorithm, there exist integers $q$ and $r$ such that $4m=kq+r$ and $0\le r < k$.
Let $r_1 = 2m-k_1 q$. Then, $r=2r_1$, $2m=k_1 q + r_1$ and $0\le r_1 < k_1$.
Since $k_1$ is not a divisor of $m$, $r_1 \ne 0$. Thus, $0 < r_1 < k_1$.

Let $i_0 \in [1,r_1]$. The sequence of integers
$\bigl( k i+ i_0 : i\in [0, q] \bigr)$ satisfies $ ki+ i_0 \equiv 1,2,\ldots, r_1 \pmod{2k_1}$,
for $i_0 =1,2,\ldots,r_1$,
and the corresponding descending sequence of integers
$\bigl( 4m+1-  (k i+ i_0 ) : i\in [0, q] \bigr)$ satisfies $4m+1- ( ki+ i_0) \equiv 2r_1, 2r_1 -1,\ldots, r_1 +1 \pmod{2k_1}$,
for $i_0 =1,2,\ldots,r_1$.
Thus, the sequence of vertices $\bigl( \{  ki+ i_0, 4m+1-  (ki+ i_0 ) \} : i\in [0, q] \bigr)$
is a path in $G(4m+1-k)$ on $q+1$ vertices, for all $i_0 \in [1,r_1]$.

Let $i_0 \in [2r_1 +1, k_1 + r_1]$. The sequence of integers
$\bigl( k i+ i_0 : i\in [0, q-1] \bigr)$ satisfies $ ki+ i_0 \equiv 2r_1 +1,2r_1 +2,\ldots, k_1 + r_1 \pmod{2k_1}$
for $i_0 =2r_1 +1,2r_1 +2,\ldots,k_1 + r_1$,
and the corresponding descending sequence of integers
$\bigl( 4m+1-  (k i+ i_0 ) : i\in [0, q-1] \bigr)$ satisfies $4m+1- ( ki+ i_0) \equiv 2k_1, 2k_1 -1,\ldots, k_1 + r_1 +1 \pmod{2k_1}$,
for $i_0 =2r_1 +1,2r_1 +2,\ldots,k_1 + r_1$.
Thus, the sequence of vertices $\bigl( \{  ki+ i_0, 4m+1-  (ki+ i_0 ) \} : i\in [0, q-1] \bigr)$
is a path in $G(4m+1-k)$ on $q$ vertices, for all $i_0 \in [2r_1 +1, k_1 + r_1]$.

Let $v$ be the number of vertices in the $r_1$ paths, where each path has $q+1$ vertices, and $k_1 - r_1$ paths,
where each path has $q$ vertices.
Then, $v = r_1 (q+1) + (k_1 - r_1)q= k_1 q + r_1 =2m$.
Hence, $G(4m+1-k)$ is a disjoint union of $r_1$ paths, where each path has $q+1$ vertices, and $k_1 - r_1$ paths,
where each path has $q$ vertices.
Since $r_1 >0$, $k_1 - r_1 >0$, and $q$ and $q+1$ have opposite parity,
$G(4m+1-k)$ has at least one path on an odd number of vertices.
Therefore, $G(4m+1-k)$ has no perfect matching.

\textbf{Case 5.}
Let $k$ be an odd positive integer such that $k \ne d$ for all positive divisors $d$ of $m$.
Let $k = 2 k_1 +1$ for some nonnegative integer $k_1$.
Since $k$ is not a divisor of $m$, we have $k \ne 1$. Thus, $k_1 > 0$.
Let $\ell = \tfrac{1}{2} (4m+1-k) = 2m - k_1$.
By the division algorithm, there exist integers $q$ and $r$ such that $4m=kq+r$ and $0\le r < k$.
Since $k$ is not a divisor of $m$, $r\ne 0$. Thus, $0 < r <k$.
Since $1\le k < 2m$, we have $q\ge 2$.
We consider the subcases where either $q$ is even or $q$ is odd.

\textbf{Subcase 5(\textit{i}).}
Assume $q$ is even.
Let $q= 2 q_1$ for some positive integer $q_1$.
Let $r_1 = 2m -k q_1$. Then, $r =2 r_1$.
The sequence of integers $\bigl( k i+\ell : i\in [-(q_1 -1), q_1] \bigr)$ satisfies $ ki+\ell \equiv r_1 + k_1 +1 \pmod{2k_1 +1}$,
and the corresponding descending sequence of integers
$\bigl( 4m+1-  (k i+\ell) : i\in [-(q_1 -1), q_1] \bigr)$ satisfies $4m+1- (k i+\ell) \equiv r_1 + k_1 +1 \pmod{2k_1 +1}$.
Thus, the sequence $\bigl( 4m+1-  (k i+\ell) : i\in [-(q_1 -1), q_1] \bigr)$ is the
sequence $\bigl( k i+\ell : i\in [-(q_1 -1), q_1] \bigr)$ in reverse order.
Hence, the sequence of vertices $\bigl( \{  k i+\ell, 4m+1-  (k i+\ell) \} : i\in [-(q_1 -1), 0] \bigr)$
is a path in $G(4m+1-k)$ on $q_1$ vertices.

Let $i_0 \in [1,r_1]$. The sequence of integers
$\bigl( k i+ i_0 : i\in [0, q] \bigr)$ satisfies $ ki+ i_0 \equiv 1,2,\ldots, r_1 \pmod{2k_1 +1}$
for $i_0 =1,2,\ldots,r_1$,
and the corresponding descending sequence of integers
$\bigl( 4m+1-  (k i+ i_0 ) : i\in [0, q] \bigr)$ satisfies $4m+1- ( ki+ i_0) \equiv 2r_1, 2r_1 -1,\ldots, r_1 +1 \pmod{2k_1 +1}$,
for $i_0 =1,2,\ldots,r_1$.
Thus, the sequence of vertices $\bigl( \{  ki+ i_0, 4m+1-  (ki+ i_0 ) \} : i\in [0, q] \bigr)$
is a path in $G(4m+1-k)$ on $q+1$ vertices, for all $i_0 \in [1,r_1]$.

Since $ 2 r_1 = r < k = 2k_1 +1$, we have $r_1\le k_1$.
Let $i_0 \in [2r_1 +1, k_1 + r_1]$. The sequence of integers
$\bigl( k i+ i_0 : i\in [0, q-1] \bigr)$ satisfies $ ki+ i_0 \equiv 2r_1 +1,2r_1 +2,\ldots, k_1 + r_1 \pmod{2k_1 +1}$
for $i_0 =2r_1 +1,2r_1 +2,\ldots,k_1 + r_1$,
and the corresponding descending sequence of integers
$\bigl( 4m+1-  (k i+ i_0 ) : i\in [0, q-1] \bigr)$ satisfies $4m+1- ( ki+ i_0) \equiv 2k_1 +1, 2k_1,\ldots, k_1 + r_1 +2 \pmod{2k_1 +1}$,
for $i_0 =2r_1 +1,2r_1 +2,\ldots,k_1 + r_1$.
Thus, the sequence of vertices $\bigl( \{  ki+ i_0, 4m+1-  (ki+ i_0 ) \} : i\in [0, q-1] \bigr)$
is a path in $G(4m+1-k)$ on $q$ vertices, for all $i_0 \in [2r_1 +1, k_1 + r_1]$.

Let $v$ be the number of vertices in the $r_1$ paths, where each path has $q+1$ vertices, $k_1 - r_1$ paths,
where each path has $q$ vertices, and one path on $q_1$ vertices.
Then, $v = r_1 (q+1) + (k_1 - r_1)q +q_1= k q_1 + r_1 =2m$.
Hence, $G(4m+1-k)$ is a disjoint union of $r_1$ paths, where each path has $q+1$ vertices,
$k_1 - r_1$ paths, where each path has $q$ vertices, and one path on $q_1$ vertices.
Since $r_1 >0$ and $q+1$ is odd,
$G(4m+1-k)$ has at least one path on an odd number of vertices.
Therefore, $G(4m+1-k)$ has no perfect matching.

\textbf{Subcase 5(\textit{ii}).}
Assume $q$ is odd.
Let $q + 1= 2 q_1$ for some positive integer $q_1$.
Then, $r = 4m -kq$ is odd.
Let $r = 2r_1 +1$ for some nonnegative integer $r_1$.
Since $r < k$, we have $r_1 < k_1$.
Recall that $\ell = \tfrac{1}{2} (4m+1 -k) = 2m -k_1$.

The sequence of integers $\bigl( k i+\ell : i\in [-(q_1 -1), q_1] \bigr)$ satisfies $ ki+\ell \equiv r_1  +1 \pmod{2k_1 +1}$
and the corresponding descending sequence of integers
$\bigl( 4m+1-  (k i+\ell) : i\in [-(q_1 -1), q_1] \bigr)$ satisfies $4m+1- (k i+\ell) \equiv r_1  +1 \pmod{2k_1 +1}$.
Thus, the sequence $\bigl( 4m+1-  (k i+\ell) : i\in [-(q_1 -1), q_1] \bigr)$ is the
sequence $\bigl( k i+\ell : i\in [-(q_1 -1), q_1] \bigr)$ in reverse order.
Hence, the sequence of vertices $\bigl( \{  k i+\ell, 4m+1-  (k i+\ell) \} : i\in [-(q_1 -1), 0] \bigr)$
is a path in $G(4m+1-k)$ on $q_1$ vertices.

Let $i_0 \in [1,r_1]$. The sequence of integers
$\bigl( k i+ i_0 : i\in [0, q] \bigr)$ satisfies $ ki+ i_0 \equiv 1,2,\ldots, r_1 \pmod{2k_1 +1}$
for $i_0 =1,2,\ldots,r_1$,
and the corresponding descending sequence of integers
$\bigl( 4m+1-  (k i+ i_0 ) : i\in [0, q] \bigr)$ satisfies $4m+1- ( ki+ i_0) \equiv 2r_1 +1, 2r_1,\ldots, r_1 +2 \pmod{2k_1 +1}$,
for $i_0 =1,2,\ldots,r_1$.
Thus, the sequence of vertices $\bigl( \{  ki+ i_0, 4m+1-  (ki+ i_0 ) \} : i\in [0, q] \bigr)$
is a path in $G(4m+1-k)$ on $q+1$ vertices, for all $i_0 \in [1,r_1]$.

Let $i_0 \in [2r_1 +2, k_1 + r_1 +1]$. The sequence of integers
$\bigl( k i+ i_0 : i\in [0, q-1] \bigr)$ satisfies $ ki+ i_0 \equiv 2r_1 +2,2r_1 +3,\ldots, k_1 + r_1 +1 \pmod{2k_1 +1}$
for $i_0 =2r_1 +2,2r_1 +3,\ldots,k_1 + r_1 +1$,
and the corresponding descending sequence of integers
$\bigl( 4m+1-  (k i+ i_0 ) : i\in [0, q-1] \bigr)$ satisfies $4m+1- ( ki+ i_0) \equiv 2k_1 +1, 2k_1,\ldots, k_1 + r_1 +2 \pmod{2k_1 +1}$,
for $i_0 =2r_1 +2,2r_1 +3,\ldots,k_1 + r_1 +1$.
Thus, the sequence of vertices $\bigl( \{  ki+ i_0, 4m+1-  (ki+ i_0 ) \} : i\in [0, q-1] \bigr)$
is a path in $G(4m+1-k)$ on $q$ vertices, for all $i_0 \in [2r_1 +2, k_1 + r_1 +1]$.

Let $v$ be the number of vertices in the $r_1$ paths, where each path has $q+1$ vertices, $k_1 - r_1$ paths,
where each path has $q$ vertices, and one path on $q_1$ vertices.
Then, $v = r_1 (q+1) + (k_1 - r_1)q +q_1= k q_1 - k_1 + r_1 =2m$.
Hence, $G(4m+1-k)$ is a disjoint union of $r_1$ paths, where each path has $q+1$ vertices,
$k_1 - r_1$ paths, where each path has $q$ vertices, and one path on $q_1$ vertices.
Since $k_1 - r_1 >0$ and $q$ is odd,
$G(4m+1-k)$ has at least one path on an odd number of vertices.
Therefore, $G(4m+1-k)$ has no perfect matching.
\vspace{0.1in}

We define the \textit{perfect matching set} of $m$, denoted by $PMS(m)$, by
\begin{equation*}
  PMS(m) =\{ a=4m+1-k : k=d \mbox{ or } k=2d, \mbox{ for some positive divisor } d \mbox{ of } m\}.
\end{equation*}
Let $a\in PMS(m)$.
Because $G(a)$ is a disjoint union of paths where each path has an even number of vertices,
$G(a)$ has exactly one perfect matching $K_a$.
Furthermore, suppose we are given two integers $a,a' \in PSM(m)$ such that $2m+1 \le a, a'\le 4m$,
and there are Klein bottle equivalent $C_4$-face-magic Klein bottle labeling $X$ and $X'$ for
which $\mathcal{L}(X) = K_{a}$ and $\mathcal{L}(X') = K_{a'}$.
Then, there exists $A\in KBLS(m,4)$ such that $X'= A(X)$.
We have $X=B(L(K_a,\alpha))$ for some $\alpha\in S_n$ for which $\alpha(1)=1$ and $B\in KBLS(m,4)$.
Hence, $K_{a'} = \mathcal{L}((AB)(L(K_a,\alpha)))$.
Thus, $(a)$ is path partition equivalent to $(a')$.
By Lemma \ref{lemmaPathPartitionEquivalence}, we have either $a'=a$ or $a'=S-a$.
Since $a'\le 4m$ and $S-a\ge 4m+2$, we have $a'\ne S-a$. Thus, $a'=a$.
For each value $a\in PSM(m)$, by Corollary \ref{corNumberOfLabelinsOnKM4},
there are $2^{m-1}\, (m-1)!$ distinct Klein bottle
nonequivalent $C_4$-face-magic Klein bottle labelings $X$ on $\mathcal{K}_{m,4}$ such that $\mathcal{L}(X)=K_a$.
We have $|PMS(m)| =2\tau(m)$, where $\tau(m)$ is the number of positive divisors of $m$.
Therefore, there are $2\tau(m) \cdot  2^{m-1}\, (m-1)!  = 2^{m}\, (m-1)! \, \tau(m)$ distinct Klein bottle
nonequivalent $C_4$-face-magic Klein bottle labelings on $\mathcal{K}_{m,4}$.
This completes the proof.
\end{proof}

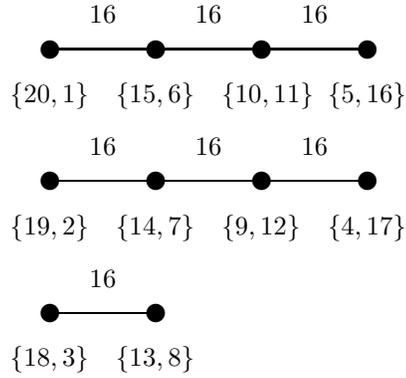
\begin{figure}
\hspace{0.0in}\begin{picture}(300,140)(-50,0)
\put(10,120){\circle*{7}}
\put(-5,100){$\{20,1 \}$}
\put(10,120){\line(1,0){40}}
\put(25,130){16}

\put(50,120){\circle*{7}}
\put(35,100){$\{15,6 \}$}
\put(50,120){\line(1,0){40}}
\put(65,130){16}

\put(90,120){\circle*{7}}
\put(75,100){$\{10,11\}$}
\put(90,120){\line(1,0){40}}
\put(105,130){16}

\put(130,120){\circle*{7}}
\put(115,100){$\{ 5,16 \}$}

\put(10,70){\circle*{7}}
\put(-5,50){$\{ 19,2 \}$}
\put(10,70){\line(1,0){40}}
\put(25,80){16}

\put(50,70){\circle*{7}}
\put(35,50){$\{ 14,7 \}$}
\put(50,70){\line(1,0){40}}
\put(65,80){16}

\put(90,70){\circle*{7}}
\put(75,50){$\{ 9,12 \}$}
\put(90,70){\line(1,0){40}}
\put(105,80){16}

\put(130,70){\circle*{7}}
\put(115,50){$\{ 4,17\}$}

\put(10,20){\circle*{7}}
\put(-5,0){$\{ 18,3 \}$}
\put(10,20){\line(1,0){40}}
\put(25,30){16}

\put(50,20){\circle*{7}}
\put(35,0){$\{ 13,8\}$}

\end{picture}
\caption{Labeling search graph $G(16)$ for the $5\times 4$ Klein bottle grid graph.}
\label{figLabelingSearchGraphG16}
\end{figure}

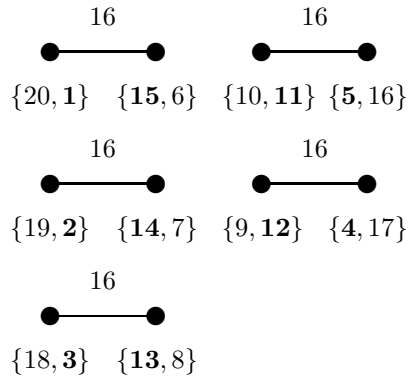
\begin{figure}
\hspace{0.0in}\begin{picture}(300,140)(-50,0)
\put(10,120){\circle*{7}}
\put(-5,100){$\{20,\textbf{1} \}$}
\put(10,120){\line(1,0){40}}
\put(25,130){16}

\put(50,120){\circle*{7}}
\put(35,100){$\{\textbf{15},6 \}$}

\put(90,120){\circle*{7}}
\put(75,100){$\{10,\textbf{11}\}$}
\put(90,120){\line(1,0){40}}
\put(105,130){16}

\put(130,120){\circle*{7}}
\put(115,100){$\{ \textbf{5},16 \}$}

\put(10,70){\circle*{7}}
\put(-5,50){$\{ 19,\textbf{2} \}$}
\put(10,70){\line(1,0){40}}
\put(25,80){16}

\put(50,70){\circle*{7}}
\put(35,50){$\{ \textbf{14},7 \}$}

\put(90,70){\circle*{7}}
\put(75,50){$\{ 9,\textbf{12} \}$}
\put(90,70){\line(1,0){40}}
\put(105,80){16}

\put(130,70){\circle*{7}}
\put(115,50){$\{ \textbf{4},17\}$}

\put(10,20){\circle*{7}}
\put(-5,0){$\{ 18,\textbf{3} \}$}
\put(10,20){\line(1,0){40}}
\put(25,30){16}

\put(50,20){\circle*{7}}
\put(35,0){$\{ \textbf{13},8\}$}

\end{picture}
\caption{$(16)$-admissible path partition $K$ of  $G(16)$ is a perfect matching on $G(16)$.}
\label{figPerfectMatchingLabelingSearchGraphG16}
\end{figure}

\begin{table}[h]
\hspace{0.0in}\begin{tabular}{|c|c|c|c|}
\hline
1   &15  &6    &20    \\[0.05in]
\hline
17  &9   &12   &4     \\[0.05in]
\hline
2   &14  &7    &19     \\[0.05in]
\hline
16  &10  &11   &5     \\[0.05in]
\hline
3   &13  &8    &18    \\[0.05in]
\hline
\end{tabular}
\caption{The standard $C_4$-face-magic Klein bottle labeling on $\mathcal{K}_{5,4}$ associated with $K$
where $K$ is the $(16)$-admissible path partition of $G(16)$ given in Figure \ref{figPerfectMatchingLabelingSearchGraphG16}.}
\label{tableStandardLabelK5by4here}
\end{table}

\begin{example}
  In Figure \ref{figLabelingSearchGraphG16} we have the labeling search graph $G(16)$ for the $5\times 4$ Klein bottle grid graph $\mathcal{K}_{5,4}$.
  In Figure \ref{figPerfectMatchingLabelingSearchGraphG16} we have a $(16)$-admissible path partition $K$
  of $G(16)$ that is a perfect matching on $G(16)$.
  The standard $C_4$-face-magic Klein bottle labeling on $\mathcal{K}_{5,4}$ associated with $K$ denoted by $L(K,\epsilon,\mu)$
  is given in Table \ref{tableStandardLabelK5by4here}.
  See Definition \ref{defnStandardKleinBottleLabeling}.
\end{example}

\section{Results on an $m \times 6$ Klein bottle grid graph} \label{S:resultsonmby6gridgraph}

\begin{notation}
For convenience, we let $S$ be the $C_4$-face-magic value of a $C_4$-face-magic Klein bottle
labeling $X$ of $\mathcal{K}_{m,6}$.
By Lemma \ref{lemmaKleinCheckerboardC4FaceValue}, we have $S=12m+2$.
\end{notation}

\begin{proposition} \label{propMiddleValuesCondition}
Let $m\ge 3$ be an odd integer and let $b$ be an integer such that either $3m+1\le b \le 4m+1$ or $8m \le b \le 9m+1$.
Then, $G(6m,b)$ has a $(6m,b)$-admissible path partition $K_b$ and thus,
there exist $(m-1)!$ Klein bottle nonequivalent $C_4$-face-magic Klein bottle labelings $X$ on $\mathcal{K}_{m,6}$
such that $\mathcal{L}(X)=K_b$.
\end{proposition}

\begin{proof}
\textbf{Case 1.} Suppose $3m+1 \le b \le 4m+1$.
Let $b=3m+k$, where $k$ is an integer such that $1\le k \le m+1$.
We determine a $(6m,b)$-admissible path partition of $G(6m,b)$.
For all integers $i$ such that $1\le i \le k-1$, we let
\begin{align}
x_{i,1,k} &= 6m-2i+1, \label{eqnVertexLabel1MiddleFirst}  \\
\tfrac{1}{2} S -x_{i,1,k} &= 2i, \label{eqnVertexLabel1MiddleOpposeFirst} \\
x_{i,2,k} &= 2i-1, \label{eqnVertexLabel2MiddleFirst} \\
\tfrac{1}{2} S - x_{i,2,k} &= 6m-2i+2, \label{eqnVertexLabel2MiddleOpposeFirst} \\
x_{i,3,k} &= 3m-2i+k+1, \mbox{ and } \label{eqnVertexLabel3MiddleFirst} \\
\tfrac{1}{2} S - x_{i,3,k} &= 3m+2i-k. \label{eqnVertexLabel3MiddleOpposeFirst}
\end{align}
For all integers $i$ such that $k\le i \le m$, we let
\begin{align}
x_{i,1,k} &= 6m-3i+k+1, \label{eqnVertexLabel1MiddleSecond} \\
\tfrac{1}{2} S -x_{i,1,k} &= 3i -k, \label{eqnVertexLabel1MiddleOpposeSecond} \\
x_{i,2,k} &= 3i-k-1,  \label{eqnVertexLabel2MiddleSecond}  \\
\tfrac{1}{2} S - x_{i,2,k} &= 6m-3i+k+2, \label{eqnVertexLabel2MiddleOpposeSecond} \\
x_{i,3,k} &= 3m-3i+2k+1, \mbox{ and }  \label{eqnVertexLabel3MiddleSecond}  \\
\tfrac{1}{2} S - x_{i,3,k} &= 3m+3i-2k. \label{eqnVertexLabel3MiddleOpposeSecond}
\end{align}
From equations (\ref{eqnVertexLabel1MiddleOpposeFirst}) and (\ref{eqnVertexLabel2MiddleFirst}),
for $1 \le i \le k-1$,
we obtain the labels $\{ q\in \mathbb{Z} : 1\le q \le 2k-2\}$.
From equations (\ref{eqnVertexLabel1MiddleOpposeSecond}), (\ref{eqnVertexLabel2MiddleSecond}), and (\ref{eqnVertexLabel3MiddleSecond}),
for $k \le i \le m$,
we obtain the labels $\{ q\in \mathbb{Z} : 2k-1 \le q \le 3m-k+1\}$.
From equation (\ref{eqnVertexLabel3MiddleFirst}), for $\stkout{k} +1 \le i \le k-1$,
and equation (\ref{eqnVertexLabel3MiddleOpposeFirst}), for $1 \le i \le \stkout{k}$,
we obtain the labels $\{ q\in \mathbb{Z} : 3m -k+2 \le q \le 3m\}$.
However, it should be noted that the details of the argument in this last case depends on whether $k$ is even or odd.
Thus for $1\le i \le m$ and $1\le j \le 3$, the vertices $\{x_{i,j,k}, \tfrac{1}{2} S - x_{i,j,k} \}$ in the graph $G(6m,b)$
are distinct.

We observe that for all $1\le i \le m$ and $1\le k \le m+1$, we have $x_{i,1,k}+x_{i,2,k} = 6m$
and $x_{i,2,k} + x_{i,3,k} = 3m+k=b$.
Hence for all $1\le i \le m$ and $1\le k \le m+1$,
the path $P_{i,k} = ( \{ x_{i,j,k}, \tfrac{1}{2} S - x_{i,j,k} \} : 1\le j \le 3)$ is
a $(6m,b)$-admissible path in $G(6m,b)$ with vertex label sequence $( x_{i,j,k} : 1\le j \le 3)$.
Therefore for all $1\le k \le m+1$, $K_{3m+k}= \{ P_{i,k} :1\le i \le m\}$ is a $(6m,b)$-admissible
path partition of $G(6m,b)$.
Since $(6m,b)$ is not a palindrome, by Lemma \ref{lemmaSpanningGraphOfPathsCriteria}, for each $3m+1\le b \le 4m+1$,
there exist $(m-1)!$ Klein bottle nonequivalent $C_4$-face-magic Klein bottle labelings $X$
of $\mathcal{K}_{m,6}$ such that $\mathcal{L}(X)=K_{b}$.

\textbf{Case 2.} Suppose $8m\le b \le 9m+1$.
Let $b=8m+k$, where $k$ is an integer such that $0\le k \le m+1$.
We determine a $(6m,b)$-admissible path partition of $G(6m,b)$.
For all integers $i$ such that $1\le i \le k-1$, we let
\begin{align}
x_{i,1,k} &= 2m-3i+k+1, \label{eqnVertexLabel1CenterFirst}  \\
\tfrac{1}{2} S -x_{i,1,k} &= 4m+3i-k, \label{eqnVertexLabel1CenterOpposeFirst} \\
x_{i,2,k} &= 4m+3i-k-1, \label{eqnVertexLabel2CenterFirst} \\
\tfrac{1}{2} S - x_{i,2,k} &= 2m-3i+k+2, \label{eqnVertexLabel2CenterOpposeFirst} \\
x_{i,3,k} &= 4m-3i+2k+1, \mbox{ and } \label{eqnVertexLabel3CenterFirst} \\
\tfrac{1}{2} S - x_{i,3,k} &= 2m+3i-2k. \label{eqnVertexLabel3CenterOpposeFirst}
\end{align}
For all integers $i$ such that $k\le i \le m$, we let
\begin{align}
x_{i,1,k} &= 2m-2i+1, \label{eqnVertexLabel1CenterSecond} \\
\tfrac{1}{2} S -x_{i,1,k} &= 4m+2i, \label{eqnVertexLabel1CenterOpposeSecond} \\
x_{i,2,k} &= 4m+2i-1,  \label{eqnVertexLabel2CenterSecond}  \\
\tfrac{1}{2} S - x_{i,2,k} &= 2m-2i+2, \label{eqnVertexLabel2CenterOpposeSecond} \\
x_{i,3,k} &= 4m-2i+k+1, \mbox{ and }  \label{eqnVertexLabel3CenterSecond}  \\
\tfrac{1}{2} S - x_{i,3,k} &= 2m+2i-k. \label{eqnVertexLabel3CenterOpposeSecond}
\end{align}
From equations (\ref{eqnVertexLabel3CenterSecond}) and (\ref{eqnVertexLabel3CenterOpposeSecond}),
for $k \le i \le m$,
we obtain the labels $\{ q\in \mathbb{Z} : 1\le q \le 2m-2k+2\}$.
From equations (\ref{eqnVertexLabel1CenterFirst}), (\ref{eqnVertexLabel2CenterOpposeFirst}), and (\ref{eqnVertexLabel3CenterOpposeFirst}),
for $1 \le i \le k-1$,
we obtain the labels $\{ q\in \mathbb{Z} : 2m-2k+3 \le q \le 2m+k-1\}$.
From equation (\ref{eqnVertexLabel3CenterSecond}), for $\stkout{m}^{+} + \, \stkout{k}^{+} \le i \le m$,
and equation (\ref{eqnVertexLabel3CenterOpposeSecond}), for $k \le i \le \stkout{m}^{+} + \, \stkout{k}^{+} -1$,
we obtain the labels $\{ q\in \mathbb{Z} : 2m +k \le q \le 3m\}$.
However, it should be noted that the details of the argument in this last case depends on whether $k$ is even or odd.
Thus for $1\le i \le m$ and $1\le j \le 3$, the vertices $\{x_{i,j,k}, \tfrac{1}{2} S - x_{i,j,k} \}$ in the graph $G(6m,b)$
are distinct.

We observe that, for all $1\le i \le m$ and $0\le k \le m+1$, we have $x_{i,1,k}+x_{i,2,k} = 6m$
and $x_{i,2,k} + x_{i,3,k} = 8m+k=b$. Hence for all $1\le i \le m$ and $0\le k \le m+1$,
the path $P_{i,k} = ( \{ x_{i,j,k}, \tfrac{1}{2} S - x_{i,j,k} \} : 1\le j \le 3)$ is
a $(6m,b)$-admissible path in $G(6m,b)$ with vertex label sequence $( x_{i,j,k} : 1\le j \le 3)$.
Therefore for each $0\le k \le m+1$, $K_{8m+k}= \{ P_{i,k} :1\le i \le m\}$ is a $(6m,b)$-admissible
path partition of $G(6m,b)$.
Since $(6m,b)$ is not a palindrome, by Lemma \ref{lemmaSpanningGraphOfPathsCriteria}, for each $8m \le b \le 9m+1$,
there exist $(m-1)!$ Klein bottle nonequivalent $C_4$-face-magic Klein bottle labelings $X$
of $\mathcal{K}_{m,6}$ such that $\mathcal{L}(X)=K_{b}$.
\end{proof}

\begin{proposition}  \label{propDiagonalValuesCondition}
Let $m\ge 3$ be an odd integer and let $k$ be an integer such that $1\le k \le m+1$.
Then, $G(2m+k,3m+k)$ has a $(2m+k,3m+k)$-admissible path partition $K_k$ and thus,
there are $(m-1)!$ Klein bottle nonequivalent $C_4$-face-magic Klein bottle labelings $X$ on $\mathcal{K}_{m,6}$
such that $\mathcal{L}(X)=K_k$.
\end{proposition}

\begin{proof}
We determine a $(2m+k,3m+k)$-admissible path partition of $G(2m+k,3m+k)$.
For all integers $i$ such that $1\le i \le k-1$, we let
\begin{align}
x_{i,1,k} &= k-i, \label{eqnVertexLabel1DiagonalFirst}  \\
\tfrac{1}{2} S -x_{i,1,k} &= 6m+i-k+1, \label{eqnVertexLabel1DiagonalOpposeFirst} \\
x_{i,2,k} &= 2m+i, \label{eqnVertexLabel2DiagonalFirst} \\
\tfrac{1}{2} S - x_{i,2,k} &= 4m-i+1, \label{eqnVertexLabel2DiagonalOpposeFirst} \\
x_{i,3,k} &= m-i+k, \mbox{ and } \label{eqnVertexLabel3DiagonalFirst} \\
\tfrac{1}{2} S - x_{i,3,k} &= 5m+i-k+1. \label{eqnVertexLabel3DiagonalOpposeFirst}
\end{align}
For all integers $i$ such that $k\le i \le m$, we let
\begin{align}
x_{i,1,k} &= 2m-i+k, \label{eqnVertexLabel1DiagonalSecond} \\
\tfrac{1}{2} S -x_{i,1,k} &= 4m+i-k+1, \label{eqnVertexLabel1DiagonalOpposeSecond} \\
x_{i,2,k} &= i,  \label{eqnVertexLabel2DiagonalSecond}  \\
\tfrac{1}{2} S - x_{i,2,k} &= 6m-i+1, \label{eqnVertexLabel2DiagonalOpposeSecond} \\
x_{i,3,k} &= 3m-i+k, \mbox{ and }  \label{eqnVertexLabel3DiagonalSecond}  \\
\tfrac{1}{2} S - x_{i,3,k} &= 3m+i-k+1. \label{eqnVertexLabel3DiagonalOpposeSecond}
\end{align}
From equation (\ref{eqnVertexLabel1DiagonalFirst}), for $1\le i \le k-1$, and equation (\ref{eqnVertexLabel2DiagonalSecond}),
for $k \le i \le m$,
we obtain the labels $\{ q\in \mathbb{Z} : 1\le q \le m\}$.
From equation (\ref{eqnVertexLabel3DiagonalFirst}), for $1\le i \le k-1$, and equation (\ref{eqnVertexLabel1DiagonalSecond}),
for $k \le i \le m$,
we obtain the labels $\{ q\in \mathbb{Z} : m+1\le q \le 2m\}$.
From equation (\ref{eqnVertexLabel2DiagonalFirst}), for $1\le i \le k-1$, and equation (\ref{eqnVertexLabel3DiagonalSecond}),
for $k \le i \le m$,
we obtain the labels $\{ q\in \mathbb{Z} : 2m+1\le q \le 3m\}$.
Thus for $1\le i \le m$ and $1\le j \le 3$, the vertices $\{x_{i,j}, \tfrac{1}{2} S - x_{i,j} \}$ in the graph $G(2m+k,3m+k)$
are distinct.

We observe that, for all $1\le i \le m$ and $1\le k \le m+1$, we have $x_{i,1,k}+x_{i,2,k} = 2m+k$
and $x_{i,2,k} + x_{i,3,k} = 3m+k$. Hence for all $1\le i \le m$ and $1\le k \le m+1$,
the path $P_{i,k} = ( \{ x_{i,j,k}, \tfrac{1}{2} S - x_{i,j,k} \} : 1\le j \le 3)$ is
a $(2m+k,3m+k)$-admissible path in $G(2m+k,3m+k)$ with vertex label sequence $( x_{i,j,k} : 1\le j \le 3)$.
Therefore, $K_k = \{ P_{i,k} :1\le i \le m\}$ is a $(2m+k,3m+k)$-admissible
path partition of $G(2m+k,3m+k)$.
Since $(2m+k,3m+k)$ is not a palindrome, by Lemma \ref{lemmaSpanningGraphOfPathsCriteria}, for each $1\le k \le m+1$,
there exist $(m-1)!$ Klein bottle nonequivalent $C_4$-face-magic Klein bottle labelings $X$
of $\mathcal{K}_{m,6}$ such that $\mathcal{L}(X)=K_{k}$.
\end{proof}

\section{Results on an $m \times n$ Klein bottle grid graph} \label{S:resultsonmbyncheckerboard}

We first introduce functions $\rho_1$ and $\rho_2$.

\begin{notation} \label{notnRhoValues}
Let $i\in \{ 1,2\}$. Suppose $n_1$ is a positive integer.
We define $\rho_i :\{1,2,\ldots,n_1\} \rightarrow \{1,2\}$ by $p_i(k)=i$ for all $k\in \{1,2,\ldots,n_1\}$.
\end{notation}

\begin{proposition} \label{propLowerSequences}
Let $m\ge 3$ be an odd integer.
Let $n\ge 6$ be an even integer, $n_1=n/2$ and $n_0=n_1 -1$.
Let $S=2(mn+1)$.
Let $(a_1,a_2,\ldots,a_{n_0})$ be one of the following three sequences:
\begin{enumerate}
\item $a_k = mn -n_0 + k$, for all $1\le k \le n_0$,
\item $a_k = m(n_1 + k) +1$, for all $1\le k \le n_0$, or
\item $a_k = 2mk + 1$, for all $1\le k \le n_0$.
\end{enumerate}
Then, $G(a_1,a_2,\ldots,a_{n_0})$ has at least $2^m$ distinct $(a_1,a_2,\ldots,a_{n_0})$-admissible path partitions and thus,
for each distinct $(a_1,a_2,\ldots,a_{n_0})$-admissible path partition $K$ of $G(a_1,a_2,\ldots,a_{n_0})$,
there exist $(m-1)!$ distinct Klein bottle nonequivalent $C_4$-face-magic Klein bottle labelings $X$ on $\mathcal{K}_{m,n}$
such that $\mathcal{L}(X)=K$.
\end{proposition}

\begin{proof}
We consider each sequence $(a_1,a_2,\ldots,a_{n_0})$ individually. The proof of the proposition for each sequence is similar.

\textbf{Case 1.} Suppose $a_k = mn -n_0 + k$, for all $1\le k \le n_0$.
We define the labels of $(a_1,a_2,\ldots,a_{n_0})$-admissible paths in $G(a_1,a_2,\ldots,a_{n_0})$ as follows:
for all $1 \le i \le m$, we let
\begin{align*}
x_{1,i,2j-1}                   &= n_1 (i-1) +j,         &\mbox{ for }  1 \le j \le \stkout{n}_1^{+},  \\
\tfrac{1}{2} S - x_{1,i,2j-1}  &= mn -n_1 (i-1) - j +1,  &\mbox{ for }  1 \le j \le \stkout{n}_1^{+},   \\
x_{1,i,2j}                     &= mn -n_1 i +j,           &\mbox{ for }  1 \le j \le \stkout{n}_1,      \\
\tfrac{1}{2} S - x_{1,i,2j}    &=  n_1 i -j+1,             &\mbox{ for }  1 \le j \le \stkout{n}_1,   \\
x_{2,i,2j-1}                   &=  mn -n_1 i +j,         &\mbox{ for }  1 \le j \le \stkout{n}_1^{+},    \\
\tfrac{1}{2} S - x_{2,i,2j-1}  &=  n_1 i -j+1,           &\mbox{ for }  1 \le j \le \stkout{n}_1^{+},     \\
x_{2,i,2j}                     &= n_1 (i-1) +j,           &\mbox{ for }  1 \le j \le \stkout{n}_1, \mbox{ and }   \\
\tfrac{1}{2} S - x_{2,i,2j}    &= mn -n_1 (i-1) - j +1,    &\mbox{ for }  1 \le j \le \stkout{n}_1.
\end{align*}
Let $\rho:\{1,2,\ldots,m\} \rightarrow \{1,2\}$. Let $P_{\rho,i}$ be the path
on the sequence of vertices $( \{ x_{\rho(i),i,j}, \tfrac{1}{2}S\allowbreak  - x_{\rho(i),i,j} \} : 1 \le j \le n_1 )$
with vertex label sequence $( x_{\rho(i),i,j} : 1 \le j \le n_1 )$.
The sets of labels
\begin{align*}
\{ x_{1,i,2j-1} : 1 \le j \le \stkout{n}_1^{+} \} \cup \{ \tfrac{1}{2} S - x_{1,i,2j} : 1 \le j \le \stkout{n}_1 \} \mbox{ and } \\
\{ x_{2,i,2j} : 1 \le j \le \stkout{n}_1 \} \cup \{ \tfrac{1}{2} S - x_{2,i,2j-1} : 1 \le j \le \stkout{n}_1^{+} \}
\end{align*}
are the set $\{ n_1 (i-1) +k: 1\le k \le n_1 \}$.
Hence for $1\le i \le m$ and $1 \le j \le n_1$, the vertices $\{ x_{\rho(i),i,j}, \tfrac{1}{2} S - x_{\rho(i),i,j} \}$
are distinct in the graph $G(a_1,a_2,\ldots,a_{n_0})$.

Furthermore for $1\le i \le m$, we have
\begin{align*}
x_{\rho(i),i,2j-1} + x_{\rho(i),i,2j} &= (  n_1 (i-1) +j ) + ( mn -n_1 i +j ) & \\
                          &= mn - n_0 +(2j-1) =a_{2j-1}, &\mbox{ for }  1 \le j \le \stkout{n}_1, \mbox{ and } \\
x_{\rho(i),i,2j} + x_{\rho(i),i,2j+1} &= x_{\rho(i),i,2j} + x_{\rho(i),i,2j-1} +1  & \\
             &=  mn - n_0 +2j =a_{2j}, &\mbox{ for }  1 \le j \le \stkout{n}_1^{-}.
\end{align*}
Thus, $P_{\rho,i}$ is an $(a_1,a_2,\ldots,a_{n_0})$-admissible path in $G(a_1,a_2,\ldots,a_{n_0})$.

We consider the example where $m=3$, $n_1=5$ and $n=10$. Then, Figures \ref{figGraphProp32Case1K1}
and \ref{figGraphProp32Case1K2} provide the $(27,28,29,30)$-admissible path partitions $K_{\rho_1}$
and $K_{\rho_2}$ in the graph $G(27,28,29,30)$.
See Notation \ref{notnRhoValues} for the definitions of $\rho_1$ and $\rho_2$.

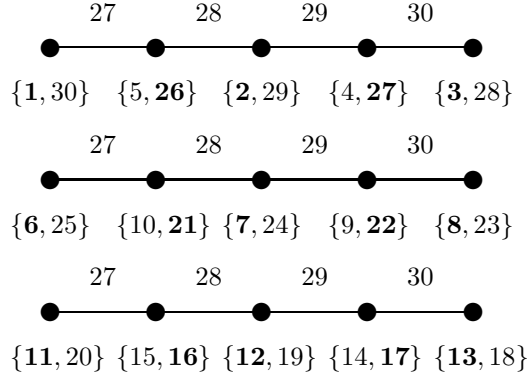
\begin{figure}
\hspace{0.0in}\begin{picture}(300,150)(-50,0)
\put(10,120){\circle*{7}}
\put(-5,100){$\{ \textbf{1},30\}$}
\put(10,120){\line(1,0){40}}
\put(25,130){27}

\put(50,120){\circle*{7}}
\put(35,100){$\{ 5,\textbf{26}\}$}
\put(50,120){\line(1,0){40}}
\put(65,130){28}

\put(90,120){\circle*{7}}
\put(75,100){$\{ \textbf{2},29\}$}
\put(90,120){\line(1,0){40}}
\put(105,130){29}

\put(130,120){\circle*{7}}
\put(115,100){$\{ 4,\textbf{27}\}$}
\put(130,120){\line(1,0){40}}
\put(145,130){30}

\put(170,120){\circle*{7}}
\put(155,100){$\{ \textbf{3},28\}$}

\put(10,70){\circle*{7}}
\put(-5,50){$\{ \textbf{6},25\}$}
\put(10,70){\line(1,0){40}}
\put(25,80){27}

\put(50,70){\circle*{7}}
\put(35,50){$\{ 10,\textbf{21}\}$}
\put(50,70){\line(1,0){40}}
\put(65,80){28}

\put(90,70){\circle*{7}}
\put(75,50){$\{ \textbf{7},24\}$}
\put(90,70){\line(1,0){40}}
\put(105,80){29}

\put(130,70){\circle*{7}}
\put(115,50){$\{ 9,\textbf{22}\}$}
\put(130,70){\line(1,0){40}}
\put(145,80){30}

\put(170,70){\circle*{7}}
\put(155,50){$\{ \textbf{8},23\}$}

\put(10,20){\circle*{7}}
\put(-5,0){$\{\textbf{11},20\}$}
\put(10,20){\line(1,0){40}}
\put(25,30){27}

\put(50,20){\circle*{7}}
\put(35,0){$\{ 15,\textbf{16}\}$}
\put(50,20){\line(1,0){40}}
\put(65,30){28}

\put(90,20){\circle*{7}}
\put(75,0){$\{ \textbf{12},19\}$}
\put(90,20){\line(1,0){40}}
\put(105,30){29}

\put(130,20){\circle*{7}}
\put(115,0){$\{ 14,\textbf{17}\}$}
\put(130,20){\line(1,0){40}}
\put(145,30){30}

\put(170,20){\circle*{7}}
\put(155,0){$\{ \textbf{13},18\}$}
\end{picture}
\caption{$(27,28,29,30)$-admissible path partition $K_{\rho_{1}}$ in $G(27,28,29,30)$.}
\label{figGraphProp32Case1K1}
\end{figure}

\begin{figure}
\hspace{0.0in}\begin{picture}(300,150)(-50,0)
\put(10,120){\circle*{7}}
\put(-5,100){$\{ 5,\textbf{26}\}$}
\put(10,120){\line(1,0){40}}
\put(25,130){27}

\put(50,120){\circle*{7}}
\put(35,100){$\{ \textbf{1},30\}$}
\put(50,120){\line(1,0){40}}
\put(65,130){28}

\put(90,120){\circle*{7}}
\put(75,100){$\{ 4,\textbf{27}\}$}
\put(90,120){\line(1,0){40}}
\put(105,130){29}

\put(130,120){\circle*{7}}
\put(115,100){$\{ \textbf{2},29\}$}
\put(130,120){\line(1,0){40}}
\put(145,130){30}

\put(170,120){\circle*{7}}
\put(155,100){$\{ 3,\textbf{28}\}$}

\put(10,70){\circle*{7}}
\put(-5,50){$\{ 10,\textbf{21}\}$}
\put(10,70){\line(1,0){40}}
\put(25,80){27}

\put(50,70){\circle*{7}}
\put(35,50){$\{ \textbf{6},25\}$}
\put(50,70){\line(1,0){40}}
\put(65,80){28}

\put(90,70){\circle*{7}}
\put(75,50){$\{ 9,\textbf{22}\}$}
\put(90,70){\line(1,0){40}}
\put(105,80){29}

\put(130,70){\circle*{7}}
\put(115,50){$\{ \textbf{7},24\}$}
\put(130,70){\line(1,0){40}}
\put(145,80){30}

\put(170,70){\circle*{7}}
\put(155,50){$\{ 8,\textbf{23}\}$}

\put(10,20){\circle*{7}}
\put(-5,0){$\{ 15,\textbf{16}\}$}
\put(10,20){\line(1,0){40}}
\put(25,30){27}

\put(50,20){\circle*{7}}
\put(35,0){$\{ \textbf{11},20\}$}
\put(50,20){\line(1,0){40}}
\put(65,30){28}

\put(90,20){\circle*{7}}
\put(75,0){$\{ 14,\textbf{17}\}$}
\put(90,20){\line(1,0){40}}
\put(105,30){29}

\put(130,20){\circle*{7}}
\put(115,0){$\{ \textbf{12},19\}$}
\put(130,20){\line(1,0){40}}
\put(145,30){30}

\put(170,20){\circle*{7}}
\put(155,0){$\{ 13,\textbf{18}\}$}
\end{picture}
\caption{$(27,28,29,30)$-admissible path partition $K_{\rho_{2}}$ in $G(27,28,29,30)$.}
\label{figGraphProp32Case1K2}
\end{figure}
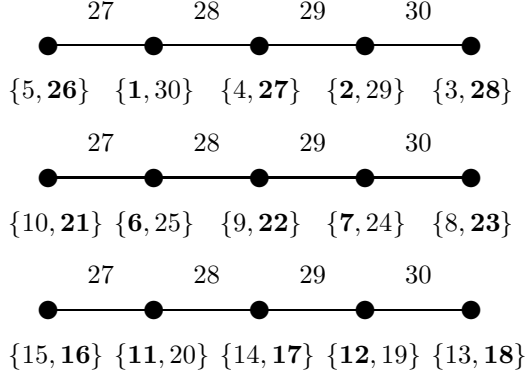

\textbf{Case 2.} Suppose $a_k = m(n_1 +k) + 1$, for all $1\le k \le n_0$.
We define the labels of $(a_1,a_2,\ldots,a_{n_0})$-admissible paths in $G(a_1,a_2,\ldots,a_{n_0})$ as follows:
for all $1 \le i \le m$, we let
\begin{align*}
x_{1,i,2j-1}                   &= m (j-1) +i,         &\mbox{ for }  1 \le j \le \stkout{n}_1^{+},  \\
\tfrac{1}{2} S - x_{1,i,2j-1}  &= m(n -j+1) -i +1,  &\mbox{ for }  1 \le j \le \stkout{n}_1^{+},   \\
x_{1,i,2j}                     &= m(n_1 +j) -i +1,           &\mbox{ for }  1 \le j \le \stkout{n}_1,      \\
\tfrac{1}{2} S - x_{1,i,2j}    &= m(n_1 -j) +i,             &\mbox{ for }  1 \le j \le \stkout{n}_1,   \\
x_{2,i,2j-1}                   &= m(n_1 +j) -i +1,         &\mbox{ for }  1 \le j \le \stkout{n}_1^{+},    \\
\tfrac{1}{2} S - x_{2,i,2j-1}  &= m(n_1 -j) +i,           &\mbox{ for }  1 \le j \le \stkout{n}_1^{+},     \\
x_{2,i,2j}                     &= m (j-1) +i,           &\mbox{ for }  1 \le j \le \stkout{n}_1, \mbox{ and }   \\
\tfrac{1}{2} S - x_{2,i,2j}    &= m(n -j+1) -i +1,    &\mbox{ for }  1 \le j \le \stkout{n}_1.
\end{align*}
Let $\rho:\{1,2,\ldots,m\} \rightarrow \{1,2\}$. Let $P_{\rho,i}$ be the path
on the sequence of vertices $( \{ x_{\rho(i),i,j}, \tfrac{1}{2}S\allowbreak  - x_{\rho(i),i,j} \} : 1 \le j \le n_1 )$
with vertex label sequence $( x_{\rho(i),i,j} : 1 \le j \le n_1 )$.
The sets of labels
\begin{align*}
\{ x_{1,i,2j-1} : 1 \le j \le \stkout{n}_1^{+} \} \cup \{ \tfrac{1}{2} S - x_{1,i,2j} : 1 \le j \le \stkout{n}_1 \} \mbox{ and } \\
\{ x_{2,i,2j} : 1 \le j \le \stkout{n}_1 \} \cup \{ \tfrac{1}{2} S - x_{2,i,2j-1} : 1 \le j \le \stkout{n}_1^{+} \}
\end{align*}
are the set $\{ mk+i: 0\le k \le n_1 -1 \}$.
Hence for $1\le i \le m$ and $1 \le j \le n_1$, the vertices $\{ x_{\rho(i),i,j}, \tfrac{1}{2} S - x_{\rho(i),i,j} \}$
are distinct in the graph $G(a_1,a_2,\ldots,a_{n_0})$.

Furthermore, for $1\le i \le m$, we have
\begin{align*}
x_{\rho(i),i,2j-1} + x_{\rho(i),i,2j} &= ( m (j-1) +i ) + (  m(n_1 +j) -i +1 ) & \\
                          &= m(n_1 +2j-1) + 1 =a_{2j-1}, &\mbox{ for }  1 \le j \le \stkout{n}_1, \mbox{ and } \\
x_{\rho(i),i,2j} + x_{\rho(i),i,2j+1} &=  x_{\rho(i),i,2j} + x_{\rho(i),i,2j-1} +m & \\
             &= m( n_1 +2j) +1 =a_{2j}, &\mbox{ for }  1 \le j \le \stkout{n}_1^{-}.
\end{align*}
Thus, $P_{\rho,i}$ is an $(a_1,a_2,\ldots,a_{n_0})$-admissible path in $G(a_1,a_2,\ldots,a_{n_0})$.

We consider the example where $m=3$, $n_1=5$ and $n=10$. Then, Figures \ref{figGraphProp32Case2K1}
and \ref{figGraphProp32Case2K2} provide the $(19,22,25,28)$-admissible path partitions $K_{\rho_1}$
and $K_{\rho_2}$ in the graph $G(19,22,25,28)$.
See Notation \ref{notnRhoValues} for the definitions of $\rho_1$ and $\rho_2$.

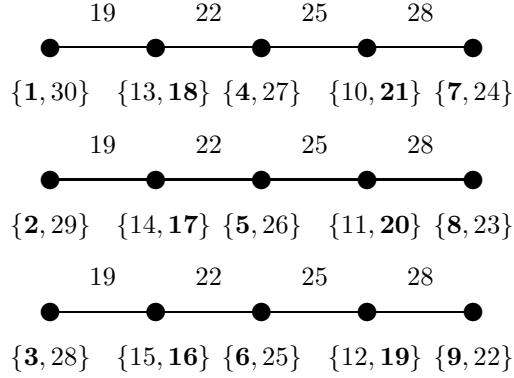
\begin{figure}
\hspace{0.0in}\begin{picture}(300,150)(-50,0)
\put(10,120){\circle*{7}}
\put(-5,100){$\{ \textbf{1},30\}$}
\put(10,120){\line(1,0){40}}
\put(25,130){19}

\put(50,120){\circle*{7}}
\put(35,100){$\{ 13,\textbf{18}\}$}
\put(50,120){\line(1,0){40}}
\put(65,130){22}

\put(90,120){\circle*{7}}
\put(75,100){$\{ \textbf{4},27\}$}
\put(90,120){\line(1,0){40}}
\put(105,130){25}

\put(130,120){\circle*{7}}
\put(115,100){$\{ 10,\textbf{21}\}$}
\put(130,120){\line(1,0){40}}
\put(145,130){28}

\put(170,120){\circle*{7}}
\put(155,100){$\{ \textbf{7},24\}$}

\put(10,70){\circle*{7}}
\put(-5,50){$\{ \textbf{2},29\}$}
\put(10,70){\line(1,0){40}}
\put(25,80){19}

\put(50,70){\circle*{7}}
\put(35,50){$\{ 14,\textbf{17}\}$}
\put(50,70){\line(1,0){40}}
\put(65,80){22}

\put(90,70){\circle*{7}}
\put(75,50){$\{ \textbf{5},26\}$}
\put(90,70){\line(1,0){40}}
\put(105,80){25}

\put(130,70){\circle*{7}}
\put(115,50){$\{ 11,\textbf{20}\}$}
\put(130,70){\line(1,0){40}}
\put(145,80){28}

\put(170,70){\circle*{7}}
\put(155,50){$\{ \textbf{8},23\}$}

\put(10,20){\circle*{7}}
\put(-5,0){$\{\textbf{3},28\}$}
\put(10,20){\line(1,0){40}}
\put(25,30){19}

\put(50,20){\circle*{7}}
\put(35,0){$\{ 15,\textbf{16}\}$}
\put(50,20){\line(1,0){40}}
\put(65,30){22}

\put(90,20){\circle*{7}}
\put(75,0){$\{ \textbf{6},25\}$}
\put(90,20){\line(1,0){40}}
\put(105,30){25}

\put(130,20){\circle*{7}}
\put(115,0){$\{ 12,\textbf{19}\}$}
\put(130,20){\line(1,0){40}}
\put(145,30){28}

\put(170,20){\circle*{7}}
\put(155,0){$\{ \textbf{9},22\}$}
\end{picture}
\caption{$(19,22,25,28)$-admissible path partition $K_{\rho_{1}}$ in $G(19,22,25,28)$.}
\label{figGraphProp32Case2K1}
\end{figure}

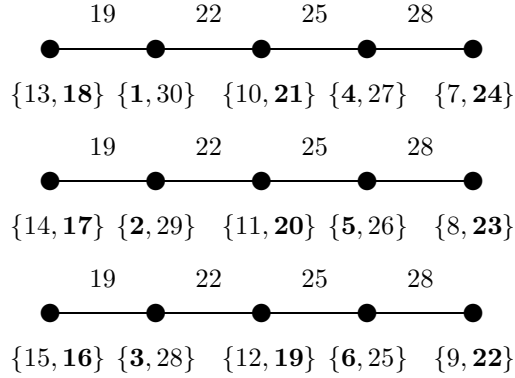
\begin{figure}
\hspace{0.0in}\begin{picture}(300,150)(-50,0)
\put(10,120){\circle*{7}}
\put(-5,100){$\{ 13,\textbf{18}\}$}
\put(10,120){\line(1,0){40}}
\put(25,130){19}

\put(50,120){\circle*{7}}
\put(35,100){$\{ \textbf{1},30\}$}
\put(50,120){\line(1,0){40}}
\put(65,130){22}

\put(90,120){\circle*{7}}
\put(75,100){$\{ 10,\textbf{21}\}$}
\put(90,120){\line(1,0){40}}
\put(105,130){25}

\put(130,120){\circle*{7}}
\put(115,100){$\{ \textbf{4},27\}$}
\put(130,120){\line(1,0){40}}
\put(145,130){28}

\put(170,120){\circle*{7}}
\put(155,100){$\{ 7,\textbf{24}\}$}

\put(10,70){\circle*{7}}
\put(-5,50){$\{ 14,\textbf{17}\}$}
\put(10,70){\line(1,0){40}}
\put(25,80){19}

\put(50,70){\circle*{7}}
\put(35,50){$\{ \textbf{2},29\}$}
\put(50,70){\line(1,0){40}}
\put(65,80){22}

\put(90,70){\circle*{7}}
\put(75,50){$\{ 11,\textbf{20}\}$}
\put(90,70){\line(1,0){40}}
\put(105,80){25}

\put(130,70){\circle*{7}}
\put(115,50){$\{ \textbf{5},26\}$}
\put(130,70){\line(1,0){40}}
\put(145,80){28}

\put(170,70){\circle*{7}}
\put(155,50){$\{ 8,\textbf{23}\}$}

\put(10,20){\circle*{7}}
\put(-5,0){$\{ 15,\textbf{16}\}$}
\put(10,20){\line(1,0){40}}
\put(25,30){19}

\put(50,20){\circle*{7}}
\put(35,0){$\{ \textbf{3},28\}$}
\put(50,20){\line(1,0){40}}
\put(65,30){22}

\put(90,20){\circle*{7}}
\put(75,0){$\{ 12,\textbf{19}\}$}
\put(90,20){\line(1,0){40}}
\put(105,30){25}

\put(130,20){\circle*{7}}
\put(115,0){$\{ \textbf{6},25\}$}
\put(130,20){\line(1,0){40}}
\put(145,30){28}

\put(170,20){\circle*{7}}
\put(155,0){$\{ 9,\textbf{22}\}$}
\end{picture}
\caption{$(19,22,25,28)$-admissible path partition $K_{\rho_{2}}$ in $G(19,22,25,28)$.}
\label{figGraphProp32Case2K2}
\end{figure}

\textbf{Case 3.} Suppose $a_k = 2mk + 1$, for all $1\le k \le n_0$.
We define the labels of $(a_1,a_2,\ldots,a_{n_0})$-admissible paths in $G(a_1,a_2,\ldots,a_{n_0})$ as follows:
for all $1 \le i \le m$, we let
\begin{align*}
x_{1,i,2j-1}                   &= 2m (j-1) +i,         &\mbox{ for }  1 \le j \le \stkout{n}_1^{+},  \\
\tfrac{1}{2} S - x_{1,i,2j-1}  &= m(n -2j+2) -i +1,  &\mbox{ for }  1 \le j \le \stkout{n}_1^{+},   \\
x_{1,i,2j}                     &= 2m j -i +1,           &\mbox{ for }  1 \le j \le \stkout{n}_1,      \\
\tfrac{1}{2} S - x_{1,i,2j}    &= m(n -2j) +i,             &\mbox{ for }  1 \le j \le \stkout{n}_1,   \\
x_{2,i,2j-1}                   &= 2m j -i +1,         &\mbox{ for }  1 \le j \le \stkout{n}_1^{+},    \\
\tfrac{1}{2} S - x_{2,i,2j-1}  &= m(n -2j) +i,           &\mbox{ for }  1 \le j \le \stkout{n}_1^{+},     \\
x_{2,i,2j}                     &= 2m (j-1) +i,           &\mbox{ for }  1 \le j \le \stkout{n}_1, \mbox{ and }   \\
\tfrac{1}{2} S - x_{2,i,2j}    &= m(n -2j+2) -i +1,    &\mbox{ for }  1 \le j \le \stkout{n}_1.
\end{align*}
Let $\rho:\{1,2,\ldots,m\} \rightarrow \{1,2\}$. Let $P_{\rho,i}$ be the path
on the sequence of vertices $( \{ x_{\rho(i),i,j}, \tfrac{1}{2}S\allowbreak  - x_{\rho(i),i,j} \} : 1 \le j \le n_1 )$
with vertex label sequence $( x_{\rho(i),i,j} : 1 \le j \le n_1 )$.
The sets of labels
\begin{align*}
\{ x_{1,i,2j-1} : 1 \le j \le \stkout{n}_1^{+} \} \cup \{ \tfrac{1}{2} S - x_{1,i,2j} : 1 \le j \le \stkout{n}_1 \} \mbox{ and } \\
\{ x_{2,i,2j} : 1 \le j \le \stkout{n}_1 \} \cup \{ \tfrac{1}{2} S - x_{2,i,2j-1} : 1 \le j \le \stkout{n}_1^{+} \}
\end{align*}
are the set $\{ 2mk+i: 0\le k \le n_1 -1 \}$.
Also, the sets of labels
\begin{align*}
\{ \tfrac{1}{2} S - x_{1,i,2j-1} : 1 \le j \le \stkout{n}_1^{+} \} \cup \{  x_{1,i,2j} : 1 \le j \le \stkout{n}_1 \} \mbox{ and } \\
\{ \tfrac{1}{2} S - x_{2,i,2j} : 1 \le j \le \stkout{n}_1 \} \cup \{ x_{2,i,2j-1} : 1 \le j \le \stkout{n}_1^{+} \}
\end{align*}
are the set $\{ 2mk-i+1: 1\le k \le n_1 \}$.
Hence for $1\le i \le m$ and $1 \le j \le n_1$, the vertices $\{ x_{\rho(i),i,j}, \tfrac{1}{2} S - x_{\rho(i),i,j} \}$
are distinct in the graph $G(a_1,a_2,\ldots,a_{n_0})$.

Furthermore, for $1\le i \le m$, we have
\begin{align*}
x_{\rho(i),i,2j-1} + x_{\rho(i),i,2j} &= ( 2m (j-1) +i ) + ( 2m j -i +1 ) & \\
                          &= 2m(2j-1) + 1 =a_{2j-1}, &\mbox{ for }  1 \le j \le \stkout{n}_1, \mbox{ and } \\
x_{\rho(i),i,2j} + x_{\rho(i),i,2j+1} &=  x_{\rho(i),i,2j} + x_{\rho(i),i,2j-1} +2m & \\
             &= 2m(2j) +1 =a_{2j}, &\mbox{ for }  1 \le j \le \stkout{n}_1^{-}.
\end{align*}
Thus, $P_{\rho,i}$ is an $(a_1,a_2,\ldots,a_{n_0})$-admissible path in $G(a_1,a_2,\ldots,a_{n_0})$.

We consider the example where $m=3$, $n_1=5$ and $n=10$. Then, Figures \ref{figGraphProp32Case3K1}
and \ref{figGraphProp32Case3K2} provide the $(7,13,19,25)$-admissible path partitions $K_{\rho_1}$
and $K_{\rho_2}$ in the graph $G(7,13,19,25)$.
See Notation \ref{notnRhoValues} for the definitions of $\rho_1$ and $\rho_2$.

\begin{figure}
\hspace{0.0in}\begin{picture}(300,150)(-50,0)
\put(10,120){\circle*{7}}
\put(-5,100){$\{ \textbf{1},30\}$}
\put(10,120){\line(1,0){40}}
\put(25,130){7}

\put(50,120){\circle*{7}}
\put(35,100){$\{ \textbf{6},25\}$}
\put(50,120){\line(1,0){40}}
\put(65,130){13}

\put(90,120){\circle*{7}}
\put(75,100){$\{ \textbf{7},24\}$}
\put(90,120){\line(1,0){40}}
\put(105,130){19}

\put(130,120){\circle*{7}}
\put(115,100){$\{ \textbf{12},19\}$}
\put(130,120){\line(1,0){40}}
\put(145,130){25}

\put(170,120){\circle*{7}}
\put(155,100){$\{ \textbf{13},18\}$}

\put(10,70){\circle*{7}}
\put(-5,50){$\{ \textbf{2},29\}$}
\put(10,70){\line(1,0){40}}
\put(25,80){7}

\put(50,70){\circle*{7}}
\put(35,50){$\{ \textbf{5},26\}$}
\put(50,70){\line(1,0){40}}
\put(65,80){13}

\put(90,70){\circle*{7}}
\put(75,50){$\{ \textbf{8},23\}$}
\put(90,70){\line(1,0){40}}
\put(105,80){19}

\put(130,70){\circle*{7}}
\put(115,50){$\{ \textbf{11},20\}$}
\put(130,70){\line(1,0){40}}
\put(145,80){25}

\put(170,70){\circle*{7}}
\put(155,50){$\{ \textbf{14},17\}$}

\put(10,20){\circle*{7}}
\put(-5,0){$\{\textbf{3},28\}$}
\put(10,20){\line(1,0){40}}
\put(25,30){7}

\put(50,20){\circle*{7}}
\put(35,0){$\{ \textbf{4},27\}$}
\put(50,20){\line(1,0){40}}
\put(65,30){13}

\put(90,20){\circle*{7}}
\put(75,0){$\{ \textbf{9},22\}$}
\put(90,20){\line(1,0){40}}
\put(105,30){19}

\put(130,20){\circle*{7}}
\put(115,0){$\{ \textbf{10},21\}$}
\put(130,20){\line(1,0){40}}
\put(145,30){25}

\put(170,20){\circle*{7}}
\put(155,0){$\{ \textbf{15},16\}$}
\end{picture}
\caption{$(7,13,19,25)$-admissible path partition $K_{\rho_{1}}$ in $G(7,13,19,25)$.}
\label{figGraphProp32Case3K1}
\end{figure}
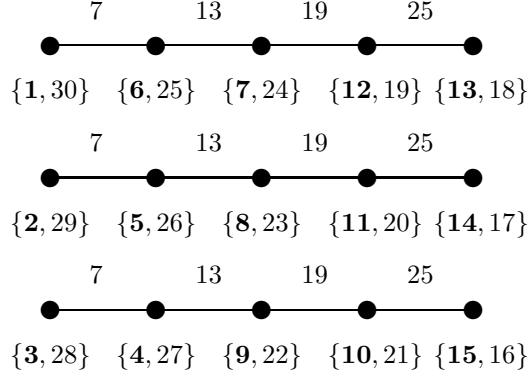

\begin{figure}
\hspace{0.0in}\begin{picture}(300,150)(-50,0)
\put(10,120){\circle*{7}}
\put(-5,100){$\{ \textbf{6},25\}$}
\put(10,120){\line(1,0){40}}
\put(25,130){7}

\put(50,120){\circle*{7}}
\put(35,100){$\{ \textbf{1},30\}$}
\put(50,120){\line(1,0){40}}
\put(65,130){13}

\put(90,120){\circle*{7}}
\put(75,100){$\{ \textbf{12},19\}$}
\put(90,120){\line(1,0){40}}
\put(105,130){19}

\put(130,120){\circle*{7}}
\put(115,100){$\{ \textbf{7},24\}$}
\put(130,120){\line(1,0){40}}
\put(145,130){25}

\put(170,120){\circle*{7}}
\put(155,100){$\{ 13,\textbf{18}\}$}

\put(10,70){\circle*{7}}
\put(-5,50){$\{ \textbf{5},26\}$}
\put(10,70){\line(1,0){40}}
\put(25,80){7}

\put(50,70){\circle*{7}}
\put(35,50){$\{ \textbf{2},29\}$}
\put(50,70){\line(1,0){40}}
\put(65,80){13}

\put(90,70){\circle*{7}}
\put(75,50){$\{ \textbf{11},20\}$}
\put(90,70){\line(1,0){40}}
\put(105,80){19}

\put(130,70){\circle*{7}}
\put(115,50){$\{ \textbf{8},23\}$}
\put(130,70){\line(1,0){40}}
\put(145,80){25}

\put(170,70){\circle*{7}}
\put(155,50){$\{ 14,\textbf{17}\}$}

\put(10,20){\circle*{7}}
\put(-5,0){$\{ \textbf{4},27\}$}
\put(10,20){\line(1,0){40}}
\put(25,30){7}

\put(50,20){\circle*{7}}
\put(35,0){$\{ \textbf{3},28\}$}
\put(50,20){\line(1,0){40}}
\put(65,30){13}

\put(90,20){\circle*{7}}
\put(75,0){$\{ \textbf{10},21\}$}
\put(90,20){\line(1,0){40}}
\put(105,30){19}

\put(130,20){\circle*{7}}
\put(115,0){$\{ \textbf{9},22\}$}
\put(130,20){\line(1,0){40}}
\put(145,30){25}

\put(170,20){\circle*{7}}
\put(155,0){$\{ 15,\textbf{16}\}$}
\end{picture}
\caption{$(7,13,19,25)$-admissible path partition $K_{\rho_{2}}$ in $G(7,13,19,25)$.}
\label{figGraphProp32Case3K2}
\end{figure}
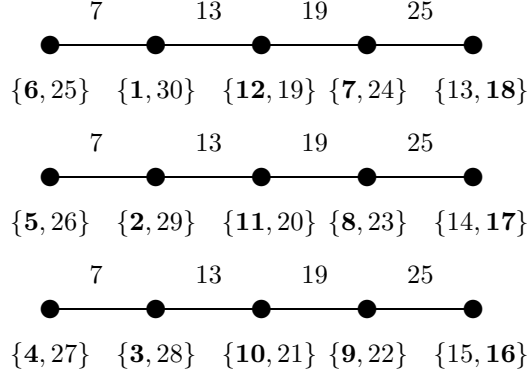
\vspace{0.1in}

Therefore in each of Cases 1, 2 and 3,
$K_{\rho} = \{ P_{\rho,i} : 1\le i \le m \}$ is a $(a_1,a_2,\ldots, \allowbreak a_{n_0})$-admissible
path partition of $G(a_1,a_2,\ldots,a_{n_0})$.
Since $(a_1,a_2,\ldots,a_{n_0})$ is not a palindrome, by Lemma \ref{lemmaSpanningGraphOfPathsCriteria}, for each $\rho$,
there exist $(m-1)!$ Klein bottle nonequivalent $C_4$-face-magic Klein bottle labelings $X$
of $\mathcal{K}_{m,n}$ such that $\mathcal{L}(X)=K_{\rho}$.
Hence, there exist at least $2^m (m-1)!$ distinct Klein bottle nonequivalent $C_4$-face-magic Klein bottle labelings $X$
of $\mathcal{K}_{m,n}$ such that $\mathcal{L}(X)$ is an $(a_1,a_2,\ldots,a_{n_0})$-admissible path
partition of $G(a_1,a_2,\ldots,a_{n_0})$.
\end{proof}

\begin{proposition} \label{propMiddleSequences}
Let $m\ge 3$ be an odd integer.
Let $n\ge 6$ be an even integer such that $n_1=n/2$ is odd.
Let $n_0=n_1 -1$ and $S=2(mn+1)$.
Let $(a_1,a_2,\ldots,a_{n_0})$ be one of the following three sequences:
\begin{enumerate}
\item $a_k = (2m-1)n_1 +2k +1$, for all $1\le k \le n_0$,
\item $a_k = m n_1 +2mk +1$, for all $1\le k \le n_0$, or
\item $a_k = 4mk +1$, for all $1\le k \le n_0$.
\end{enumerate}
Then, $G(a_1,a_2,\ldots,a_{n_0})$ has at least $2^m$ distinct $(a_1,a_2,\ldots,a_{n_0})$-admissible path partitions and thus,
for each distinct $(a_1,a_2,\ldots,a_{n_0})$-admissible path partition $K$ of $G(a_1,a_2,\ldots,a_{n_0})$,
there exist $(m-1)!$ distinct Klein bottle nonequivalent $C_4$-face-magic Klein bottle labelings $X$ on $\mathcal{K}_{m,n}$
such that $\mathcal{L}(X)=K$.
\end{proposition}

\begin{proof}
We consider each sequence $(a_1,a_2,\ldots,a_{n_0})$ individually. The proof of the proposition for each sequence is similar.

\textbf{Case 1.} Suppose $a_k = (2m-1)n_1 +2k +1$, for all $1\le k \le n_0$.
We define the labels of $(a_1,a_2,\ldots,a_{n_0})$-admissible paths in $G(a_1,a_2,\ldots,a_{n_0})$ as follows:
for all $1 \le i \le m$, we let
\begin{align*}
x_{1,i,2j-1}                   &= n_1 (i-1) +2j -1,         &\mbox{ for }  1 \le j \le \stkout{n}_1^{+},  \\
\tfrac{1}{2} S - x_{1,i,2j-1}  &= mn -n_1 (i-1) - 2j +2,  &\mbox{ for }  1 \le j \le \stkout{n}_1^{+},   \\
x_{1,i,2j}                     &= mn -n_1 i +2j,           &\mbox{ for }  1 \le j \le \stkout{n}_1,      \\
\tfrac{1}{2} S - x_{1,i,2j}    &=  n_1 i -2j+1,             &\mbox{ for }  1 \le j \le \stkout{n}_1,   \\
x_{2,i,2j-1}                   &=  mn -n_1 i +2j-1,         &\mbox{ for }  1 \le j \le \stkout{n}_1^{+},    \\
\tfrac{1}{2} S - x_{2,i,2j-1}  &=  n_1 i -2j+2,           &\mbox{ for }  1 \le j \le \stkout{n}_1^{+},     \\
x_{2,i,2j}                     &= n_1 (i-1) +2j,           &\mbox{ for }  1 \le j \le \stkout{n}_1, \mbox{ and }   \\
\tfrac{1}{2} S - x_{2,i,2j}    &= mn -n_1 (i-1) - 2j +1,    &\mbox{ for }  1 \le j \le \stkout{n}_1.
\end{align*}
Let $\rho:\{1,2,\ldots,m\} \rightarrow \{1,2\}$. Let $P_{\rho,i}$ be the path
on the sequence of vertices $( \{ x_{\rho(i),i,j}, \tfrac{1}{2}S\allowbreak  - x_{\rho(i),i,j} \} : 1 \le j \le n_1 )$
with vertex label sequence $( x_{\rho(i),i,j} : 1 \le j \le n_1 )$.
The sets of labels
\begin{align*}
\{ x_{1,i,2j-1} : 1 \le j \le \stkout{n}_1^{+} \} \cup \{ \tfrac{1}{2} S - x_{1,i,2j} : 1 \le j \le \stkout{n}_1 \} \mbox{ and } \\
\{ x_{2,i,2j} : 1 \le j \le \stkout{n}_1 \} \cup \{ \tfrac{1}{2} S - x_{2,i,2j-1} : 1 \le j \le \stkout{n}_1^{+} \}
\end{align*}
are the set $\{ n_1 (i-1) +k: 1\le k \le n_1 \}$.
Hence for $1\le i \le m$ and $1 \le j \le n_1$, the vertices $\{ x_{\rho(i),i,j}, \tfrac{1}{2} S - x_{\rho(i),i,j} \}$
are distinct in the graph $G(a_1,a_2,\ldots,a_{n_0})$.

Furthermore, for $1\le i \le m$, we have
\begin{align*}
x_{\rho(i),i,2j-1} + x_{\rho(i),i,2j} &= (  n_1 (i-1) +2j -1 ) + ( mn -n_1 i +2j )  \\
                          &= (2m-1)n_1 +2(2j-1) +1 =a_{2j-1}, \mbox{ \  for }  1 \le j \le \stkout{n}_1, \mbox{ and } \\
x_{\rho(i),i,2j} + x_{\rho(i),i,2j+1} &= x_{\rho(i),i,2j} + x_{\rho(i),i,2j-1} +2   \\
             &=  (2m-1)n_1 + 2(2j) +1 =a_{2j}, \mbox{ \ \ \ \ for }  1 \le j \le \stkout{n}_1^{-}.
\end{align*}
Thus, $P_{\rho,i}$ is an $(a_1,a_2,\ldots,a_{n_0})$-admissible path in $G(a_1,a_2,\ldots,a_{n_0})$.

We consider the example where $m=3$, $n_1=5$ and $n=10$. Then, Figures \ref{figGraphProp33Case1K1}
and \ref{figGraphProp33Case1K2} provide the $(28,30,32,34)$-admissible path partitions $K_{\rho_1}$
and $K_{\rho_2}$ in the graph $G(28,30,32,34)$.
See Notation \ref{notnRhoValues} for the definitions of $\rho_1$ and $\rho_2$.

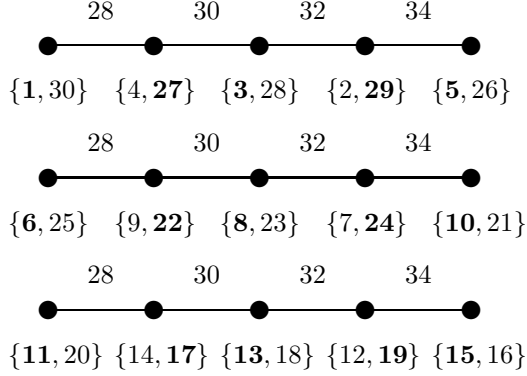
\begin{figure}
\hspace{0.0in}\begin{picture}(300,150)(-50,0)
\put(10,120){\circle*{7}}
\put(-5,100){$\{ \textbf{1},30\}$}
\put(10,120){\line(1,0){40}}
\put(25,130){28}

\put(50,120){\circle*{7}}
\put(35,100){$\{ 4,\textbf{27}\}$}
\put(50,120){\line(1,0){40}}
\put(65,130){30}

\put(90,120){\circle*{7}}
\put(75,100){$\{ \textbf{3},28\}$}
\put(90,120){\line(1,0){40}}
\put(105,130){32}

\put(130,120){\circle*{7}}
\put(115,100){$\{ 2,\textbf{29}\}$}
\put(130,120){\line(1,0){40}}
\put(145,130){34}

\put(170,120){\circle*{7}}
\put(155,100){$\{ \textbf{5},26\}$}

\put(10,70){\circle*{7}}
\put(-5,50){$\{ \textbf{6},25\}$}
\put(10,70){\line(1,0){40}}
\put(25,80){28}

\put(50,70){\circle*{7}}
\put(35,50){$\{ 9,\textbf{22}\}$}
\put(50,70){\line(1,0){40}}
\put(65,80){30}

\put(90,70){\circle*{7}}
\put(75,50){$\{ \textbf{8},23\}$}
\put(90,70){\line(1,0){40}}
\put(105,80){32}

\put(130,70){\circle*{7}}
\put(115,50){$\{ 7,\textbf{24}\}$}
\put(130,70){\line(1,0){40}}
\put(145,80){34}

\put(170,70){\circle*{7}}
\put(155,50){$\{ \textbf{10},21\}$}

\put(10,20){\circle*{7}}
\put(-5,0){$\{\textbf{11},20\}$}
\put(10,20){\line(1,0){40}}
\put(25,30){28}

\put(50,20){\circle*{7}}
\put(35,0){$\{ 14,\textbf{17}\}$}
\put(50,20){\line(1,0){40}}
\put(65,30){30}

\put(90,20){\circle*{7}}
\put(75,0){$\{ \textbf{13},18\}$}
\put(90,20){\line(1,0){40}}
\put(105,30){32}

\put(130,20){\circle*{7}}
\put(115,0){$\{ 12,\textbf{19}\}$}
\put(130,20){\line(1,0){40}}
\put(145,30){34}

\put(170,20){\circle*{7}}
\put(155,0){$\{ \textbf{15},16\}$}
\end{picture}
\caption{$(28,30,32,34)$-admissible path partition $K_{\rho_{1}}$ in $G(28,30,32,34)$.}
\label{figGraphProp33Case1K1}
\end{figure}

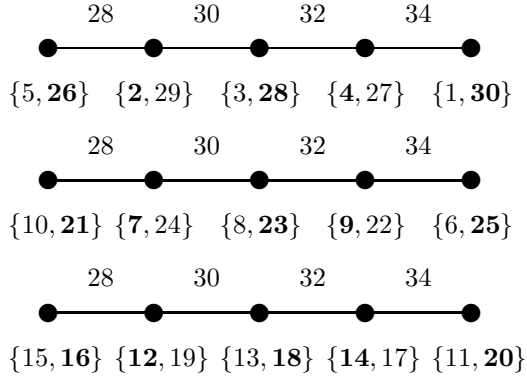
\begin{figure}
\hspace{0.0in}\begin{picture}(300,150)(-50,0)
\put(10,120){\circle*{7}}
\put(-5,100){$\{ 5,\textbf{26}\}$}
\put(10,120){\line(1,0){40}}
\put(25,130){28}

\put(50,120){\circle*{7}}
\put(35,100){$\{ \textbf{2},29\}$}
\put(50,120){\line(1,0){40}}
\put(65,130){30}

\put(90,120){\circle*{7}}
\put(75,100){$\{ 3,\textbf{28}\}$}
\put(90,120){\line(1,0){40}}
\put(105,130){32}

\put(130,120){\circle*{7}}
\put(115,100){$\{ \textbf{4},27\}$}
\put(130,120){\line(1,0){40}}
\put(145,130){34}

\put(170,120){\circle*{7}}
\put(155,100){$\{ 1,\textbf{30}\}$}

\put(10,70){\circle*{7}}
\put(-5,50){$\{ 10,\textbf{21}\}$}
\put(10,70){\line(1,0){40}}
\put(25,80){28}

\put(50,70){\circle*{7}}
\put(35,50){$\{ \textbf{7},24\}$}
\put(50,70){\line(1,0){40}}
\put(65,80){30}

\put(90,70){\circle*{7}}
\put(75,50){$\{ 8,\textbf{23}\}$}
\put(90,70){\line(1,0){40}}
\put(105,80){32}

\put(130,70){\circle*{7}}
\put(115,50){$\{ \textbf{9},22\}$}
\put(130,70){\line(1,0){40}}
\put(145,80){34}

\put(170,70){\circle*{7}}
\put(155,50){$\{ 6,\textbf{25}\}$}

\put(10,20){\circle*{7}}
\put(-5,0){$\{15,\textbf{16}\}$}
\put(10,20){\line(1,0){40}}
\put(25,30){28}

\put(50,20){\circle*{7}}
\put(35,0){$\{ \textbf{12},19\}$}
\put(50,20){\line(1,0){40}}
\put(65,30){30}

\put(90,20){\circle*{7}}
\put(75,0){$\{ 13,\textbf{18}\}$}
\put(90,20){\line(1,0){40}}
\put(105,30){32}

\put(130,20){\circle*{7}}
\put(115,0){$\{ \textbf{14},17\}$}
\put(130,20){\line(1,0){40}}
\put(145,30){34}

\put(170,20){\circle*{7}}
\put(155,0){$\{ 11,\textbf{20}\}$}
\end{picture}
\caption{$(28,30,32,34)$-admissible path partition $K_{\rho_{2}}$ in $G(28,30,32,34)$.}
\label{figGraphProp33Case1K2}
\end{figure}

\textbf{Case 2.} Suppose $a_k = m n_1 +2mk +1$, for all $1\le k \le n_0$.
We define the labels of $(a_1,a_2,\ldots,a_{n_0})$-admissible paths in $G(a_1,a_2,\ldots,a_{n_0})$ as follows:
for all $1 \le i \le m$, we let
\begin{align*}
x_{1,i,2j-1}                   &= m (2j-2) +i,         &\mbox{ for }  1 \le j \le \stkout{n}_1^{+},  \\
\tfrac{1}{2} S - x_{1,i,2j-1}  &= m(n -2j+2) -i +1,  &\mbox{ for }  1 \le j \le \stkout{n}_1^{+},   \\
x_{1,i,2j}                     &= m(n_1 + 2j) -i +1,           &\mbox{ for }  1 \le j \le \stkout{n}_1,      \\
\tfrac{1}{2} S - x_{1,i,2j}    &= m(n_1 - 2j) +i,             &\mbox{ for }  1 \le j \le \stkout{n}_1,   \\
x_{2,i,2j-1}                   &= m(n_1 +2j -1) -i +1,         &\mbox{ for }  1 \le j \le \stkout{n}_1^{+},    \\
\tfrac{1}{2} S - x_{2,i,2j-1}  &= m(n_1 - 2j +1) +i,           &\mbox{ for }  1 \le j \le \stkout{n}_1^{+},     \\
x_{2,i,2j}                     &= m (2j-1) +i,           &\mbox{ for }  1 \le j \le \stkout{n}_1, \mbox{ and }   \\
\tfrac{1}{2} S - x_{2,i,2j}    &= m(n -2j+1) -i +1,    &\mbox{ for }  1 \le j \le \stkout{n}_1.
\end{align*}
Let $\rho:\{1,2,\ldots,m\} \rightarrow \{1,2\}$. Let $P_{\rho,i}$ be the path
on the sequence of vertices $( \{ x_{\rho(i),i,j}, \tfrac{1}{2}S\allowbreak  - x_{\rho(i),i,j} \} : 1 \le j \le n_1 )$
with vertex label sequence $( x_{\rho(i),i,j} : 1 \le j \le n_1 )$.
The sets of labels
\begin{align*}
\{ x_{1,i,2j-1} : 1 \le j \le \stkout{n}_1^{+} \} \cup \{ \tfrac{1}{2} S - x_{1,i,2j} : 1 \le j \le \stkout{n}_1 \} \mbox{ and } \\
\{ x_{2,i,2j} : 1 \le j \le \stkout{n}_1 \} \cup \{ \tfrac{1}{2} S - x_{2,i,2j-1} : 1 \le j \le \stkout{n}_1^{+} \}
\end{align*}
are the set $\{ mk+i: 0\le k \le n_1 -1 \}$.
Hence for $1\le i \le m$ and $1 \le j \le n_1$, the vertices $\{ x_{\rho(i),i,j}, \tfrac{1}{2} S - x_{\rho(i),i,j} \}$
are distinct in the graph $G(a_1,a_2,\ldots,a_{n_0})$.

Furthermore, for $1\le i \le m$, we have
\begin{align*}
x_{\rho(i),i,2j-1} + x_{\rho(i),i,2j} &= m n_1 +2m (2j-1) +1 =a_{2j-1}, &\mbox{ for }  1 \le j \le \stkout{n}_1, \mbox{ and } \\
x_{\rho(i),i,2j} + x_{\rho(i),i,2j+1} &=  x_{\rho(i),i,2j} + x_{\rho(i),i,2j-1} +2m  =a_{2j}, &\mbox{ for }  1 \le j \le \stkout{n}_1^{-}.
\end{align*}
Thus, $P_{\rho,i}$ is an $(a_1,a_2,\ldots,a_{n_0})$-admissible path in $G(a_1,a_2,\ldots,a_{n_0})$.

We consider the example where $m=3$, $n_1=5$ and $n=10$. Then, Figures \ref{figGraphProp33Case2K1}
and \ref{figGraphProp33Case2K2} provide the $(22,28,34,40)$-admissible path partitions $K_{\rho_1}$
and $K_{\rho_2}$ in the graph $G(22,28,34,40)$.
See Notation \ref{notnRhoValues} for the definitions of $\rho_1$ and $\rho_2$.

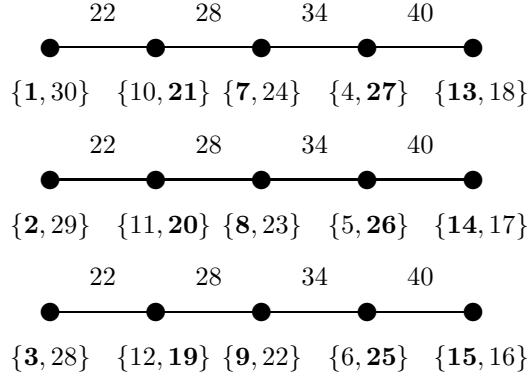
\begin{figure}
\hspace{0.0in}\begin{picture}(300,150)(-50,0)
\put(10,120){\circle*{7}}
\put(-5,100){$\{ \textbf{1},30\}$}
\put(10,120){\line(1,0){40}}
\put(25,130){22}

\put(50,120){\circle*{7}}
\put(35,100){$\{ 10,\textbf{21}\}$}
\put(50,120){\line(1,0){40}}
\put(65,130){28}

\put(90,120){\circle*{7}}
\put(75,100){$\{ \textbf{7},24\}$}
\put(90,120){\line(1,0){40}}
\put(105,130){34}

\put(130,120){\circle*{7}}
\put(115,100){$\{ 4,\textbf{27}\}$}
\put(130,120){\line(1,0){40}}
\put(145,130){40}

\put(170,120){\circle*{7}}
\put(155,100){$\{ \textbf{13},18\}$}

\put(10,70){\circle*{7}}
\put(-5,50){$\{ \textbf{2},29\}$}
\put(10,70){\line(1,0){40}}
\put(25,80){22}

\put(50,70){\circle*{7}}
\put(35,50){$\{ 11,\textbf{20}\}$}
\put(50,70){\line(1,0){40}}
\put(65,80){28}

\put(90,70){\circle*{7}}
\put(75,50){$\{ \textbf{8},23\}$}
\put(90,70){\line(1,0){40}}
\put(105,80){34}

\put(130,70){\circle*{7}}
\put(115,50){$\{ 5,\textbf{26}\}$}
\put(130,70){\line(1,0){40}}
\put(145,80){40}

\put(170,70){\circle*{7}}
\put(155,50){$\{ \textbf{14},17\}$}

\put(10,20){\circle*{7}}
\put(-5,0){$\{\textbf{3},28\}$}
\put(10,20){\line(1,0){40}}
\put(25,30){22}

\put(50,20){\circle*{7}}
\put(35,0){$\{ 12,\textbf{19}\}$}
\put(50,20){\line(1,0){40}}
\put(65,30){28}

\put(90,20){\circle*{7}}
\put(75,0){$\{ \textbf{9},22\}$}
\put(90,20){\line(1,0){40}}
\put(105,30){34}

\put(130,20){\circle*{7}}
\put(115,0){$\{ 6,\textbf{25}\}$}
\put(130,20){\line(1,0){40}}
\put(145,30){40}

\put(170,20){\circle*{7}}
\put(155,0){$\{ \textbf{15},16\}$}
\end{picture}
\caption{$(22,28,34,40)$-admissible path partition $K_{\rho_{1}}$ in $G(22,28,34,40)$.}
\label{figGraphProp33Case2K1}
\end{figure}

\begin{figure}
\hspace{0.0in}\begin{picture}(300,150)(-50,0)
\put(10,120){\circle*{7}}
\put(-5,100){$\{ 13,\textbf{18}\}$}
\put(10,120){\line(1,0){40}}
\put(25,130){22}

\put(50,120){\circle*{7}}
\put(35,100){$\{ \textbf{4},27\}$}
\put(50,120){\line(1,0){40}}
\put(65,130){28}

\put(90,120){\circle*{7}}
\put(75,100){$\{ 7,\textbf{24}\}$}
\put(90,120){\line(1,0){40}}
\put(105,130){34}

\put(130,120){\circle*{7}}
\put(115,100){$\{ \textbf{10},21\}$}
\put(130,120){\line(1,0){40}}
\put(145,130){40}

\put(170,120){\circle*{7}}
\put(155,100){$\{ 1,\textbf{30}\}$}

\put(10,70){\circle*{7}}
\put(-5,50){$\{ 14,\textbf{17}\}$}
\put(10,70){\line(1,0){40}}
\put(25,80){22}

\put(50,70){\circle*{7}}
\put(35,50){$\{ \textbf{5},26\}$}
\put(50,70){\line(1,0){40}}
\put(65,80){28}

\put(90,70){\circle*{7}}
\put(75,50){$\{ 8,\textbf{23}\}$}
\put(90,70){\line(1,0){40}}
\put(105,80){34}

\put(130,70){\circle*{7}}
\put(115,50){$\{ \textbf{11},20\}$}
\put(130,70){\line(1,0){40}}
\put(145,80){40}

\put(170,70){\circle*{7}}
\put(155,50){$\{ 2,\textbf{29}\}$}

\put(10,20){\circle*{7}}
\put(-5,0){$\{15,\textbf{16}\}$}
\put(10,20){\line(1,0){40}}
\put(25,30){22}

\put(50,20){\circle*{7}}
\put(35,0){$\{ \textbf{6},25\}$}
\put(50,20){\line(1,0){40}}
\put(65,30){28}

\put(90,20){\circle*{7}}
\put(75,0){$\{ 9,\textbf{22}\}$}
\put(90,20){\line(1,0){40}}
\put(105,30){34}

\put(130,20){\circle*{7}}
\put(115,0){$\{ \textbf{12},19\}$}
\put(130,20){\line(1,0){40}}
\put(145,30){40}

\put(170,20){\circle*{7}}
\put(155,0){$\{ 3,\textbf{28}\}$}
\end{picture}
\caption{$(22,28,34,40)$-admissible path partition $K_{\rho_{2}}$ in $G(22,28,34,40)$.}
\label{figGraphProp33Case2K2}
\end{figure}
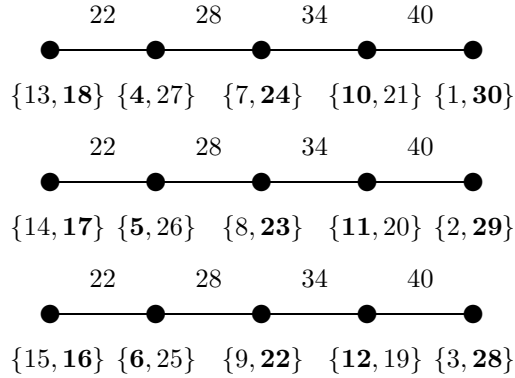

\textbf{Case 3.} Suppose $a_k = 4mk +1$, for all $1\le k \le n_0$.
We define the labels of $(a_1,a_2,\ldots,a_{n_0})$-admissible paths in $G(a_1,a_2,\ldots,a_{n_0})$ as follows:
for all $1 \le i \le m$, we let
\begin{align*}
x_{1,i,2j-1}                   &= 4m (j-1) +i,         &\mbox{ for }  1 \le j \le \stkout{n}_1^{+},  \\
\tfrac{1}{2} S - x_{1,i,2j-1}  &= m(n -4j+4) -i +1,  &\mbox{ for }  1 \le j \le \stkout{n}_1^{+},   \\
x_{1,i,2j}                     &= 4m j -i +1,           &\mbox{ for }  1 \le j \le \stkout{n}_1,      \\
\tfrac{1}{2} S - x_{1,i,2j}    &= m(n -4j) +i,             &\mbox{ for }  1 \le j \le \stkout{n}_1,   \\
x_{2,i,2j-1}                   &= m (4j -2) -i +1,         &\mbox{ for }  1 \le j \le \stkout{n}_1^{+},    \\
\tfrac{1}{2} S - x_{2,i,2j-1}  &= m(n -4j +2) +i,           &\mbox{ for }  1 \le j \le \stkout{n}_1^{+},     \\
x_{2,i,2j}                     &= m (4j-2) +i,           &\mbox{ for }  1 \le j \le \stkout{n}_1, \mbox{ and }   \\
\tfrac{1}{2} S - x_{2,i,2j}    &= m(n -4j+2) -i +1,    &\mbox{ for }  1 \le j \le \stkout{n}_1.
\end{align*}
Let $\rho:\{1,2,\ldots,m\} \rightarrow \{1,2\}$. Let $P_{\rho,i}$ be the path
on the sequence of vertices $( \{ x_{\rho(i),i,j}, \tfrac{1}{2}S\allowbreak  - x_{\rho(i),i,j} \} : 1 \le j \le n_1 )$
with vertex label sequence $( x_{\rho(i),i,j} : 1 \le j \le n_1 )$.
The sets of labels
\begin{align*}
\{ x_{1,i,2j-1} : 1 \le j \le \stkout{n}_1^{+} \} \cup \{ \tfrac{1}{2} S - x_{1,i,2j} : 1 \le j \le \stkout{n}_1 \} \mbox{ and } \\
\{ x_{2,i,2j} : 1 \le j \le \stkout{n}_1 \} \cup \{ \tfrac{1}{2} S - x_{2,i,2j-1} : 1 \le j \le \stkout{n}_1^{+} \}
\end{align*}
are the set $\{ 2mk+i: 0\le k \le n_1 -1 \}$.
Also, the sets of labels
\begin{align*}
\{ \tfrac{1}{2} S - x_{1,i,2j-1} : 1 \le j \le \stkout{n}_1^{+} \} \cup \{  x_{1,i,2j} : 1 \le j \le \stkout{n}_1 \} \mbox{ and } \\
\{ x_{2,i,2j} : 1 \le j \le \stkout{n}_1 \} \cup \{ \tfrac{1}{2} S -  x_{2,i,2j-1} : 1 \le j \le \stkout{n}_1^{+} \}
\end{align*}
are the set $\{ 2mk-i+1: 1\le k \le n_1 \}$.
Hence for $1\le i \le m$ and $1 \le j \le n_1$, the vertices $\{ x_{\rho(i),i,j}, \tfrac{1}{2} S - x_{\rho(i),i,j} \}$
are distinct in the graph $G(a_1,a_2,\ldots,a_{n_0})$.

Furthermore, for $1\le i \le m$, we have
\begin{align*}
x_{\rho(i),i,2j-1} + x_{\rho(i),i,2j} &= 4m(2j-1) + 1 =a_{2j-1}, &\mbox{ for }  1 \le j \le \stkout{n}_1, \mbox{ and } \\
x_{\rho(i),i,2j} + x_{\rho(i),i,2j+1} &=  x_{\rho(i),i,2j} + x_{\rho(i),i,2j-1} +4m =a_{2j}, &\mbox{ for }  1 \le j \le \stkout{n}_1^{-}.
\end{align*}
Thus, $P_{\rho,i}$ is an $(a_1,a_2,\ldots,a_{n_0})$-admissible path in $G(a_1,a_2,\ldots,a_{n_0})$.

We consider the example where $m=3$, $n_1=5$ and $n=10$. Then, Figures \ref{figGraphProp33Case3K1}
and \ref{figGraphProp33Case3K2} provide the $(13,25,37,49)$-admissible path partitions $K_{\rho_1}$
and $K_{\rho_2}$ in the graph $G(13,25,37,49)$.
See Notation \ref{notnRhoValues} for the definitions of $\rho_1$ and $\rho_2$.

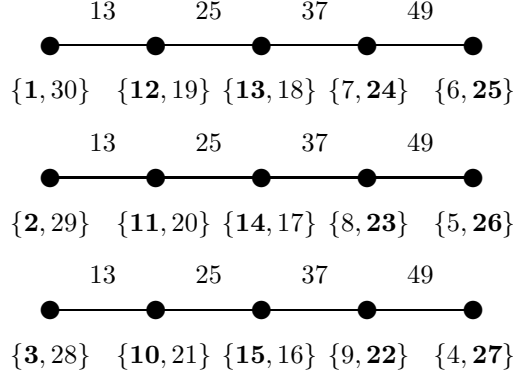
\begin{figure}
\hspace{0.0in}\begin{picture}(300,150)(-50,0)
\put(10,120){\circle*{7}}
\put(-5,100){$\{ \textbf{1},30\}$}
\put(10,120){\line(1,0){40}}
\put(25,130){13}

\put(50,120){\circle*{7}}
\put(35,100){$\{ \textbf{12},19\}$}
\put(50,120){\line(1,0){40}}
\put(65,130){25}

\put(90,120){\circle*{7}}
\put(75,100){$\{ \textbf{13},18\}$}
\put(90,120){\line(1,0){40}}
\put(105,130){37}

\put(130,120){\circle*{7}}
\put(115,100){$\{ 7,\textbf{24}\}$}
\put(130,120){\line(1,0){40}}
\put(145,130){49}

\put(170,120){\circle*{7}}
\put(155,100){$\{ 6,\textbf{25}\}$}

\put(10,70){\circle*{7}}
\put(-5,50){$\{ \textbf{2},29\}$}
\put(10,70){\line(1,0){40}}
\put(25,80){13}

\put(50,70){\circle*{7}}
\put(35,50){$\{ \textbf{11},20\}$}
\put(50,70){\line(1,0){40}}
\put(65,80){25}

\put(90,70){\circle*{7}}
\put(75,50){$\{ \textbf{14},17\}$}
\put(90,70){\line(1,0){40}}
\put(105,80){37}

\put(130,70){\circle*{7}}
\put(115,50){$\{ 8,\textbf{23}\}$}
\put(130,70){\line(1,0){40}}
\put(145,80){49}

\put(170,70){\circle*{7}}
\put(155,50){$\{ 5,\textbf{26}\}$}

\put(10,20){\circle*{7}}
\put(-5,0){$\{\textbf{3},28\}$}
\put(10,20){\line(1,0){40}}
\put(25,30){13}

\put(50,20){\circle*{7}}
\put(35,0){$\{ \textbf{10},21\}$}
\put(50,20){\line(1,0){40}}
\put(65,30){25}

\put(90,20){\circle*{7}}
\put(75,0){$\{ \textbf{15},16\}$}
\put(90,20){\line(1,0){40}}
\put(105,30){37}

\put(130,20){\circle*{7}}
\put(115,0){$\{ 9,\textbf{22}\}$}
\put(130,20){\line(1,0){40}}
\put(145,30){49}

\put(170,20){\circle*{7}}
\put(155,0){$\{ 4,\textbf{27}\}$}
\end{picture}
\caption{$(13,25,37,49)$-admissible path partition $K_{\rho_{1}}$ in $G(13,25,37,49)$.}
\label{figGraphProp33Case3K1}
\end{figure}

\begin{figure}
\hspace{0.0in}\begin{picture}(300,150)(-50,0)
\put(10,120){\circle*{7}}
\put(-5,100){$\{ \textbf{6},25\}$}
\put(10,120){\line(1,0){40}}
\put(25,130){13}

\put(50,120){\circle*{7}}
\put(35,100){$\{ \textbf{7},24\}$}
\put(50,120){\line(1,0){40}}
\put(65,130){25}

\put(90,120){\circle*{7}}
\put(75,100){$\{ 13,\textbf{18}\}$}
\put(90,120){\line(1,0){40}}
\put(105,130){37}

\put(130,120){\circle*{7}}
\put(115,100){$\{ 12,\textbf{19}\}$}
\put(130,120){\line(1,0){40}}
\put(145,130){49}

\put(170,120){\circle*{7}}
\put(155,100){$\{ 1,\textbf{30}\}$}

\put(10,70){\circle*{7}}
\put(-5,50){$\{ \textbf{5},26\}$}
\put(10,70){\line(1,0){40}}
\put(25,80){13}

\put(50,70){\circle*{7}}
\put(35,50){$\{ \textbf{8},23\}$}
\put(50,70){\line(1,0){40}}
\put(65,80){25}

\put(90,70){\circle*{7}}
\put(75,50){$\{ 14,\textbf{17}\}$}
\put(90,70){\line(1,0){40}}
\put(105,80){37}

\put(130,70){\circle*{7}}
\put(115,50){$\{ 11,\textbf{20}\}$}
\put(130,70){\line(1,0){40}}
\put(145,80){49}

\put(170,70){\circle*{7}}
\put(155,50){$\{ 2,\textbf{29}\}$}

\put(10,20){\circle*{7}}
\put(-5,0){$\{ \textbf{4},27\}$}
\put(10,20){\line(1,0){40}}
\put(25,30){22}

\put(50,20){\circle*{7}}
\put(35,0){$\{ \textbf{9},22\}$}
\put(50,20){\line(1,0){40}}
\put(65,30){28}

\put(90,20){\circle*{7}}
\put(75,0){$\{ 15,\textbf{16}\}$}
\put(90,20){\line(1,0){40}}
\put(105,30){34}

\put(130,20){\circle*{7}}
\put(115,0){$\{ 10,\textbf{21}\}$}
\put(130,20){\line(1,0){40}}
\put(145,30){40}

\put(170,20){\circle*{7}}
\put(155,0){$\{ 3,\textbf{28}\}$}
\end{picture}
\caption{$(13,25,37,49)$-admissible path partition $K_{\rho_{2}}$ in $G(13,25,37,49)$.}
\label{figGraphProp33Case3K2}
\end{figure}
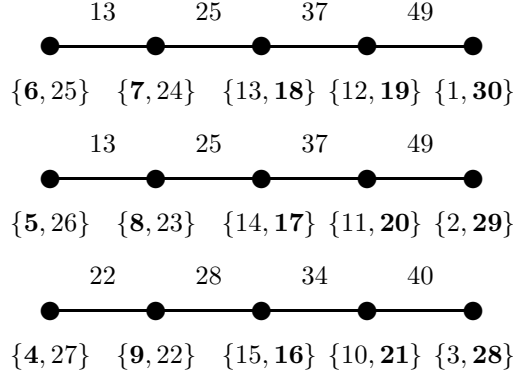
\vspace{0.1in}

Therefore in each of Cases 1, 2 and 3,
$K_{\rho} = \{ P_{\rho,i} : 1\le i \le m \}$ is a $(a_1,a_2,\ldots, \allowbreak a_{n_0})$-admissible
path partition of $G(a_1,a_2,\ldots,a_{n_0})$.
Since $(a_1,a_2,\ldots,a_{n_0})$ is not a palindrome, by Lemma \ref{lemmaSpanningGraphOfPathsCriteria}, for each $\rho$,
there exist $(m-1)!$ Klein bottle nonequivalent $C_4$-face-magic Klein bottle labelings $X$
of $\mathcal{K}_{m,n}$ such that $\mathcal{L}(X)=K_{\rho}$.
Hence, there exist at least $2^m (m-1)!$ distinct Klein bottle nonequivalent $C_4$-face-magic Klein bottle labelings $X$
of $\mathcal{K}_{m,n}$ such that $\mathcal{L}(X)$ is an $(a_1,a_2,\ldots,a_{n_0})$-admissible path
partition of $G(a_1,a_2,\ldots,a_{n_0})$.
\end{proof}

\begin{proposition} \label{propPalindromicSequences}
Let $m\ge 3$ be an odd integer.
Let $n\ge 6$ be an even integer such that $n_1=n/2$ is even.
Let $n_0=n_1 -1$, $\stkout{n}_1=\lfloor n_1/2\rfloor$, and $S=2(mn+1)$.
Let $(a_1,a_2,\ldots,a_{n_0})$ be one of the following four sequences:
\begin{enumerate}
\item $a_{2k-1} = mn$, for all $1 \le k \le \stkout{n}_1$, and $a_{2k} = mn+2$, for all $1 \le k \le \stkout{n}_1 -1$,
\item $a_{2k-1} = mn-1$, for all $1 \le k \le \stkout{n}_1$, and $a_{2k} = mn+3$, for all $1 \le k \le \stkout{n}_1 -1$,
\item $a_{2k-1} = m(n-1)+1$, for all $1 \le k \le \stkout{n}_1$, and $a_{2k} = m(n+1)=1$, for all $1 \le k \le \stkout{n}_1 -1$, or
\item $a_{2k-1} = m(n-2)+1$, for all $1 \le k \le \stkout{n}_1$, and $a_{2k} = m(n+2)+1$, for all $1 \le k \le \stkout{n}_1 -1$.
\end{enumerate}
Then, $(a_1,a_2,\ldots,a_{n_0})$ is a palindrome and there exists an $(a_1,a_2,\ldots,a_{n_0})$-admissible path
partition $K$ of $G(a_1,a_2,\ldots,a_{n_0})$. Furthermore, there exist $2^{m-1} (m-1)!$ distinct
Klein bottle nonequivalent $C_4$-face-magic Klein bottle labelings $X$ on $\mathcal{K}_{m,n}$ such that $\mathcal{L}(X)=K$.
\end{proposition}

\begin{proof}
We consider each sequence $(a_1,a_2,\ldots,a_{n_0})$ individually. The proof of the proposition for each sequence is similar.

\textbf{Case 1.} Suppose $a_{2k-1} = mn$, for all $1\le k \le \stkout{n}_1$.
and $a_{2k} = mn +2$, for all $1\le k \le \stkout{n}_1^{-}$.
We define the labels of $(a_1,a_2,\ldots,a_{n_0})$-admissible paths in $G(a_1,a_2,\ldots,a_{n_0})$ as follows:
for all $1 \le i \le m$, we let
\begin{align*}
x_{i,2j-1}                   &= n_1 (i-1) +2j -1,         &\mbox{ for }  1 \le j \le \stkout{n}_1,  \\
\tfrac{1}{2} S - x_{i,2j-1}  &=  n_1 (2m -i +1) - 2j + 2,  &\mbox{ for }  1 \le j \le \stkout{n}_1,   \\
x_{i,2j}                     &= n_1 (2m -i +1) - 2j + 1,           &\mbox{ for }  1 \le j \le \stkout{n}_1,  \mbox{ and }       \\
\tfrac{1}{2} S - x_{i,2j}    &=  n_1 (i-1) + 2j,             &\mbox{ for }  1 \le j \le \stkout{n}_1.  \\
\end{align*}
Let $P_{i}$ be the path on the sequence of vertices $( \{ x_{i,j}, \tfrac{1}{2}S - x_{i,j} \} : 1 \le j \le n_1 )$
with vertex label sequence $( x_{i,j} : 1 \le j \le n_1 )$.
The set of labels
\begin{align*}
\{ x_{i,2j-1} : 1 \le j \le \stkout{n}_1 \} \cup \{ \tfrac{1}{2} S - x_{i,2j} : 1 \le j \le \stkout{n}_1 \}
\end{align*}
is the set $\{ n_1 (i-1) +k: 1\le k \le n_1 \}$.
Hence for $1\le i \le m$ and $1 \le j \le n_1$, the vertices $\{ x_{i,j}, \tfrac{1}{2} S - x_{i,j} \}$
are distinct in the graph $G(a_1,a_2,\ldots,a_{n_0})$.

Furthermore, for $1\le i \le m$, we have
\begin{align*}
x_{i,2j-1} + x_{i,2j} &= (  n_1 (i-1) +2j -1 ) + ( n_1 ( 2m -i +1) -2j +1 )  \\
                          &= mn =a_{2j-1}, \mbox{ \ \ \ \ for }  1 \le j \le \stkout{n}_1, \mbox{ and } \\
x_{i,2j} + x_{i,2j+1} &= x_{i,2j} + x_{i,2j-1} +2   \\
             &=  mn + 2 =a_{2j}, \mbox{ \ \ \ \ for }  1 \le j \le \stkout{n}_1.
\end{align*}
Thus, $P_{i}$ is an $(a_1,a_2,\ldots,a_{n_0})$-admissible path in $G(a_1,a_2,\ldots,a_{n_0})$.
Also, we observe that $a_{2k-1} = a_{n_1 -2k +1} = mn$ for all $1\le k \le \stkout{n}_1$, and
$a_{2k} = a_{n_1 -2k } = mn +2$ for all $1\le k \le \stkout{n}_1^{-}$.
Hence, $(a_1,a_2,\ldots,a_{n_0})$ is a palindrome.

We consider the example where $m=3$, $n_1=6$ and $n=12$. Then, Figure \ref{figGraphProp34Case1}
provides the $(36,38,36,38,36)$-admissible path partition $K$ in the graph $G(36,38,36,38,36)$.

\begin{figure}
\hspace{0.0in}\begin{picture}(300,150)(-50,0)
\put(-10,120){\circle*{7}}
\put(-25,100){$\{ \textbf{1},36\}$}
\put(-10,120){\line(1,0){40}}
\put(5,130){36}

\put(30,120){\circle*{7}}
\put(15,100){$\{ 2,\textbf{35}\}$}
\put(30,120){\line(1,0){40}}
\put(45,130){38}

\put(70,120){\circle*{7}}
\put(55,100){$\{ \textbf{3},34\}$}
\put(70,120){\line(1,0){40}}
\put(85,130){36}

\put(110,120){\circle*{7}}
\put(95,100){$\{ 4,\textbf{33}\}$}
\put(110,120){\line(1,0){40}}
\put(125,130){38}

\put(150,120){\circle*{7}}
\put(135,100){$\{ \textbf{5},32\}$}
\put(150,120){\line(1,0){40}}
\put(165,130){36}

\put(190,120){\circle*{7}}
\put(175,100){$\{ 6,\textbf{31}\}$}

\put(-10,70){\circle*{7}}
\put(-25,50){$\{ \textbf{7},30\}$}
\put(-10,70){\line(1,0){40}}
\put(5,80){36}

\put(30,70){\circle*{7}}
\put(15,50){$\{ 8,\textbf{29}\}$}
\put(30,70){\line(1,0){40}}
\put(45,80){38}

\put(70,70){\circle*{7}}
\put(55,50){$\{ \textbf{9},28\}$}
\put(70,70){\line(1,0){40}}
\put(85,80){36}

\put(110,70){\circle*{7}}
\put(95,50){$\{ 10,\textbf{27}\}$}
\put(110,70){\line(1,0){40}}
\put(125,80){38}

\put(150,70){\circle*{7}}
\put(135,50){$\{ \textbf{11},26\}$}
\put(150,70){\line(1,0){40}}
\put(165,80){36}

\put(190,70){\circle*{7}}
\put(175,50){$\{ 12,\textbf{25}\}$}

\put(-10,20){\circle*{7}}
\put(-25,0){$\{ \textbf{13},24\}$}
\put(-10,20){\line(1,0){40}}
\put(5,30){36}

\put(30,20){\circle*{7}}
\put(15,0){$\{ 14,\textbf{23}\}$}
\put(30,20){\line(1,0){40}}
\put(45,30){38}

\put(70,20){\circle*{7}}
\put(55,0){$\{ \textbf{15},22\}$}
\put(70,20){\line(1,0){40}}
\put(85,30){36}

\put(110,20){\circle*{7}}
\put(95,0){$\{ 16,\textbf{21}\}$}
\put(110,20){\line(1,0){40}}
\put(125,30){38}

\put(150,20){\circle*{7}}
\put(135,0){$\{ \textbf{17},20\}$}
\put(150,20){\line(1,0){40}}
\put(165,30){36}

\put(190,20){\circle*{7}}
\put(175,0){$\{ 18,\textbf{19}\}$}
\end{picture}
\caption{$(36,38,36,38,36)$-admissible path partition $K$ in $G(36,38,36,38,36)$.}
\label{figGraphProp34Case1}
\end{figure}
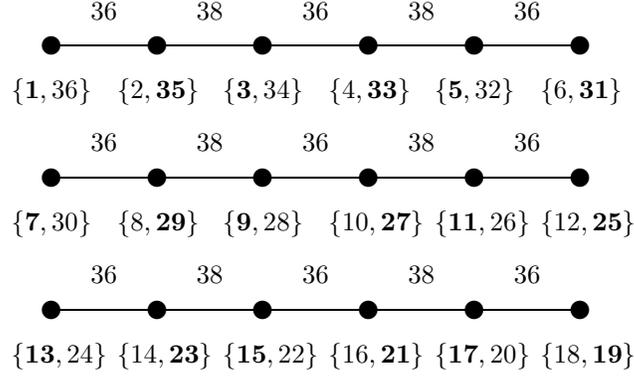

\textbf{Case 2.} Suppose $a_{2k-1} = mn -1$ for all $1\le k \le \stkout{n}_1$,
and $a_{2k} = mn +3$ for all $1\le k \le \stkout{n}_1^{-}$.
Since $m\ge 3$ is odd, we let $m = 2 m_1 +1$ for some positive integer $m_1$.
We define the labels of $(a_1,a_2,\ldots,a_{n_0})$-admissible paths in $G(a_1,a_2,\ldots,a_{n_0})$ as follows:
for all $1 \le i \le m_1$, we let
\begin{align*}
x_{2i-1,2j-1}                   &= n (i-1) +4j -3,         &\mbox{ for }  1 \le j \le \stkout{n}_1,  \\
\tfrac{1}{2} S - x_{2i-1,2j-1}  &=  n (m -i +1) - 4j + 4,  &\mbox{ for }  1 \le j \le \stkout{n}_1,   \\
x_{2i-1,2j}                     &= n (m -i +1) - 4j + 2,           &\mbox{ for }  1 \le j \le \stkout{n}_1,       \\
\tfrac{1}{2} S - x_{2i-1,2j}    &=  n (i-1) + 4j -1,             &\mbox{ for }  1 \le j \le \stkout{n}_1,  \\
x_{2i,2j-1}                   &= n (i-1) +4j -2,         &\mbox{ for }  1 \le j \le \stkout{n}_1,  \\
\tfrac{1}{2} S - x_{2i,2j-1}  &=  n (m -i +1) - 4j + 3,  &\mbox{ for }  1 \le j \le \stkout{n}_1,   \\
x_{2i,2j}                     &= n (m -i +1) - 4j + 1,           &\mbox{ for }  1 \le j \le \stkout{n}_1,  \mbox{ and }       \\
\tfrac{1}{2} S - x_{2i,2j}    &=  n (i-1) + 4j,             &\mbox{ for }  1 \le j \le \stkout{n}_1.  \\
\end{align*}
Also, we let
\begin{align*}
x_{m,2j-1}                   &= m_1 n +4j -3,         &\mbox{ for }  1 \le j \le \stkout{n}_1,  \\
\tfrac{1}{2} S - x_{m,2j-1}  &=  (m_1 +1) n - 4j + 4,  &\mbox{ for }  1 \le j \le \stkout{n}_1,   \\
x_{m,2j}                     &=  (m_1 +1) n  - 4j + 2,           &\mbox{ for }  1 \le j \le \stkout{n}_1,    \mbox{ and }     \\
\tfrac{1}{2} S - x_{m,2j}    &=   m_1 n + 4j -1,             &\mbox{ for }  1 \le j \le \stkout{n}_1.
\end{align*}
Let $P_{i}$ be the path on the sequence of vertices $( \{ x_{i,j}, \tfrac{1}{2}S - x_{i,j} \} : 1 \le j \le n_1 )$
with vertex label sequence $( x_{i,j} : 1 \le j \le n_1 )$.
The set of labels
\begin{align*}
\{ x_{2i-1,2j-1} : 1 \le j \le \stkout{n}_1 \} &\cup \{ \tfrac{1}{2} S - x_{2i-1,2j} : 1 \le j \le \stkout{n}_1 \} \\
     &\cup \{ \tfrac{1}{2} S - x_{2i,2j-1} : 1 \le j \le \stkout{n}_1 \} \cup \{ x_{2i,2j} : 1 \le j \le \stkout{n}_1 \}
\end{align*}
is the set $\{ n (i-1) +k: 1\le k \le n \}$.
Let $n_2 = \stkout{n}_1$.
Also, the set of labels
\begin{align*}
\{ x_{m,2j-1} : 1 \le j \le \stkout{n}_2^{+} \} &\cup \{ \tfrac{1}{2} S - x_{m,2j} : 1 \le j \le \stkout{n}_2 \}  \\
   &\cup \{ \tfrac{1}{2} S - x_{m,2j-1} : \stkout{n}_2^{+} < j \le n_2 \} \cup \{ x_{m,2j} :  \stkout{n}_2 \le j \le n_2 \}
\end{align*}
is the set $\{ m_1 n +k: 1\le k \le n_1 \}$.
Hence for $1\le i \le m$ and $1 \le j \le n_1$, the vertices $\{ x_{i,j}, \tfrac{1}{2} S - x_{i,j} \}$
are distinct in the graph $G(a_1,a_2,\ldots,a_{n_0})$.

Furthermore for $1\le i \le m_1$, we have
\begin{align*}
x_{2i-1,2j-1} + x_{2i-1,2j} &= (  n (i-1) +4j -3 ) + ( n ( m -i +1) -4j +2 )  \\
                          &= mn -1 =a_{2j-1}, \mbox{ \ \ \ \ for }  1 \le j \le \stkout{n}_1, \\
x_{2i-1,2j} + x_{2i-1,2j+1} &= x_{2i-1,2j} + x_{2i-1,2j-1} +4   \\
             &=  mn + 3 =a_{2j}, \mbox{ \ \ \ \ for }  1 \le j \le \stkout{n}_1, \\
x_{2i,2j-1} + x_{2i,2j} &= (  n (i-1) +4j -2 ) + ( n ( m -i +1) -4j +1 )  \\
                          &= mn -1 =a_{2j-1}, \mbox{ \ \ \ \ for }  1 \le j \le \stkout{n}_1, \mbox{ and } \\
x_{2i,2j} + x_{2i,2j+1} &= x_{2i,2j} + x_{2i,2j-1} +4   \\
             &=  mn + 3 =a_{2j}, \mbox{ \ \ \ \ for }  1 \le j \le \stkout{n}_1.
\end{align*}
In addition, we have
\begin{align*}
x_{m,2j-1} + x_{m,2j} &= (  m_1 n +4j -3 ) + ( (m_1 +1) n -4j +2 ) & \\
                          &= mn -1 =a_{2j-1}, &\mbox{ for }  1 \le j \le \stkout{n}_1, \mbox{ and } \\
x_{m,2j} + x_{m,2j+1} &= x_{m,2j} + x_{m,2j-1} +4  & \\
             &=  mn + 3 =a_{2j}, &\mbox{ for }  1 \le j \le \stkout{n}_1.
\end{align*}
Thus, $P_{i}$ is an $(a_1,a_2,\ldots,a_{n_0})$-admissible path in $G(a_1,a_2,\ldots,a_{n_0})$.
Also, we observe that $a_{2k-1} = a_{n_1 -2k +1} = mn$ for all $1\le k \le \stkout{n}_1$, and
$a_{2k} = a_{n_1 -2k } = mn +2$ for all $1\le k \le \stkout{n}_1^{-}$.
Hence, $(a_1,a_2,\ldots,a_{n_0})$ is a palindrome.

We consider the example where $m=3$, $n_1=6$ and $n=12$. Then, Figure \ref{figGraphProp34Case2}
provides the $(35,39,35,39,35)$-admissible path partition $K$ in the graph $G(35,39,35,39,35)$.

\begin{figure}
\hspace{0.0in}\begin{picture}(300,150)(-50,0)
\put(-10,120){\circle*{7}}
\put(-25,100){$\{ \textbf{1},36\}$}
\put(-10,120){\line(1,0){40}}
\put(5,130){35}

\put(30,120){\circle*{7}}
\put(15,100){$\{ 3,\textbf{34}\}$}
\put(30,120){\line(1,0){40}}
\put(45,130){39}

\put(70,120){\circle*{7}}
\put(55,100){$\{ \textbf{5},32\}$}
\put(70,120){\line(1,0){40}}
\put(85,130){35}

\put(110,120){\circle*{7}}
\put(95,100){$\{ 7,\textbf{30}\}$}
\put(110,120){\line(1,0){40}}
\put(125,130){39}

\put(150,120){\circle*{7}}
\put(135,100){$\{ \textbf{9},28\}$}
\put(150,120){\line(1,0){40}}
\put(165,130){35}

\put(190,120){\circle*{7}}
\put(175,100){$\{ 11,\textbf{26}\}$}

\put(-10,70){\circle*{7}}
\put(-25,50){$\{ \textbf{2},35\}$}
\put(-10,70){\line(1,0){40}}
\put(5,80){35}

\put(30,70){\circle*{7}}
\put(15,50){$\{ 4,\textbf{33}\}$}
\put(30,70){\line(1,0){40}}
\put(45,80){39}

\put(70,70){\circle*{7}}
\put(55,50){$\{ \textbf{6},31\}$}
\put(70,70){\line(1,0){40}}
\put(85,80){35}

\put(110,70){\circle*{7}}
\put(95,50){$\{ 8,\textbf{29}\}$}
\put(110,70){\line(1,0){40}}
\put(125,80){39}

\put(150,70){\circle*{7}}
\put(135,50){$\{ \textbf{10},27\}$}
\put(150,70){\line(1,0){40}}
\put(165,80){35}

\put(190,70){\circle*{7}}
\put(175,50){$\{ 12,\textbf{25}\}$}

\put(-10,20){\circle*{7}}
\put(-25,0){$\{ \textbf{13},24\}$}
\put(-10,20){\line(1,0){40}}
\put(5,30){35}

\put(30,20){\circle*{7}}
\put(15,0){$\{ 15,\textbf{22}\}$}
\put(30,20){\line(1,0){40}}
\put(45,30){39}

\put(70,20){\circle*{7}}
\put(55,0){$\{ \textbf{17},20\}$}
\put(70,20){\line(1,0){40}}
\put(85,30){35}

\put(110,20){\circle*{7}}
\put(95,0){$\{ \textbf{18},19\}$}
\put(110,20){\line(1,0){40}}
\put(125,30){39}

\put(150,20){\circle*{7}}
\put(135,0){$\{ 16,\textbf{21}\}$}
\put(150,20){\line(1,0){40}}
\put(165,30){35}

\put(190,20){\circle*{7}}
\put(175,0){$\{ \textbf{14},23\}$}
\end{picture}
\caption{$(35,39,35,39,35)$-admissible path partition $K$ in $G(35,39,35,39,35)$.}
\label{figGraphProp34Case2}
\end{figure}
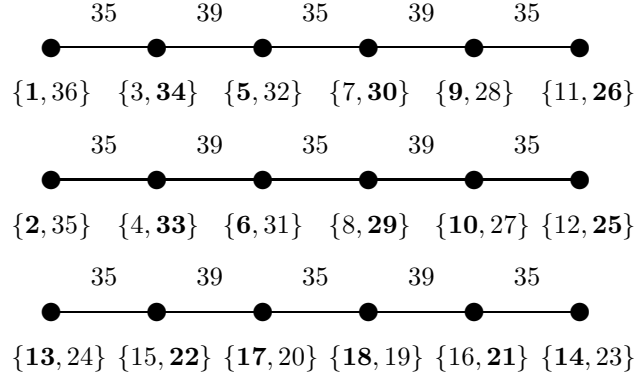

\textbf{Case 3.} Suppose $a_{2k-1} = m(n-1) +1$ for all $1\le k \le \stkout{n}_1$,
and $a_{2k} = m(n+1) +1$ for all $1\le k \le \stkout{n}_1^{-}$.
We define the labels of $(a_1,a_2,\ldots,a_{n_0})$-admissible paths in $G(a_1,a_2,\ldots,a_{n_0})$ as follows:
for all $1 \le i \le m$, we let
\begin{align*}
x_{i,2j-1}                   &= m (2j-2) +i,         &\mbox{ for }  1 \le j \le \stkout{n}_1,  \\
\tfrac{1}{2} S - x_{i,2j-1}  &= m (n -2j +2) - i + 1,  &\mbox{ for }  1 \le j \le \stkout{n}_1,   \\
x_{i,2j}                     &= m (n -2j +1) - i + 1,           &\mbox{ for }  1 \le j \le \stkout{n}_1,  \mbox{ and }       \\
\tfrac{1}{2} S - x_{i,2j}    &=  m (2j-1) + i,             &\mbox{ for }  1 \le j \le \stkout{n}_1.  \\
\end{align*}
Let $P_{i}$ be the path on the sequence of vertices $( \{ x_{i,j}, \tfrac{1}{2}S - x_{i,j} \} : 1 \le j \le n_1 )$
with vertex label sequence $( x_{i,j} : 1 \le j \le n_1 )$.
The set of labels
\begin{align*}
\{ x_{i,2j-1} : 1 \le j \le \stkout{n}_1 \} \cup \{ \tfrac{1}{2} S - x_{i,2j} : 1 \le j \le \stkout{n}_1 \}
\end{align*}
is the set $\{ mk + i: 0\le k \le n_1 -1 \}$.
Hence for $1\le i \le m$ and $1 \le j \le n_1$, the vertices $\{ x_{i,j}, \tfrac{1}{2} S - x_{i,j} \}$
are distinct in the graph $G(a_1,a_2,\ldots,a_{n_0})$.

Furthermore for $1\le i \le m$, we have
\begin{align*}
x_{i,2j-1} + x_{i,2j} &= (  m (2j-2) +i ) + ( m (n -2j +1) - i + 1 ) & \\
                          &= m(n-1) +1 =a_{2j-1}, &\mbox{ for }  1 \le j \le \stkout{n}_1, \mbox{ and } \\
x_{i,2j} + x_{i,2j+1} &= x_{i,2j} + x_{i,2j-1} +2m  & \\
             &=  m(n+1) + 1 =a_{2j}, &\mbox{ for }  1 \le j \le \stkout{n}_1.
\end{align*}
Thus, $P_{i}$ is an $(a_1,a_2,\ldots,a_{n_0})$-admissible path in $G(a_1,a_2,\ldots,a_{n_0})$.
Also, we observe that $a_{2k-1} = a_{n_1 -2k +1} = m(n-1) +1$, for all $1\le k \le \stkout{n}_1$, and
$a_{2k} = a_{n_1 -2k } = m(n+1) +1$, for all $1\le k \le \stkout{n}_1^{-}$.
Hence, $(a_1,a_2,\ldots,a_{n_0})$ is a palindrome.

We consider the example where $m=3$, $n_1=6$ and $n=12$. Then, Figure \ref{figGraphProp34Case3}
provides the $(34,40,34,40,34)$-admissible path partition $K$ in the graph $G(34,40,34,40,34)$.

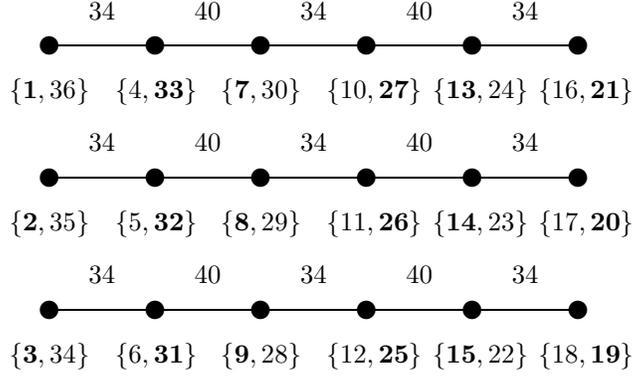
\begin{figure}
\hspace{0.0in}\begin{picture}(300,150)(-50,0)
\put(-10,120){\circle*{7}}
\put(-25,100){$\{ \textbf{1},36\}$}
\put(-10,120){\line(1,0){40}}
\put(5,130){34}

\put(30,120){\circle*{7}}
\put(15,100){$\{ 4,\textbf{33}\}$}
\put(30,120){\line(1,0){40}}
\put(45,130){40}

\put(70,120){\circle*{7}}
\put(55,100){$\{ \textbf{7},30\}$}
\put(70,120){\line(1,0){40}}
\put(85,130){34}

\put(110,120){\circle*{7}}
\put(95,100){$\{ 10,\textbf{27}\}$}
\put(110,120){\line(1,0){40}}
\put(125,130){40}

\put(150,120){\circle*{7}}
\put(135,100){$\{ \textbf{13},24\}$}
\put(150,120){\line(1,0){40}}
\put(165,130){34}

\put(190,120){\circle*{7}}
\put(175,100){$\{ 16,\textbf{21}\}$}

\put(-10,70){\circle*{7}}
\put(-25,50){$\{ \textbf{2},35\}$}
\put(-10,70){\line(1,0){40}}
\put(5,80){34}

\put(30,70){\circle*{7}}
\put(15,50){$\{ 5,\textbf{32}\}$}
\put(30,70){\line(1,0){40}}
\put(45,80){40}

\put(70,70){\circle*{7}}
\put(55,50){$\{ \textbf{8},29\}$}
\put(70,70){\line(1,0){40}}
\put(85,80){34}

\put(110,70){\circle*{7}}
\put(95,50){$\{ 11,\textbf{26}\}$}
\put(110,70){\line(1,0){40}}
\put(125,80){40}

\put(150,70){\circle*{7}}
\put(135,50){$\{ \textbf{14},23\}$}
\put(150,70){\line(1,0){40}}
\put(165,80){34}

\put(190,70){\circle*{7}}
\put(175,50){$\{ 17,\textbf{20}\}$}

\put(-10,20){\circle*{7}}
\put(-25,0){$\{ \textbf{3},34\}$}
\put(-10,20){\line(1,0){40}}
\put(5,30){34}

\put(30,20){\circle*{7}}
\put(15,0){$\{ 6,\textbf{31}\}$}
\put(30,20){\line(1,0){40}}
\put(45,30){40}

\put(70,20){\circle*{7}}
\put(55,0){$\{ \textbf{9},28\}$}
\put(70,20){\line(1,0){40}}
\put(85,30){34}

\put(110,20){\circle*{7}}
\put(95,0){$\{ 12,\textbf{25}\}$}
\put(110,20){\line(1,0){40}}
\put(125,30){40}

\put(150,20){\circle*{7}}
\put(135,0){$\{ \textbf{15},22\}$}
\put(150,20){\line(1,0){40}}
\put(165,30){34}

\put(190,20){\circle*{7}}
\put(175,0){$\{ 18,\textbf{19}\}$}
\end{picture}
\caption{$(34,40,34,40,34)$-admissible path partition $K$ in $G(34,40,34,40,34)$.}
\label{figGraphProp34Case3}
\end{figure}

\textbf{Case 4.} Suppose $a_{2k-1} = m(n-2) +1$ for all $1\le k \le \stkout{n}_1$,
and $a_{2k} = m(n+2) +1$ for all $1\le k \le \stkout{n}_1^{-}$.
We define the labels of $(a_1,a_2,\ldots,a_{n_0})$-admissible paths in $G(a_1,a_2,\ldots,a_{n_0})$ as follows:
for all $1 \le i \le m$, we let
\begin{align*}
x_{i,2j-1}                   &= m (4j-4) +i,         &\mbox{ for }  1 \le j \le \stkout{n}_1,  \\
\tfrac{1}{2} S - x_{i,2j-1}  &= m (n -4j +4) - i + 1,  &\mbox{ for }  1 \le j \le \stkout{n}_1,   \\
x_{i,2j}                     &= m (n -4j +2) - i + 1,           &\mbox{ for }  1 \le j \le \stkout{n}_1,  \mbox{ and }       \\
\tfrac{1}{2} S - x_{i,2j}    &=  m (4j-2) + i,             &\mbox{ for }  1 \le j \le \stkout{n}_1.  \\
\end{align*}
Let $P_{i}$ be the path on the sequence of vertices $( \{ x_{i,j}, \tfrac{1}{2}S - x_{i,j} \} : 1 \le j \le n_1 )$
with vertex label sequence $( x_{i,j} : 1 \le j \le n_1 )$.
The set of labels
\begin{align*}
\{ x_{i,2j-1} : 1 \le j \le \stkout{n}_1 \} \cup \{ \tfrac{1}{2} S - x_{i,2j} : 1 \le j \le \stkout{n}_1 \}
\end{align*}
is the set $\{ 2mk + i: 0\le k \le n_1 -1 \}$.
Also, the set of labels
\begin{align*}
\{ \tfrac{1}{2} S - x_{i,2j-1} : 1 \le j \le \stkout{n}_1^{+} \} \cup \{  x_{i,2j} : 1 \le j \le \stkout{n}_1 \}
\end{align*}
is the set $\{ 2mk-i+1: 1\le k \le n_1 \}$.
Hence for $1\le i \le m$ and $1 \le j \le n_1$, the vertices $\{ x_{i,j}, \tfrac{1}{2} S - x_{i,j} \}$
are distinct in the graph $G(a_1,a_2,\ldots,a_{n_0})$.

Furthermore for $1\le i \le m$, we have
\begin{align*}
x_{i,2j-1} + x_{i,2j} &= (  m (4j-4) +i ) + ( m (n -4j +2) - i + 1 ) & \\
                          &= m(n-2) +1 =a_{2j-1}, &\mbox{ for }  1 \le j \le \stkout{n}_1, \mbox{ and } \\
x_{i,2j} + x_{i,2j+1} &= x_{i,2j} + x_{i,2j-1} +4m  & \\
             &=  m(n+2) + 1 =a_{2j}, &\mbox{ for }  1 \le j \le \stkout{n}_1.
\end{align*}
Thus, $P_{i}$ is an $(a_1,a_2,\ldots,a_{n_0})$-admissible path in $G(a_1,a_2,\ldots,a_{n_0})$.
Also, we observe that $a_{2k-1} = a_{n_1 -2k +1} = m(n-2) +1$ for all $1\le k \le \stkout{n}_1$, and
$a_{2k} = a_{n_1 -2k } = m(n+2) +1$ for all $1\le k \le \stkout{n}_1^{-}$.
Hence, $(a_1,a_2,\ldots,a_{n_0})$ is a palindrome.

We consider the example where $m=3$, $n_1=6$ and $n=12$. Then, Figure \ref{figGraphProp34Case4}
provides the $(31,43,31,43,31)$-admissible path partition $K$ in the graph $G(31,43,31,43,31)$.

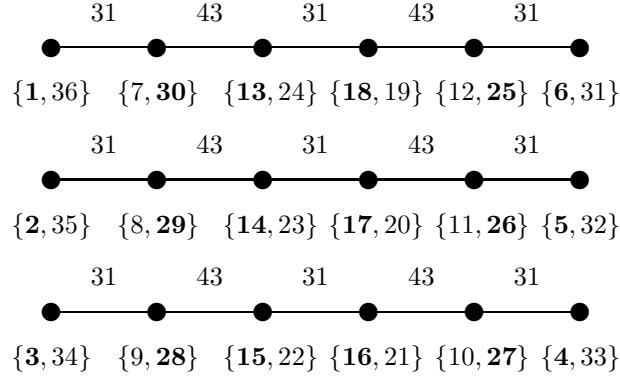
\begin{figure}
\hspace{0.0in}\begin{picture}(300,150)(-50,0)
\put(-10,120){\circle*{7}}
\put(-25,100){$\{ \textbf{1},36\}$}
\put(-10,120){\line(1,0){40}}
\put(5,130){31}

\put(30,120){\circle*{7}}
\put(15,100){$\{ 7,\textbf{30}\}$}
\put(30,120){\line(1,0){40}}
\put(45,130){43}

\put(70,120){\circle*{7}}
\put(55,100){$\{ \textbf{13},24\}$}
\put(70,120){\line(1,0){40}}
\put(85,130){31}

\put(110,120){\circle*{7}}
\put(95,100){$\{ \textbf{18},19\}$}
\put(110,120){\line(1,0){40}}
\put(125,130){43}

\put(150,120){\circle*{7}}
\put(135,100){$\{ 12,\textbf{25}\}$}
\put(150,120){\line(1,0){40}}
\put(165,130){31}

\put(190,120){\circle*{7}}
\put(175,100){$\{ \textbf{6},31\}$}

\put(-10,70){\circle*{7}}
\put(-25,50){$\{ \textbf{2},35\}$}
\put(-10,70){\line(1,0){40}}
\put(5,80){31}

\put(30,70){\circle*{7}}
\put(15,50){$\{ 8,\textbf{29}\}$}
\put(30,70){\line(1,0){40}}
\put(45,80){43}

\put(70,70){\circle*{7}}
\put(55,50){$\{ \textbf{14},23\}$}
\put(70,70){\line(1,0){40}}
\put(85,80){31}

\put(110,70){\circle*{7}}
\put(95,50){$\{ \textbf{17},20\}$}
\put(110,70){\line(1,0){40}}
\put(125,80){43}

\put(150,70){\circle*{7}}
\put(135,50){$\{ 11,\textbf{26}\}$}
\put(150,70){\line(1,0){40}}
\put(165,80){31}

\put(190,70){\circle*{7}}
\put(175,50){$\{ \textbf{5},32\}$}

\put(-10,20){\circle*{7}}
\put(-25,0){$\{ \textbf{3},34\}$}
\put(-10,20){\line(1,0){40}}
\put(5,30){31}

\put(30,20){\circle*{7}}
\put(15,0){$\{ 9,\textbf{28}\}$}
\put(30,20){\line(1,0){40}}
\put(45,30){43}

\put(70,20){\circle*{7}}
\put(55,0){$\{ \textbf{15},22\}$}
\put(70,20){\line(1,0){40}}
\put(85,30){31}

\put(110,20){\circle*{7}}
\put(95,0){$\{ \textbf{16},21\}$}
\put(110,20){\line(1,0){40}}
\put(125,30){43}

\put(150,20){\circle*{7}}
\put(135,0){$\{ 10,\textbf{27}\}$}
\put(150,20){\line(1,0){40}}
\put(165,30){31}

\put(190,20){\circle*{7}}
\put(175,0){$\{ \textbf{4},33\}$}
\end{picture}
\caption{$(31,43,31,43,31)$-admissible path partition $K$ in $G(31,43,31,43,31)$.}
\label{figGraphProp34Case4}
\end{figure}
\vspace{0.1in}

Therefore in each of Cases 1, 2 3, and 4,
$K = \{ P_{i} : 1\le i \le m \}$ is a $(a_1,a_2,\ldots, \allowbreak a_{n_0})$-admissible
path partition of $G(a_1,a_2,\ldots,a_{n_0})$.
Since $(a_1,a_2,\ldots,a_{n_0})$ is a palindrome, by Lemma \ref{lemmaSpanningGraphOfPathsCriteria},
there exist $2^{m-1} (m-1)!$ distinct Klein bottle nonequivalent $C_4$-face-magic Klein bottle labelings $X$
of $\mathcal{K}_{m,n}$ such that $\mathcal{L}(X)=K$.
Hence, there exist at least $2^{m-1} (m-1)!$ distinct Klein bottle nonequivalent $C_4$-face-magic Klein bottle labelings $X$
of $\mathcal{K}_{m,n}$ such that $\mathcal{L}(X)$ is an $(a_1,a_2,\ldots,a_{n_0})$-admissible path
partition of $G(a_1,a_2,\ldots,a_{n_0})$.
\end{proof}

\begin{theorem}
Let $m\ge 3$ be an odd integer and let $n \ge 6$ be an even integer.
\begin{enumerate}
\item If $n \equiv 2\pmod{4}$, then the minimum number of distinct Klein bottle nonequivalent $C_4$-face-magic Klein bottle
labelings $X$ on $\mathcal{K}_{m,n}$ is $(6\cdot 2^m)(m-1)!$.
\item If $n \equiv 0\pmod{4}$, then the minimum number of distinct Klein bottle nonequivalent $C_4$-face-magic Klein bottle
labelings $X$ on $\mathcal{K}_{m,n}$ is $(5\cdot 2^m)(m-1)!$.
\end{enumerate}
\end{theorem}

\begin{proof}
Let $A_i = (a_{i,1},a_{i,2},\ldots,a_{i,n_0})$ be the sequence where $a_{1,k}= mn -n_0 +k$,
$a_{2,k}= m(n_1 + k) +1$, and $a_{3,k} = 2mk + 1$, for all $1\le k \le n_0$.
Also, let $A_i = (a_{i,1},a_{i,2},\ldots,a_{i,n_0})$ be the sequence where $a_{4,k}= (2m-1)n_1 + 2k + 1$,
$a_{5,k}= m n_1 + 2mk + 1$, and $a_{6,k} = 4mk + 1$, for all $1\le k \le n_0$.
Furthermore, let $\rho: NAPP(m,n)\longrightarrow NAPP(m,n)$ (see Definition \ref{defnAPPSet} for the definition of $NAPP(m,n)$)
and $\sigma: NAPP(m,n)\longrightarrow NAPP(m,n)$
be functions defined by
\begin{align*}
\rho(a_1,a_2,\ldots,a_{n_0 -1},a_{n_0}) &= (a_{n_0}, a_{n_0 -1},\ldots,a_2,a_1) \mbox{ \ \ \ \ and } \\
\sigma(a_1,a_2,\ldots,a_{n_0 -1},a_{n_0}) &= (S- a_1,S- a_2,\ldots,S- a_{n_0 -1},S- a_{n_0}).
\end{align*}
Suppose $A, A' \in NAPP(m,n)$ are sequences such that $A$ is path partition equivalent to $A'$.
By Lemma \ref{lemmaPathPartitionEquivalence}, either $A'=A$, $A'=\rho(A)$, $A'=\sigma(A)$, or $A'=(\sigma\rho)(A)$.
Of these four possibilities for $A'$, let $A'$ be the minimal element in the set
$\{ A, \rho(A), \sigma(A), (\sigma\rho)(A)\}$ under the lexicographic order.
Let
\begin{align*}
MNAPP(m,n) = &\biggl\{ A \in NAPP(m,n) : A \mbox{ is the minimal element in the set } \\
  &\{ A, \rho(A), \sigma(A), (\sigma\rho)(A) \} \mbox{ under the lexicographic order}\biggr\}
\end{align*}
be the minimal set of sequences $A\in NAPP(m,n)$ needed
in order to count all Klein bottle nonequivalent $C_4$-face-magic Klein bottle labelings on $\mathcal{K}_{m,n}$ in which $A$ is not a palindrome.
Let $\widetilde{A}_1$ and $\widetilde{A}_2$ be two distinct sequences in $MNAPP(m,n)$.
For $i=1,2$, let $X_i$ be a $C_4$-face-magic Klein bottle labeling on $\mathcal{K}_{m,n}$
such that $\mathcal{L}(X_i)$ is an $\widetilde{A}_i$-admissible path partition of $G(\widetilde{A}_i)$.
Then, $\widetilde{A}_2$ is not any of the pairs $\widetilde{A}_1$, $\rho(\widetilde{A}_1)$,
$\sigma(\widetilde{A}_1)$, or $(\sigma \rho)(\widetilde{A}_1)$.
By Lemma \ref{lemmaPathPartitionEquivalence}, $\widetilde{A}_1$ is not path partition equivalent to $\widetilde{A}_2$.
Thus, $X_1$ and $X_2$ are Klein bottle nonequivalent $C_4$-face-magic Klein bottle labelings on $\mathcal{K}_{m,n}$.
Let $N(A)$ be the number of distinct $A$-admissible path partitions of $G(A)$.

We consider the cases $n \equiv 0\pmod{4}$ and $n \equiv 2\pmod{4}$ individually.

\textbf{Case 1.} Suppose $n \equiv 2\pmod{4}$.
Then, $n_1$ is odd.
Let $A=(a_1,a_2,\ldots,a_{n_0})$ be a sequence such that there is an $A$-admissible path partition of $G(A)$.
For the purposes of contradiction, assume $A$ is a palindrome.
Since $x_{1,\stkout{n}_1} + x_{1,\stkout{n}_1 +1} = a_{\stkout{n}_1} = a_{\stkout{n}_1 +1} = x_{1,\stkout{n}_1 +1} + x_{1,\stkout{n}_1 +2}$,
we have $x_{1,\stkout{n}_1} = x_{1,\stkout{n}_1 +2}$.
This is a contradiction. Therefore, $A$ is not a palindrome.

By Lemma \ref{lemmaSpanningGraphOfPathsCriteria}, there are $N(A) \cdot (m-1)!$ distinct Klein bottle
nonequivalent $C_4$-face-magic Klein bottle labelings $X$ on $\mathcal{K}_{m,n}$ such that
$\mathcal{L}(X)$ is an $A$-admissible path partition of $G(A)$.
Thus, the number of distinct Klein bottle nonequivalent $C_4$-face-magic Klein bottle labelings on $\mathcal{K}_{m,n}$,
denoted by $NKBL(m,n)$,  is
\begin{equation*}
  NKBL(m,n) = \sum_{A\in MNAPP(m,n)} N(A) \cdot (m-1)!.
\end{equation*}

By Proposition \ref{propLowerSequences}, we have $N(A_i)\ge 2^m$ for all $1\le i \le 3$.
By Proposition \ref{propMiddleSequences}, we have $N(A_i)\ge 2^m$ for all $4\le i \le 6$.
Thus, the number of distinct Klein bottle nonequivalent $C_4$-face-magic Klein bottle labelings on $\mathcal{K}_{m,n}$ satisfies
\begin{align*}
  NKBL(m,n) &= \sum_{A\in MNAPP(m,n)} N(A) \cdot (m-1)!  \\
             &\ge \sum_{i=1}^{6} N(A_i) \cdot (m-1)! \ge (6\cdot 2^m)\cdot (m-1)!.
\end{align*}

\textbf{Case 2.} Suppose $n \equiv 0\pmod{4}$.
Then, $n_1$ is even.

Let $A_7=(a_{7,1},a_{7,2},\ldots,a_{7,n_0})$ be the sequence where
$a_{7,2k-1} = mn$ for all $1 \le k \le \stkout{n}_1$, and $a_{7,2k} = mn+2$ for all $1 \le k \le \stkout{n}_1 -1$.
Let $A_8=(a_{8,1},a_{8,2},\ldots,a_{8,n_0})$ be the sequence where
$a_{8,2k-1} = mn-1$ for all $1 \le k \le \stkout{n}_1$, and $a_{8,2k} = mn+3$ for all $1 \le k \le \stkout{n}_1 -1$.
Let $A_9=(a_{9,1},a_{9,2},\ldots,a_{9,n_0})$ be the sequence where
$a_{9,2k-1} = m(n-1)+1$ for all $1 \le k \le \stkout{n}_1$, and $a_{9,2k} = m(n+1)=1$ for all $1 \le k \le \stkout{n}_1 -1$.
Let $A_{10}=(a_{10,1},a_{10,2},\ldots,a_{10,n_0})$ be the sequence where
$a_{10,2k-1} = m(n-2)+1$ for all $1 \le k \le \stkout{n}_1$, and $a_{10,2k} = m(n+2)+1$ for all $1 \le k \le \stkout{n}_1 -1$.

Suppose $A, A' \in PAPP(m,n)$ (see Definition \ref{defnAPPSet} for the definition of $PAPP(m,n)$) are sequences such that $A$ is path partition equivalent to $A'$.
Define $\sigma: PAPP(m,n)\longrightarrow PAPP(m,n)$ by
\begin{align*}
\sigma(a_1,a_2,\ldots,a_{n_0 -1},a_{n_0}) &= (S- a_1,S- a_2,\ldots,S- a_{n_0 -1},S- a_{n_0}).
\end{align*}
By Lemma \ref{lemmaPathPartitionEquivalence}, either $A'=A$, or $A'=\sigma(A)$.
Of these two possibilities for $A'$, let $A'$ be the minimal element in the set $\{ A, \sigma(A)\}$ under the lexicographic order.
Let
\begin{align*}
MPAPP(m,n) = &\biggl\{ A \in PAPP(m,n) : A \mbox{ is the minimal element in the set } \\
  &\{ A, \sigma(A) \} \mbox{ under the lexicographic order}\biggr\}
\end{align*}
be the minimal set of sequences $A\in PAPP(m,n)$ needed
in order to count all Klein bottle nonequivalent $C_4$-face-magic Klein bottle labelings on $\mathcal{K}_{m,n}$ in which $A$ is a palindrome.
Let $\widetilde{A}_1$ and $\widetilde{A}_2$ be two distinct sequences in $MPAPP(m,n)$.
For $i=1,2$, let $X_i$ be a $C_4$-face-magic Klein bottle labeling on $\mathcal{K}_{m,n}$
such that $\mathcal{L}(X_i)$ is an $\widetilde{A}_i$-admissible path partition of $G(\widetilde{A}_i)$.
Then, $\widetilde{A}_2$ is not any of the pairs $\widetilde{A}_1$, or $\sigma(\widetilde{A}_1)$.
By Lemma \ref{lemmaPathPartitionEquivalence}, $\widetilde{A}_1$ is not path partition equivalent to $\widetilde{A}_2$.
Thus $X_1$ and $X_2$ are Klein bottle nonequivalent $C_4$-face-magic Klein bottle labelings on $\mathcal{K}_{m,n}$.
Let $N(A)$ be the number of distinct $A$-admissible path partitions of $G(A)$.
On one hand, if $A$ is not a palindrome, by Lemma \ref{lemmaSpanningGraphOfPathsCriteria}, there are $N(A) \cdot (m-1)!$ distinct Klein bottle
nonequivalent $C_4$-face-magic Klein bottle labelings $X$ on $\mathcal{K}_{m,n}$ such that
$\mathcal{L}(X)$ is an $A$-admissible path partition of $G(A)$.
On the other hand, if $A$ is a palindrome, by Lemma \ref{lemmaSpanningGraphOfPathsCriteria}, there are $N(A) \cdot 2^{m-1} (m-1)!$ distinct Klein bottle
nonequivalent $C_4$-face-magic Klein bottle labelings $X$ on $\mathcal{K}_{m,n}$ such that
$\mathcal{L}(X)$ is an $A$-admissible path partition of $G(A)$.
Thus, the number of distinct Klein bottle nonequivalent $C_4$-face-magic Klein bottle labelings on $\mathcal{K}_{m,n}$,
denoted by $NKBL(m,n)$,  is
\begin{equation*}
  NKBL(m,n) = \sum_{{\scriptscriptstyle A\in MNAPP(m,n)}} N(A) \cdot (m-1)! + \sum_{{\scriptscriptstyle A\in MPAPP(m,n)}} N(A) \cdot 2^{m-1} (m-1)!.
\end{equation*}

By Proposition \ref{propLowerSequences}, we have $N(A_i)\ge 2^m$ for all $1\le i \le 3$.
By Proposition \ref{propPalindromicSequences}, we have $N(A_i)\ge 1$ for all $7\le i \le 10$.
Thus, the number of distinct Klein bottle nonequivalent $C_4$-face-magic Klein bottle labelings on $\mathcal{K}_{m,n}$ satisfies
\begin{align*}
  NKBL(m,n) = \sum_{{\scriptscriptstyle A\in MNAPP(m,n)}} N(A) \cdot (m-1)!  + \sum_{{\scriptscriptstyle A\in MPAPP(m,n)}} N(A) \cdot 2^{m-1} (m-1)! \\
             \ge \sum_{i=1}^{3} N(A_i) \cdot (m-1)! + \sum_{i=7}^{10} N(A_i) \cdot 2^{m-1} (m-1)! \ge (5\cdot 2^m)\cdot (m-1)!.
\end{align*}
\end{proof}

\end{document}